\documentclass[11pt,a4paper]{amsproc}
\usepackage[T1]{fontenc}
\usepackage[utf8]{inputenc}
\usepackage{lmodern}
\usepackage[a4paper,margin=25mm,headheight=14pt,footskip=10mm]{geometry}
\usepackage{amsmath,amssymb,amsthm,mathtools,enumitem,booktabs,longtable}
\usepackage{microtype}
\usepackage{xcolor,fancyhdr}
\usepackage{tikz,tikz-cd}
\usetikzlibrary{decorations.pathreplacing,calc}

\newtheoremstyle{plain}{}{}{\itshape}{}{\bfseries}{.}{5pt plus 1pt minus 1pt}{}
\newtheoremstyle{definition}{}{}{\normalfont}{}{\bfseries}{.}{5pt plus 1pt minus 1pt}{}
\newtheoremstyle{remark}{}{}{\normalfont}{}{\itshape}{.}{5pt plus 1pt minus 1pt}{}
\makeatletter
\renewenvironment{proof}[1][\proofname]{\par
  \pushQED{\qed}%
  \normalfont \topsep6\p@\@plus6\p@\relax
  \trivlist
  \item[\hskip\labelsep\itshape #1\@addpunct{.}]\ignorespaces
}{%
  \popQED\endtrivlist\@endpefalse
}
\makeatother

\usepackage[backref=section]{hyperref}
\usepackage{bookmark}
\definecolor{citegreen}{rgb}{0,0.5,0.15}
\definecolor{urlgreen}{rgb}{0,0.35,0.15}
\hypersetup{colorlinks=true,
  linkcolor=blue,
  citecolor=citegreen,
  urlcolor=urlgreen,
  linkbordercolor={0 0.4 1},
  citebordercolor={0 0.9 0.37},
  unicode=true,
  pdftitle={Controlled analytic properties and the Baum--Connes conjecture},
  pdfauthor={Vadim Alekseev and Martin Finn-Sell},
  pdfsubject={}}

\makeatletter
\let\originalbackout\Hy@backout
\renewcommand*{\Hy@backout}[1]{%
  \begingroup
    \ifnum\value{subsection}>0
      \edef\@currentlabel{\thesubsection}%
      \edef\@currentHref{subsection.\theHsubsection}%
    \else\ifnum\value{section}>0
      \edef\@currentlabel{\thesection}%
      \edef\@currentHref{\Hy@chapapp.\theHsection}%
    \else
      \def\@currentlabel{0}%
      \def\@currentHref{Doc-Start}%
    \fi\fi
    \originalbackout{#1}%
  \endgroup
}
\makeatother

\renewcommand*{\backref}[1]{}
\renewcommand*{\backrefalt}[4]{%
  \ifnum#1>0
    \begingroup
      \hypersetup{linkcolor=blue}\color{blue}%
      \ensuremath{\uparrow}\,#2%
    \endgroup
  \fi
}

\numberwithin{equation}{section}
\theoremstyle{plain}
\newtheorem{theorem}{Theorem}[section]
\newtheorem{proposition}[theorem]{Proposition}
\newtheorem{lemma}[theorem]{Lemma}
\newtheorem*{restatedlemmainner}{\restatedlemmaheading}
\newenvironment{restatedlemma}[2][]{%
  \def\restatedlemmaheading{Lemma~\ref{#2}}%
  \begin{restatedlemmainner}[#1]%
}{\end{restatedlemmainner}}

\newtheorem{introtheorem}{Theorem}

\theoremstyle{definition}
\newtheorem{definition}[theorem]{Definition}
\newtheorem{construction}[theorem]{Construction}
\theoremstyle{remark}
\newtheorem{remark}[theorem]{Remark}
\newtheorem*{remark*}{Remark}
\newcommand{\partref}[2]{\ref{#1}~\textup{(#2)}}
\newcommand{\C}{\mathbb C}
\newcommand{\N}{\mathbb N}
\newcommand{\Z}{\mathbb Z}
\newcommand{\K}{\mathcal K}
\newcommand{\B}{\mathcal B}
\newcommand{\Cr}{C^*_r}
\newcommand{\prop}{\operatorname{Prop}}
\newcommand{\supp}{\operatorname{supp}}
\newcommand{\diam}{\operatorname{diam}}
\newcommand{\girth}{\operatorname{girth}}
\newcommand{\Tr}{\operatorname{Tr}}
\newcommand{\pos}{\operatorname{pos}}
\newcommand{\eps}{\varepsilon}
\newcommand{\RD}{\textup{(RD)}}

\newcommand{\ind}{\mathbf 1}
\newcommand{\norm}[1]{\lVert#1\rVert}
\newcommand{\papertitle}{Controlled analytic properties and the Baum--Connes conjecture}
\newcommand{\papershorttitle}{Controlled analytic properties and the Baum--Connes conjecture}
\makeatletter
\def\@adminfootnotes{%
  \let\@makefnmark\relax \let\@thefnmark\relax
  \ifx\@empty\@date\else \@footnotetext{\@setdate}\fi
  \ifx\@empty\@subjclass\else \@footnotetext{\@setsubjclass}\fi
  \ifx\@empty\@keywords\else \@footnotetext{\@setkeywords}\fi
  \ifx\@empty\thankses\else \@footnotetext{%
    \def\par{\let\par\@par}\@setthanks}%
  \fi
}
\makeatother

\setlist{nosep,leftmargin=2em}
\allowdisplaybreaks[2]
\title[\papershorttitle]{\papertitle}
\author{Vadim Alekseev}
\address{Vadim Alekseev, Institut f\"ur Geometrie, TU Dresden, 01062 Dresden, Germany}
\email{vadim.alekseev@tu-dresden.de}
\author{Martin Finn-Sell}
\address{Martin Finn-Sell}
\email{martin.r.sell@gmail.com}

\date{26 November 2026}
\begin{document}
\begin{abstract}
We develop a framework for transferring the reduced Baum--Connes assembly map
through asymptotically faithful coverings $G_m\to G$, using finite propagation
representatives in quantitative $K$-theory. Upward control comes from operator norm localisation of the approximants and holds unconditionally for lacunary hyperbolic groups. This allows to deduce the injectivity of Baum--Connes assembly map for every lacunary hyperbolic group, with all separable coefficients. Downward control comes from rapid decay property of the limit which allows to estimate norms of matrices over the group ring if one controls the lengths of entries and the injectivity radius. We deduce that a lacunary hyperbolic group with rapid decay and fast-growing injectivity radii satisfies the ordinary reduced Baum--Connes conjecture in both degrees. In an appendix we prove rapid decay for sufficiently sparse graphical small
cancellation presentations with controlled overlap cardinality, and realise
this overlap condition by a finite-coordinate refinement of Osajda's labelling. For a suitable subsequence of a large-girth expander family this gives groups with isometrically embedded expanders which satisfy the Baum--Connes conjecture.
\end{abstract}
\maketitle
\thispagestyle{empty}

\setcounter{tocdepth}{1}
\tableofcontents
\vspace*{-2em}
\section{Introduction}\label{sec:intro}
This paper concerns the Baum--Connes conjecture for a specific subclass of lacunary hyperbolic
groups. This class has appeared throughout the literature at the limits of geometric group theory,
and has a thorough treatment in the work of Ol'shanski\u{\i}--Osin--Sapir \cite{OOS09}. What led the
authors to consider this class is the following question: ``What are discrete groups that
contain expanders in their Cayley graph like?''

Small cancellation provides a flexible source of such groups. Gromov's construction
\cite{Gro03}, discussed and expanded by Arzhantseva--Delzant \cite{AD08}, uses geometric
small cancellation. Graphical small cancellation gives a different construction, in which
suitable labellings of finite graphs produce isometric embeddings of those graphs in the
resulting Cayley graph \cite{Oll06,Gru15,Osa20}. Gruber--Sisto distinguish their random Gromov monsters from infinitely
presented graphical small cancellation groups by quasi-isometry invariants \cite{GS21}.
We work with graphical presentations, and with sufficiently sparse subsequences of their
defining graphs. It is this additional sparsity that gives the lacunary hyperbolic
approximations used below.

The labelling method also applies to graphs with a wall structure. Arzhantseva--Osajda
\cite{AO14} use this to construct groups without property A which nevertheless have the
Haagerup property. For those groups the Higson--Kasparov theorem \cite{HK01} already gives
the Baum--Connes conjecture with coefficients. Our interest here is in groups containing
isometrically embedded expanders. Such groups do not coarsely embed into Hilbert space,
and the expander construction gives counterexamples to the Baum--Connes conjecture with coefficients
\cite{HLS02,Osa20}. This leaves a different question: whether their ordinary reduced
assembly map, with trivial coefficients $\C$, is an isomorphism.

The method we introduce is inspired by two sources: the work of Willett--Yu \cite{WY12}
on the coarse Baum--Connes conjecture for large-girth graph sequences, and the work of
Yu \cite{Yu98} on groups of finite asymptotic dimension. The former uses an
asymptotically faithful covering sequence, recalled in Section~\ref{sec:afc}, to lift
finite propagation operators. The latter uses a quantitative version of $K$-theory for
algebras equipped with a propagation filtration. This theory has been developed
systematically by Oyono-Oyono--Yu \cite{OOY15}.

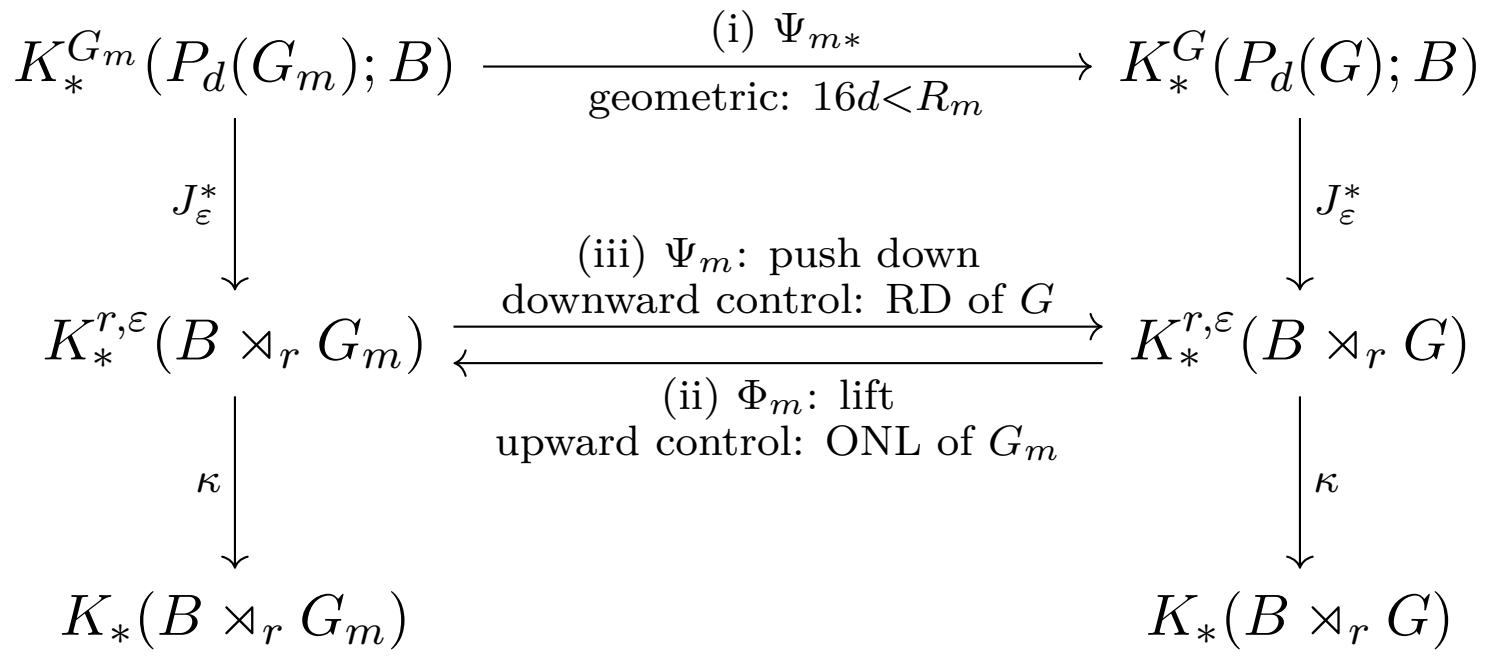
\begin{figure}[!ht]
\centering
\begin{tikzcd}[column sep=10.5em,row sep=3em]
K^{G_m}_*(P_d(G_m);B)\arrow[r,"\text{(i) }\Psi_{m*}","\scriptstyle\text{geometric: }16d<R_m"']\arrow[d,"J^*_\eps"'] &
K^{G}_*(P_d(G);B)\arrow[d,"J^*_\eps"]\\
K_*^{r,\eps}(B\rtimes_rG_m)\arrow[r,shift left=1.3,"\substack{\text{(iii) }\Psi_m\text{: push down}\\ \text{downward control: RD of }G}"]\arrow[d,"\kappa"'] &
K_*^{r,\eps}(B\rtimes_rG)\arrow[l,shift left=1.3,"\substack{\text{(ii) }\Phi_m\text{: lift}\\ \text{upward control: ONL of }G_m}"]\arrow[d,"\kappa"]\\
K_*(B\rtimes_rG_m) & K_*(B\rtimes_rG)
\end{tikzcd}
\setlength{\captionindent}{0pt}
\setlength{\abovecaptionskip}{6pt}
\makeatletter
\renewcommand{\@captionfont}{\normalfont\linespread{1}\footnotesize}
\makeatother
\caption{The maps in the proofs, for $m$ large compared with $d$ and $r$. The vertical
composites are the assembly maps $\mu_d$ (Section~\ref{sec:bc}): the polynomial index
$J^*_\eps$ of a normalised cycle lands in quantitative $K$-theory of the Roe algebra
$\K(L_d)\otimes(B\rtimes_rG)$ (the factor $\K(L_d)$ is suppressed) at a propagation $r$
depending only on $d$ and $\eps$, and $\kappa$ is the functional-calculus map.
The horizontal maps have the following controls:
(i) transport of classes is defined when $16d<R_m$ and commutes with
$J^*_\eps$ (Lemmas~\ref{lem:transport} and~\ref{lem:transport-classes});
(ii) lifting is multiplicative below the isometry radius, and its norm is controlled
from above by operator norm localisation of the approximants, giving weak
equivalence of the asymptotic lift with the regular representation
(Theorem~\ref{thm:weak-equivalence}); this holds for every lacunary hyperbolic group;
(iii) for $B=\C$, pushing down is a $*$-homomorphism on group rings; the
self-adjoint defects of a finite homotopy of propagation $t$ have propagation at
most $2t$ and active matrix size $N$. Their norm control comes from rapid decay
of $G$ when $2t\lceil\log(2eNt)\rceil\le R_m$ (Proposition~\ref{prop:rd-pushdown}). Tolerances change by
bounded factors along (ii) and (iii).}
\label{fig:directions}
\end{figure}

Our strategy is to cut the $K$-theory group into pieces that can be lifted through maps
which are multiplicative on a given scale, to solve the assembly problem upstairs, and
to push the solution back down. Figure~\ref{fig:directions} shows the maps involved.
We construct these maps on elements of bounded propagation and control their norms
at the scales needed for the argument. Quantitative $K$-theory records these scales
and the defects of the representatives and homotopies being transported.

\clearpage
To prove injectivity, suppose that $\mu_G[F]=0$. A null-homotopy of
$J^*_\eps(F)$ of some finite propagation lifts by (ii) to a null-homotopy of the
polynomial index of the lifted cycle $F_m$. Hyperbolic injectivity kills $[F_m]$ in the
top left corner, and (i) carries this vanishing back to $[F]$.

For surjectivity we take $B=\C$. A class in
$K_*(\Cr G)$ is represented in the middle row and lifted by (ii). Hyperbolic
surjectivity gives a cycle $F_m$ and a homotopy, of some propagation $t$ and matrix
size $N$, from the lift to $J^*_\eps(F_m)$. This homotopy is pushed down by (iii), and
$F_m$ by (i).

We first establish the representation-theoretic estimate that controls lifting
in both arguments. Lifting coefficients to $G_m$ gives a $*$-homomorphism
$\Phi_\infty\colon\C G\to\prod_m\Cr G_m/\bigoplus_m\Cr G_m$, defined at each
propagation by the coordinates with $m$ large. Let $\lambda_\infty$ be the
corresponding representation of $\C G$.

\begin{introtheorem}[Theorem~\ref{thm:weak-equivalence}]\label{thm:intro-weak-equivalence}
Let $\{G_m\to G\}_m$ be an asymptotically faithful covering sequence of marked groups.
Then $\lambda_\infty$ weakly contains the regular representation $\lambda_G$. If the
approximants have operator norm localisation with asymptotic controls, in particular
if $G$ is lacunary hyperbolic, then $\lambda_\infty$ and $\lambda_G$ are weakly
equivalent. In that case, for every separable $G$-$C^*$-algebra $B$, lifting extends
to an injective $*$-homomorphism
\[
 B\rtimes_rG\longrightarrow\prod_mB\rtimes_rG_m\Big/\bigoplus_mB\rtimes_rG_m,
\]
and for each fixed propagation it induces maps of quantitative $K$-theory groups into
those of $B\rtimes_rG_m$, for all sufficiently large $m$.
\end{introtheorem}

The last assertion allows us to work over individual approximants. A homotopy of finite
propagation lifts coordinate by coordinate to homotopies over
$B\rtimes_rG_m$ for all large $m$, with controlled defects. For lacunary hyperbolic
groups the operator norm localisation of the approximants is unconditional
(Theorem~\ref{thm:onl-lacunary}): it comes from finite asymptotic dimension of the
hyperbolic approximants.

\begin{introtheorem}[Theorem~\ref{thm:main-inj}]\label{thm:intro-injectivity}
For every lacunary hyperbolic group $G$ and every separable $G$-$C^*$-algebra $B$, the
reduced Baum--Connes assembly map with coefficients in $B$ is injective in both degrees.
\end{introtheorem}

The strong Novikov conjecture for lacunary hyperbolic groups is known to hold by
\cite[\S1.1]{OOS09}. The proof combines Theorem~\ref{thm:intro-weak-equivalence} with transport of a vanishing class at the
larger Rips scale where hyperbolic assembly is injective.

For surjectivity we also need to control pushdown in (iii). Its norm estimate will
use rapid decay of the limit group and depend on the propagation $t$ and total matrix
size $N$ of the homotopy. We therefore compare the growth of the radii $R_m$ with
$t$ and $N$; this comparison
is the fast-growing injectivity radii condition \eqref{eq:fir} below. Throughout the paper, quantitative
homotopies are recorded as actual paths, with their propagation, matrix size, and
defects specified. This is particularly useful in odd degree, where the standard
quantitative equivalence relation enlarges the homotopy parameters. The following
framework isolates what is needed: upward control, a pushdown bound on self-adjoint
defects of propagation at most $h(m,N)$ and matrix size $N$, and solutions upstairs
whose propagation and size fit under $h$.

\begin{introtheorem}[Theorem~\ref{thm:qs-framework}]\label{thm:intro-surjectivity}
Let $\{G_m\to G\}_m$ be an asymptotically faithful covering sequence with operator norm
localisation with asymptotic controls. Suppose that the approximants have surjectivity
with asymptotic controls in the sense of
Definition~\ref{def:asym-surj}. Then the assembly map for $G$ is quantitatively
surjective with size-dependent controls.
\end{introtheorem}

Theorem~\ref{thm:intro-surjectivity} gives a sufficient condition for surjectivity of the assembly map.
The relevant definitions and the proof are given in Section~\ref{sec:bc}.
To deduce the theorem, we fix the propagation,
matrix size and tolerance of the inputs. The hypothesis then supplies solutions
on arbitrarily late approximants. For each such approximant, one Rips scale and
one pair of propagation and size bounds work for every input with the fixed
parameters. These bounds allow
both the cycles and the homotopies to descend to $G$.

Section~\ref{sec:lacunary} applies this sufficient condition to lacunary hyperbolic groups,
leading to Theorem~\ref{thm:intro-rate} below. We now describe how the surjectivity hypothesis of
Theorem~\ref{thm:intro-surjectivity} is verified in that setting. Rapid decay of $G$ supplies the pushdown control
$(K,h)$, and the fast-growing injectivity radii condition ensures that the
witnesses supplied by hyperbolic assembly fit within this control.
To obtain the norm estimate, let $\pi\colon H\to G$ be a marked quotient with
isometry radius $R$, and suppose
that $G$ has ball-form rapid decay with constants $(C,s)$
(Definition~\ref{def:rd}). Proposition~\ref{prop:rd-pushdown} shows that, for
$b=b^*\in M_N(\C H)$ supported in $B_H(e,r)$ with $r\ge1$,
\[
 r\lceil\log(eNr)\rceil\le R
 \quad\Longrightarrow\quad
 \|\pi(b)\|\le K\|b\|,\qquad K=Ce^{2s+1/2}.
\]
We apply this estimate to the defects $p^2-p$ of quasi-projections and to both
defects $u^*u-1$ and $uu^*-1$ of quasi-unitaries. All these defects are
self-adjoint. If a homotopy has propagation at most $t_m$ and total active matrix
size at most $N_m$, its defects have propagation at most $2t_m$ and the same size.
Thus their norms increase by at most a factor of $K$ under pushdown, provided
their propagation stays within the radius allowed by the preceding estimate.

The corresponding cycle on $P_{k_m}(G_m)$ descends by the transport lemma when
$16k_m<R_m$. These two bounds are required in
Definition~\ref{def:asym-surj}. We therefore need the finite witnesses supplied
by hyperbolic assembly to satisfy them on arbitrarily late approximants.
The fast-growing injectivity radii condition below ensures this for every fixed
set of input parameters.

Fix a finite marking $S$. There is a finite
nondecreasing function $F_S\colon\N\times\N\to\N$, depending only on the size of the
marking, with the following role. Let $G=\varinjlim G_m$ be a lacunary presentation
with marking $S$ (Definition~\ref{def:lacunary}): the groups $G_m$ are
$\delta_m$-hyperbolic in the four-point sense, the limit maps $G_m\to G$ have isometry
radii $R_m\to\infty$, and $\delta_m=o(R_m)$. The presentation, with its data
$(G_m,\delta_m,R_m)_m$, has \emph{fast-growing injectivity radii} if
\begin{equation}\tag{FIR}\label{eq:fir}
 R_m\ge mF_S(\delta_m,m)\qquad\text{for infinitely many }m.
\end{equation}
A lacunary hyperbolic group \emph{has fast-growing injectivity radii} if it has a
lacunary presentation with this property; the presentation is then part of the data
(Definition~\ref{def:fir}). Definition~\ref{def:rate-function}
makes $F_S$ precise. Its value bounds the geometric and analytic requirements uniformly over all
$\delta$-hyperbolic groups with this marking and inputs of bounded complexity.
The complex $P_{d_\delta}(G_m)$, with
$d_\delta=32\delta+20$, is a model for the classifying space for proper actions
\cite{MS02}. Cycles and classes can be transported at this scale once
$16d_\delta<R_m$.
At this scale, an assembly preimage is witnessed by a finite
propagation homotopy. Since the Baum--Connes conjecture holds for hyperbolic
groups, such a homotopy exists. Compactness at fixed input size then gives
uniform finite bounds on its propagation $t$ and total matrix size $N$
(Lemma~\ref{lem:finite-witness}). There are only finitely many marked groups with a
fixed hyperbolicity bound (Lemma~\ref{lem:finite-markings}).
The rapid decay estimate then transports the self-adjoint defects of that
homotopy when $2t\lceil\log(eN\,2t)\rceil\le R_m$. Matrix size enters through its
logarithm.
Taking a maximum over these finite data defines $F_S$. Size-dependent existence
of these homotopies is itself equivalent to ordinary assembly surjectivity
(Proposition~\ref{prop:qs-ordinary}). The point of $F_S$ is its uniformity over the
approximating groups and its compatibility with the pushdown estimate.

\begin{introtheorem}[Theorem~\ref{thm:main-rate}]\label{thm:intro-rate}
Let $G$ be a lacunary hyperbolic group with property \RD, presented as the limit of
$\delta_m$-hyperbolic groups $G_m$ with a fixed finite marking $S$. If this presentation
has fast-growing injectivity radii \eqref{eq:fir}, then the ordinary
reduced Baum--Connes assembly map for $G$ is an isomorphism in both degrees.
\end{introtheorem}

Injectivity is Theorem~\ref{thm:intro-injectivity}. For surjectivity, we obtain solutions on the approximating
groups since the Baum--Connes conjecture holds for hyperbolic groups.
Condition \eqref{eq:fir} makes them fit within the pushdown radius of the rapid
decay estimate, so Theorem~\ref{thm:intro-surjectivity} applies.

Groups with fast-growing injectivity radii arise from small cancellation by
choosing each new relator, or relator graph, long enough compared with the group built
so far. For the expander examples one also needs rapid decay. Appendix~\ref{app:rd}
proves it for a specified finite-coordinate refinement of Osajda's labellings.

\begin{introtheorem}[Theorem~\ref{thm:c16} and Theorem~\ref{thm:main-expanders}]\label{thm:intro-examples}
\begin{enumerate}[label=(\roman*),ref=\roman*]
\item\label{thm:intro-classical} There are classical $C'(1/6)$ groups with fast-growing injectivity radii.
Every such group has quantitatively surjective assembly with size-dependent controls
in both degrees.
\item\label{thm:intro-expanders} There are graphical $C'(1/20)$ groups containing isometrically embedded
expanders, having ball-form rapid decay with exponent $13/2$ and fast-growing
injectivity radii. Their ordinary reduced assembly maps are isomorphisms in both degrees.
\end{enumerate}
\end{introtheorem}

The first part is already known as a statement about the Baum--Connes conjecture: classical $C'(1/6)$
groups have the Haagerup property \cite{AO15}, so satisfy the conjecture with
coefficients \cite{HK01}. It is useful here as an illustration of the transfer
argument. Rapid decay for all classical $C'(1/6)$ groups, including infinitely
presented ones, is now available from Bodart--Osin \cite{BO26}; it extends the earlier
$C'(1/10)$ theorem of Arzhantseva--Dru\c{t}u \cite{AD12}.

For the second part, Appendix~\ref{app:rd} first proves a uniform convolution estimate under a bound on
the number of vertices in overlaps of defining graphs, and then realises that bound
by adjoining one finite label coordinate. Rapid decay and \eqref{eq:fir} are both
sparsity conditions on the same sequence of girths, and one subsequence meets both.
For these groups Theorem~\ref{thm:intro-injectivity} gives injectivity with all
separable coefficients, while the expander obstruction to surjectivity with some
coefficients remains \cite[Corollaries 3.3--3.4]{Osa20}.

Theorem~\partref{thm:intro-examples}{\ref{thm:intro-expanders}} is particularly striking in view of the irrational-trace
construction of OpenAI \cite[Theorem~1.1]{OAI26}. Starting with suitable
isometrically embedded regular expanders in a graphical group $G$, that
construction produces an extension $V\rtimes G$, where $V$ is a quotient of
the lamp module $\mathbb F_p^{(G)}$. Thus $V\rtimes G$ is a quotient of the
wreath product $\mathbb F_p\wr G$. Its reduced group $C^*$-algebra contains
a projection of irrational canonical trace, whereas assembly classes have
rational trace, so its ordinary assembly map is not surjective. The underlying
groups $G$ can be chosen to satisfy the assumptions of Theorem~\partref{thm:intro-examples}{\ref{thm:intro-expanders}}, as explained
in Remark~\ref{rem:irrational-traces}.
Consequently, these groups can themselves satisfy the ordinary
Baum--Connes conjecture, even though their isometrically embedded expanders supply the geometry for irrational traces and thus failure of the Baum--Connes conjecture for these semidirect products.

\subsection*{Appendices}
Two appendices provide the auxiliary results used in the main argument.

Appendix~\ref{app:bc}, \emph{Folklore on the Baum--Connes assembly map}, contains several standard facts about the assembly map in precise forms:
every equivariant Kasparov class is represented on one standard module by an
exactly equivariant operator, using Meyer \cite{Mey00} and Gabe--Szab\'o
\cite{GS25} (Lemma~\ref{lem:ample-rep}); the resulting index map is
Kasparov assembly, including separable coefficients
(Lemma~\ref{lem:roe-coeff}, after Roe \cite{Roe02}); the free quotient by a
normal subgroup induces a map on classes and homotopies
(Lemma~\ref{lem:transport-classes}); and equalities in ordinary
$K$-theory are witnessed by finite propagation paths, with bounds uniform on
compact sets of inputs.

Appendix~\ref{app:rd}, \emph{Rapid decay and the refined graphical
construction}, proves rapid decay for sufficiently sparse graphical small
cancellation presentations with controlled overlap cardinality. The overlap
bound is constructed by a finite-coordinate refinement of Osajda's labelling.
This supplies the rapid decay hypothesis of Theorem~\ref{thm:intro-rate} for the expander examples
in Theorem~\partref{thm:intro-examples}{\ref{thm:intro-expanders}}.

The finite-matrix pushdown estimate and the construction in
Appendix~\ref{app:rd} are proved in full.

\subsection*{Conventions}
All groups $G$ in this paper are discrete and finitely generated. A marking $S$ is a fixed finite list of
generators; word length uses $S$ and its inverses, and marked homomorphisms preserve
this list. All logarithms are natural. The radii and
hyperbolicity constants can be enlarged or decreased to positive integers as
appropriate. The Baum--Connes assembly map is usdd in the form
\begin{equation}\label{eq:assembly}
 \mu_{G,B}\colon K_*^{\mathrm{top}}(G;B)
 =\varinjlim_d KK^G_*(C_0(P_d(G)),B)
 \longrightarrow K_*(B\rtimes_rG),\qquad *=0,1.
\end{equation}
Here $P_d(G)$ is the Rips complex of finitely supported probability measures whose
supports have diameter at most $d$. The case $B=\C$, with trivial action, is denoted
by $\mu_G$ and called ordinary reduced assembly.

\subsection*{History of the project and AI assistance}
The general ideas leading to this project emerged around 2018; unfortunately the first attempts to prove the main results for lacunary hyperbolic groups contained a gap. Between 2022 and 2024 the authors have closed this gap, proving Theorems~\ref{thm:intro-weak-equivalence} through~\partref{thm:intro-examples}{\ref{thm:intro-classical}}, but lacked a construction in Theorem~\partref{thm:intro-examples}{\ref{thm:intro-expanders}} which would give interesting examples with isometrically embedded expanders. The construction in Appendix~\ref{app:rd}, leading to Theorem~\partref{thm:intro-examples}{\ref{thm:intro-expanders}}, was produced with AI assistance using geometric ideas supplied by the authors; AI tools were further used to check the correctness of the proofs and to generally improve the manuscript. The large language models used in the process were Claude Opus 5.5, Claude Fable 5.2 and GPT-6 Astra.

\section{Asymptotically faithful covering spaces}\label{sec:afc}

\begin{definition}\label{def:radius}
Let $\pi\colon X\to Y$ be a surjective map of metric spaces. A number $R>0$ is an
\emph{isometry radius} if every ball of radius $R$ in $X$ maps isometrically onto the
corresponding ball in $Y$. The supremum of these numbers is the \emph{injectivity
radius} of $\pi$. We work with chosen isometry radii below this supremum whenever an
endpoint could matter.
\end{definition}

\begin{definition}
Two metric families $\mathcal X=\{X_m\}_m$ and $\mathcal Y=\{Y_m\}_m$ form an
\emph{asymptotically faithful covering family} if there are surjective local
isometries $\pi_m\colon X_m\to Y_m$ with isometry radii $R_m\to\infty$.
When deck-invariant operators are used, the coverings are regular: the deck group
acts transitively on every fibre.
\end{definition}

This notion was introduced by Willett--Yu \cite{WY12} for box spaces and sequences
of graphs of large girth. In this paper the target family is constant, and the
maps of interest are marked group quotients. These are regular coverings in the
local metric sense just defined. Their induced maps of Rips complexes are genuine
topological coverings at the scales specified in Lemma~\ref{lem:rips-covering}.

For a marked quotient $\pi\colon H\to G$, if every nonidentity element of the
kernel has length at least $L$, every radius $R$ with $4R<L$ is an isometry radius.
Indeed, two points in $B_H(e,R)$ have distance at most $2R$; a shorter word joining
their images would give a nontrivial kernel word of length at most $4R$. Surjectivity
onto the ball follows by lifting a geodesic word.

\begin{definition}\label{def:four-point}
A marked group is \emph{$\delta$-hyperbolic} if its vertex metric satisfies
\begin{equation}\label{eq:fourpoint}
 d(x,y)+d(z,t)\le
 \max\{d(x,z)+d(y,t),d(x,t)+d(y,z)\}+2\delta
\end{equation}
for every $x,y,z,t$. We take $\delta\in\N$, $\delta\ge1$.
\end{definition}

This is the convention of Meintrup--Schick \cite{MS02}, under which the Rips scale
$32\delta+20$ used below is valid. Equivalently, Gromov products at every basepoint
satisfy $(x|z)\ge\min\{(x|y),(y|z)\}-\delta$.

\begin{definition}\label{def:lacunary}
Let $G_1\to G_2\to\cdots$ be marked surjections, with direct limit $G$, and let
$\pi_m\colon G_m\to G$ be the limit maps. The presentation is \emph{lacunary} if
each $G_m$ is $\delta_m$-hyperbolic, the maps $\pi_m$ have isometry radii
$R_m\to\infty$, and $\delta_m=o(R_m)$. A group is \emph{lacunary hyperbolic} if
it has such a presentation.
\end{definition}

Existence of this presentation is equivalent to the usual definition, that some
asymptotic cone is an $\mathbb R$-tree \cite[Theorem~3.3]{OOS09}. The common marking
matters: both hyperbolicity and rapid decay constants depend on word length. All
our estimates concern maps to the limit. Consecutive maps can be used to control
these radii by the following observation.

\begin{lemma}\label{lem:inj-limit}
Let $G_1\to G_2\to\cdots$ be marked surjections with limit $G$. If $\iota_k$ is
an isometry radius for $G_k\to G_{k+1}$, then every
$R<\inf_{k\ge m}\iota_k$ is an isometry radius for $G_m\to G$.
\end{lemma}
\begin{proof}
Every composition from $G_m$ onwards is isometric on balls of radius $R$. A
relation in the direct limit already holds at some finite stage, so the distance
between the images of two fixed elements is eventually their distance in the
limit. Thus no distance in such a ball can decrease on passing to $G$.
\end{proof}

\section{Operator norm localisation with asymptotic controls}\label{sec:onl}
In this section we prove a basic representation-theoretic result concerning
asymptotically faithful covering sequences. First we recall the algebras attached
to a group or metric space.

\begin{definition}
Let $X$ be a uniformly discrete metric space of bounded geometry. An operator
$T$ on $\ell^2(X)\otimes\mathcal H$ has \emph{propagation at most $r$} if its
matrix entries vanish at pairs of points more than $r$ apart. We write
$\mathcal A[X]_r$ for these bounded operators, allowing an arbitrary separable
Hilbert fibre $\mathcal H$, and $\mathcal A[X]=\bigcup_r\mathcal A[X]_r$.
For $\mathcal H=\C$ its closure is the uniform Roe algebra. Requiring compact
matrix entries on a fixed infinite-dimensional fibre gives the algebraic Roe
algebra and its closure $C^*X$. The same definitions with invariant operators
give the equivariant Roe algebra. For a free, proper, cocompact action, its
standard stabilised version is Morita equivalent to the reduced group
$C^*$-algebra \cite{WY12}.
\end{definition}

For a group we also use the filtration
\[
 \C G_r=\{a\in\C G:\supp a\subseteq B_G(e,r)\}.
\]
It is the propagation filtration of the regular representation, and extends
to matrices by allowing a finite-dimensional Hilbert fibre.

\begin{definition}\label{def:filtered}
A \emph{filtered $C^*$-algebra} is a $C^*$-algebra $A$ with closed self-adjoint
linear subspaces $A_r$, increasing in $r$, such that
$A_rA_{r'}\subseteq A_{r+r'}$ and $\bigcup_rA_r$ is dense in $A$. In the unital
case $1\in A_r$ for all $r$. A map on propagation at most $R$ is an
\emph{$R$-quasi-$*$-homomorphism} if it is linear, preserves adjoints, and
preserves every product for which the sum of the propagations is less than $R$.
Continuity is understood on each specified propagation subspace.
\end{definition}

\begin{definition}[Lifting maps]\label{def:lifting}
Let $\pi\colon X\to Y$ be a regular covering with deck group $N$ and isometry
radius $R$. For $r<R$ and $T\in\mathcal A[Y]_r$, put
\[
 \Phi(T)_{x,x'}=
 \begin{cases}T_{\pi(x),\pi(x')},&d(x,x')\le r,\\0,&d(x,x')>r.
 \end{cases}
\]
This defines an $N$-invariant operator; its boundedness follows, for example,
from the finite row bound. It is a bijection onto the $N$-invariant operators
at this propagation, with inverse $\Psi$. These are the \emph{lifting} and
\emph{pushdown} maps. They preserve adjoints and products whenever the sum of
propagations is less than $R$.
\end{definition}

Regularity is used for the inverse: corresponding matrix entries on two sheets
must agree. On group rings $\Phi$ relabels coefficients by their unique short
lifts, and $\Psi$ is the quotient homomorphism. A vector supported in a ball
that lifts isometrically is lifted to one chosen sheet. This preserves its
$\ell^2$-norm.

Given an asymptotically faithful family, write $\Phi_m,\Psi_m$ for these maps.
They are \emph{asymptotically isometric} if
$\|\Phi_m(T)\|\to\|T\|$ for every fixed finite propagation operator $T$.
The analysis of their continuity is at the heart of this paper.

\begin{definition}\label{def:onl}
A discrete bounded-geometry space has operator norm localisation (ONL) with constant $c\in(0,1)$ and nondecreasing control $f$ if every propagation-$r$ operator $T$ on $\ell^2(X)\otimes\mathcal H$ admits a unit vector $\eta$ with
\[
 \diam\supp\eta\le f(r),\qquad \|T\eta\|\ge c\|T\|.
\]
The constants are uniform over the Hilbert fibre.
\end{definition}
\begin{lemma}\label{lem:onl-constant}
Suppose $X$ has ONL with constant $c$ and control $f$, and let $n\ge1$.
\begin{enumerate}
\item For self-adjoint $T$ of propagation at most $R$ there is a unit vector $\eta$ with $\diam\supp\eta\le2(n-1)R+f(nR)$ and $c^{1/n}\|T\|\le\|T\eta\|$.
\item $X$ has ONL with constant $c^{1/n}$ and control $g(R)=4(n-1)R+f(2nR)$.
\end{enumerate}
In particular $X$ has ONL with any constant $c'\in(0,1)$.
\end{lemma}

\begin{proof}
(1) is the argument of \cite[Proposition~2.4]{CTWY08}, which needs $\|T^n\|=\|T\|^n$ and so applies to self-adjoint $T$. Assume
$\|T\|=1$. Since $T^n$ has propagation at most $nR$, ONL gives $\varphi$ of support diameter at most $f(nR)$ with
$\|T^n\varphi\|\ge c\|\varphi\|$, so $\|T(T^j\varphi)\|\ge c^{1/n}\|T^j\varphi\|$ for some $0\le j\le n-1$. The support of $T^j\varphi$ lies in the
$jR$-neighbourhood of $\supp\varphi$, so its diameter is at most $2(n-1)R+f(nR)$. (The proof of \cite[Proposition~2.4]{CTWY08} states
$(n-1)R+f(nR)$; the factor $2$ is needed, and changes only constants.) (2) For arbitrary $T$ of propagation at most $R$, apply (1) to the self-adjoint $T^*T$, of propagation at most $2R$. This gives $\eta$
with $c^{1/n}\|T\|^2\|\eta\|\le\|T^*T\eta\|\le\|T\|\|T\eta\|$.
\end{proof}

\begin{definition}\label{def:asym-onl}
An asymptotically faithful family has \emph{ONL with asymptotic controls} if its approximants have a common ONL constant $c$ and nondecreasing controls $f_m$ such that
\[
 \rho_m:=\sup\{r:r+f_m(r)\le R_m\}\longrightarrow\infty.
\]
\end{definition}
\begin{lemma}\label{lem:cm-to-one}
Under the preceding hypothesis there are constants $c_m\in(0,1)$, $c_m\to1$, and radii
$S_m\to\infty$ such that $X_m$ has ONL with constant $c_m$ and control $g_m$, and $S_m+g_m(S_m)\le R_m$.
\end{lemma}

\begin{proof}
Choose $\rho'_m$ with $\rho'_m+f_m(\rho'_m)\le R_m$ and $\rho'_m\to\infty$, which is possible since $\rho_m\to\infty$, and integers $n_m\to\infty$ with $S_m:=\rho'_m/(4n_m)\to\infty$. By
Lemma~\partref{lem:onl-constant}{2}, $X_m$ has ONL with constant $c_m=c^{1/n_m}\to1$ and control $g_m(k)=4(n_m-1)k+f_m(2n_mk)$. Since $f_m$ is
non-decreasing, $S_m+g_m(S_m)\le4n_mS_m+f_m(4n_mS_m)=\rho'_m+f_m(\rho'_m)\le R_m$.
\end{proof}

\begin{theorem}\label{thm:asym-isometric}
Let $\{\pi_m\colon X_m\to Y\}$ be a regular asymptotically faithful covering family with ONL with asymptotic controls. Then the lifting maps $\Phi_m\colon\mathcal A[Y]\to\mathcal A[X_m]^{N_m}$ are asymptotically isometric.
More precisely, for $T\in\mathcal A[Y]_r$ and $m$ with $r\le S_m$,
\[
c_m\|\Phi_m(T)\|\le\|T\|.
\]
\end{theorem}

\begin{proof}
Let $T\in\mathcal A[Y]_r$ and $\eps>0$. There is a finitely supported unit vector $\xi\in\ell^2(Y)\otimes\mathcal H$ with
$\|T\xi\|\ge(1-\eps)\|T\|$. For $m$ large, $\xi$ and $T\xi$ lie in a ball of radius less than
$R_m/2$, so
\[
(1-\eps)\|T\|\le\|T\xi\|=\|\Phi_m(T)\Phi_m(\xi)\|\le\|\Phi_m(T)\|.
\]
Conversely, for $r\le S_m$, ONL of $X_m$ gives a unit vector $\eta_m$ with
$\diam\supp\eta_m\le g_m(r)$ and $c_m\|\Phi_m(T)\|\le\|\Phi_m(T)\eta_m\|$. Since
$r+g_m(r)\le S_m+g_m(S_m)\le R_m$, the vectors $\eta_m$ and $\Phi_m(T)\eta_m$ lie on one sheet, so
$\|\Phi_m(T)\eta_m\|=\|T\Psi_m(\eta_m)\|\le\|T\|$. As $c_m\to1$, the claim follows.
\end{proof}

Here $N_m$ is the deck group; for a marked quotient it is $\ker\pi_m$. The inequality is uniform over the entire propagation subspace, so it applies to every point of a norm-continuous path and to all matrix amplifications.

We write $\Phi_m^{(r)}$ and $\Psi_m^{(r)}$ for restrictions to propagation at
most $r$. The proof gives $\|\Phi_m^{(r)}\|\le c_m^{-1}$ for $r\le S_m$,
uniformly over matrix amplifications and Hilbert fibres.

\subsection{Equivariant lifting maps and representation theory}

Throughout this subsection $G$ is a finitely generated group with an asymptotically faithful
covering family $\{G_m\to G\}_m$, and we consider elements of $\C G_r$, i.e.\ the left-equivariant
operators of propagation at most $r$. For $m$ large, $\Phi_m^{(r)}$ is a bijection
$\C G_r\to(\C G_m)_r$ with inverse $\Psi_m^{(r)}$; this inverse is the restriction of the quotient
homomorphism $\pi_m$.

\begin{definition}
For $a\in\C G$ put $\|a\|_{(r)}:=\sup\{\|av\|_2\mid v\in\ell^2(B_r),\ \|v\|_2=1\}$.
\end{definition}

\begin{proposition}\label{prop:local-norm}
For every finitely generated group $G$ and $a\in\C G$, $\|a\|=\lim_{r\to\infty}\|a\|_{(r)}$.
\end{proposition}

\begin{proof}
The function $r\mapsto\|a\|_{(r)}$ is increasing and bounded by $\|a\|$. Given $\eps>0$, choose a finitely
supported unit vector $v$ with $\|av\|_2\ge(1-\eps)\|a\|$; its finite support lies in $B_r$ for some $r$. Then $(1-\eps)\|a\|\le\|a\|_{(r)}\le\|a\|$.
\end{proof}

\begin{remark}[Equivariant localisation in every group]\label{rem:weak-onl}
For each fixed $r$ the convergence in Proposition~\ref{prop:local-norm} is uniform on the unit sphere
of $\C G_r$: this sphere is compact, since $\C G_r$ is finite dimensional, and the functions
$a\mapsto\|a\|_{(R)}/\|a\|$ are continuous, increasing in $R$ and converge pointwise to $1$, so Dini's
theorem applies. Hence every finitely generated group satisfies: for every $c\in(0,1)$ there is a
function $f$ such that for all $a\in\C G_r$ there is a unit vector $v\in\ell^2(G)$ supported in
$B_{f(r)}$ with $c\|a\|\le\|av\|_2$. The same holds for $M_n(\C G)$ for each fixed $n$, with $f$ depending
on $n$.
\end{remark}

\begin{definition}[The lifting limsup representation]
Define a norm on $\C G$ by $\|a\|_\infty=\limsup_m\|\Phi_m(a)\|_{\Cr G_m}$, and let $\lambda_\infty$ be a faithful representation of its completion.
This is a $C^*$-norm: coefficient lifts give an algebraic $*$-homomorphism into
$\prod_m\Cr G_m/\bigoplus_m\Cr G_m$, whose norm is the displayed limsup.
The sequences are bounded by the coefficient $\ell^1$ norm. This norm is faithful on
$\C G$, by Theorem~\partref{thm:weak-equivalence}{1}.
\end{definition}

Theorem~\ref{thm:weak-equivalence} in Section~\ref{sec:weak-equivalence} states the
comparison of $\lambda_\infty$ with $\lambda$, with coefficients, and its consequence for
quantitative $K$-theory. We state it there, immediately before its use in the proof of
injectivity. The lower bound $\lambda\prec\lambda_\infty$ holds for every asymptotically
faithful family. For the upper bound we use ONL of the approximants with asymptotic
controls. In the hyperbolic case these controls come from finite asymptotic
dimension with linear control, as established in Theorem~\ref{thm:onl-lacunary}.

\subsection{Pushdown control}
The norms of the inverse maps $\Psi_m$ can behave very differently from
those of $\Phi_m$. We formulate the estimate needed for assembly separately.

\begin{definition}[Pushdown control]\label{def:pushdown}
Let $\{G_m\to G\}$ have isometry radii $R_m$. For a class
$\mathcal A_{m,N}$ of self-adjoint finite propagation operators on the
approximants, a \emph{pushdown control} consists of $K\ge1$ and radii
$0\le h(m,N)\le R_m/2$, with $h(m,N)\to\infty$ for each fixed $N$, such that
\[
 \|\Psi_m(T)\|\le K\|T\|
 \quad(T\in\mathcal A_{m,N},\ \prop(T)\le h(m,N)).
\]
The parameter $N$ will be the total active matrix size. The classes used
below contain the defects $p^2-p$, $u^*u-1$, and $uu^*-1$ of the relevant
homotopies. When the estimate holds for all operators and all Hilbert
fibres, the parameter $N$ is omitted.
\end{definition}

\begin{lemma}\label{lem:onl-pushdown}
Suppose the target space $Y$ has ONL with constant $c$ and control $f$.
Then $K=c^{-1}$ and
\[
 h(m)=\max\bigl(\{0\}\cup\{r\in\N:r+f(r)<R_m/2\}\bigr)
\]
give pushdown control for all deck-invariant finite propagation operators,
uniformly over the Hilbert fibre.
\end{lemma}
\begin{proof}
For $T$ of propagation $r\le h(m)$, ONL of $Y$ gives a unit vector $\eta$
of support diameter at most $f(r)$ with
$c\|\Psi_m(T)\|\le\|\Psi_m(T)\eta\|$. Lift $\eta$ to one sheet.
The propagation neighbourhood also lies on that sheet, so
\[
 c\|\Psi_m(T)\|\le\|\Psi_m(T)\eta\|
 =\|T\Phi_m(\eta)\|\le\|T\|.
\]
For each fixed $r$, the defining inequality for $h(m)$ holds eventually.
\end{proof}

\begin{remark}\label{rem:equivariant-only}
For a lacunary hyperbolic limit, a common pushdown bound for \emph{all}
finite propagation operators would imply ONL of the limit, hence property
A \cite{Sak14}. Indeed, fix $r$, take one sufficiently late approximant,
and apply its ONL to $\Phi_m(T)$ for every propagation-$r$ operator $T$.
Push the localising vector to the limit. Together with
$\|T\|\le K\|\Phi_m(T)\|$, this gives ONL with, say, constant $1/(2K)$
and a finite control depending on $r$. The choice of approximant is uniform over propagation-$r$ operators.
For the expander examples, Proposition~\ref{prop:rd-pushdown} supplies the
pushdown estimate on finite group-ring matrices that is used in assembly.
\end{remark}

\section{Quantitative assembly through coverings}\label{sec:bc}

\subsection{Quantitative \texorpdfstring{$K$}{K}-theory}
We use the conventions of \cite{OOY15}. An
\emph{$(r,\eps)$-quasi-projection} is a self-adjoint element $p\in M_n(A_r)$
with $\|p^2-p\|<\eps$. An \emph{$(r,\eps)$-quasi-unitary} is an element
$u\in M_n(A_r)$ with
\[
 \|u^*u-1\|<\eps,\qquad \|uu^*-1\|<\eps.
\]
For $\eps<1/4$, spectral projection and polar decomposition give ordinary
$K$-theory classes. The corresponding elements differ from $p$ and $u$ by at
most $2\eps$.

\begin{definition}[{\cite[Definition~1.13]{OOY15}}]\label{def:qk}
For unital $A$, the group $K_0^{r,\eps}(A)$ consists of pairs $(p,l)$ modulo
stabilisation and homotopy: $(p,l)$ and $(q,l')$ are equivalent if, after zero
padding, $p\oplus1_{k+l'}$ and $q\oplus1_{k+l}$ are joined by
$(r,\eps)$-quasi-projections for some $k$. The group $K_1^{r,\eps}(A)$ is formed
from quasi-unitaries, stabilised by identity blocks, with equivalence given by
a path of \emph{$(2r,3\eps)$-quasi-unitaries}. For nonunital $A$ one uses the
unitisation and, in degree zero, the kernel of the scalar map. Increasing $r$
and $\eps$ gives structure maps, and functional calculus gives maps to $K_*(A)$.
\end{definition}

In degree one, the equivalence relation at parameters $(r,\eps)$ uses
paths at parameters $(2r,3\eps)$. We therefore record the actual paths and
their bounds throughout the argument. In either degree, each defect has
propagation at most twice that of the path. Closed defect bounds $\le\eps$ are convenient in
compactness arguments; a strict enlargement of the tolerance gives the
usual open convention.

\subsection{Roe algebras of Rips complexes}\label{sec:cycles}
We use the notation
\[
 K^G_*(P_d(G);B)=KK^G_*(C_0(P_d(G)),B),
\]
and omit $B$ when $B=\C$.
Let $G$ be finitely generated, with the word metric of the fixed marking, and let
$P_d(G)$ be its Rips complex: the finitely supported probability measures on $G$ with support of diameter
at most $d$, with the topology of \cite[\S6.1]{OOY15}. We fix once and for all a countable dense subset $W\subseteq(0,1)$
that is algebraically independent over $\mathbb Q$. Let $Z_d\subseteq P_d(G)$ be the set of measures $\sum_{i=1}^kw_i\delta_{g_i}$, with $g_i$
distinct, whose weights $w_1,\dots,w_{k-1}$, in some order, are distinct elements of $W$ with sum less than $1$, and
$w_k=1-\sum_{i<k}w_i$. By algebraic independence the $k$ weights are pairwise distinct. So an element of $G$ fixing such a
measure fixes each support point, and $G$ acts freely on $Z_d$. Moreover $Z_d$ is countable, dense and $G$-invariant,
and $Z_d\subseteq Z_{d'}$ for $d\le d'$. Let $D_d\subseteq Z_d$ consist of the measures whose heaviest point is $e$. This is a
fundamental domain, its measures are supported in $B(e,d)$, and $D_d\subseteq D_{d'}$. These choices depend only on
$(G,S)$, and push-forward along a homomorphism that is injective on supports preserves weights. Hence it maps $Z_d$
into $Z_d$ and $D_d$ into $D_d$. Put $L_d=\ell^2(D_d)\otimes H_0$. Then
\[
\ell^2(Z_d)\otimes H_0=\ell^2(G)\otimes L_d,
\]
with $G$ acting by left translation on the first factor. A bounded $G$-equivariant operator of finite
propagation is a finite sum $T=\sum_\gamma\rho_\gamma\otimes T_\gamma$, with $\rho$ the right regular
representation and $T_\gamma\in\B(L_d)$; we measure its propagation by $\max\{|\gamma|:T_\gamma\ne0\}$. This is
the word pseudometric obtained from the orbit coordinate of $Z_d$, and agrees with the propagation for the metric of $P_d(G)$ up
to constants depending on $d$. Such $T$ is locally compact exactly when all $T_\gamma$ are compact, since
$D_d$ has compact closure. Hence:

\begin{lemma}\label{lem:roe-rips}
The equivariant Roe algebra of $P_d(G)$ is $C^*(P_d(G))^G=\Cr G\otimes\K(L_d)$, filtered by
$\C G_r\odot\K(L_d)$. For $d\le d'$ the map induced by $P_d(G)\subseteq P_{d'}(G)$ is $\mathrm{id}\otimes$ (the
inclusion $\K(L_d)\subseteq\K(L_{d'})$).\qed
\end{lemma}

\subsection{Assembly as an index map}
We represent $K^G_0(P_d(G))$ on the module $\ell^2(Z_d)\otimes H_0$, which is ample because $Z_d$ is dense and every point evaluation occurs with infinite multiplicity through $H_0$. Let $\phi\colon P_d(G)\to[0,1]$ be $\phi(\nu)=\nu(e)$. It is continuous (it is $1$-Lipschitz for the metric of \cite[\S6.1]{OOY15}), its support is the closed star of $e$, contained in the
compact set $V_d$ of measures with support in $B(e,d)$ and hence in the set $W_d$ of measures with support in $B(e,2d)$, and
$\sum_{g\in G}(g\phi)(\nu)=\sum_g\nu(g)=1$; it is the barycentric coordinate $\lambda_e$ of \cite[\S6.1]{OOY15}. It is canonical, the same for every $d$, and preserved by
push-forward along maps injective on supports.

\begin{definition}
A \emph{normalised cycle} on $P_d(G)$ is an operator $F=\sum_{g\in G}g\sqrt\phi\,T\sqrt\phi\,g^{-1}$, with $T$ a
bounded operator of norm at most $1$, such that $F^*F-1$ and $FF^*-1$ are locally compact, and $[F,f]$ is compact for every $f\in C_0(P_d(G))$.
\end{definition}

A normalised cycle is $G$-equivariant. With the isometry $V\xi=(\sqrt{g\phi}\,\xi)_g$ into
$\bigoplus_g\ell^2(Z_d)\otimes H_0$ one has $F=V^*\bigl(\bigoplus_g gTg^{-1}\bigr)V$, so $\|F\|\le1$. Each term is supported on $(gW_d)^2$, whose points have supports in $B(g,2d)$, so $F$ has
propagation at most $6d$. Every class in $K^G_0(P_d(G))$ is represented by a normalised cycle: take a
representative $F'$ on this module with $\|F'\|\le1$ (Lemma~\ref{lem:ample-rep}), and average, $F=\sum_gg\sqrt\phi F'\sqrt\phi g^{-1}$; then $F-F'$ is
locally compact \cite[Definition~4.2]{WY12}, \cite[\S3]{Roe02}.

\newcommand{\standardmodulelemmastatement}{%
Let $H_0$ be a separable Hilbert space, graded with even and odd parts both infinite-dimensional, and let $(B,\beta)$ be a
separable $G$-$C^*$-algebra. Let $\mathcal H_B=\ell^2(Z_d)\otimes H_0\otimes B$, with the multiplication representation $M$ of $C_0(P_d(G))$ and the
diagonal action $\gamma_g=g\otimes1\otimes\beta_g$. Every class in $KK^G_*(C_0(P_d(G)),B)$ is represented by a cycle $(\mathcal H_B,M,F)$ with $F$ an
exactly equivariant self-adjoint contraction. In particular, every class in $K^G_*(P_d(G))$ is represented by a normalised cycle.
}
\begin{lemma}[Representatives on the standard module]\label{lem:ample-rep}
\standardmodulelemmastatement
\end{lemma}

The ungraded normalised cycles of the definition correspond to graded ones as follows. Identify the two parts of $H_0$ and
write an odd operator as $\left[\begin{smallmatrix}0&F_-\\F_+&0\end{smallmatrix}\right]$; the normalised cycle is $F_+$. Averaging commutes with this dictionary.
The proof is given in Appendix~\ref{app:bc-rep}.

For bounded $F$ let
\[
I(F)=\begin{pmatrix}FF^*+(1-FF^*)FF^*&F(1-F^*F)+(1-FF^*)F(1-F^*F)\\(1-F^*F)F^*&(1-F^*F)^2\end{pmatrix},
\]
the conjugate of $\left[\begin{smallmatrix}1&0\\0&0\end{smallmatrix}\right]$ by an invertible matrix polynomial in
$F$; it is an idempotent. For a normalised cycle, $\|I(F)\|\le6$, $I(F)-\left[\begin{smallmatrix}1&0\\0&0\end{smallmatrix}\right]$ lies
in $M_2(C^*(P_d(G))^G)$, and has propagation at most $30d$. The index map
\[
\mu_d\colon K^G_0(P_d(G))\to K_0(C^*(P_d(G))^G)\cong K_0(\Cr G),\qquad[F]\mapsto[I(F)]-\left[\begin{smallmatrix}1&0\\0&0\end{smallmatrix}\right],
\]
is the assembly map in the index picture of Willett--Yu \cite[Definition~4.2]{WY12}: the same idempotent $I(F)$, for
the same averaged operator $F$, with the equivariant version described after their Definition~4.3 (for free actions). In their book it is
the evaluation-at-one assembly map \cite[Definition~7.1.1]{WY20}, composed with the map from Kasparov's cycles
given by averaging over partitions of unity and taking the boundary map \cite[Construction~6.7.3]{WY20}. The
same construction with coefficients in a $G$-$C^*$-algebra $A$, on the Hilbert $A$-module $\ell^2(Z_d)\otimes H_0\otimes A$, gives
$\mu_d\colon KK^G_0(C_0(P_d(G)),A)\to K_0(C^*(P_d(G),A)^G)\cong K_0(A\rtimes_rG)$. Here $C^*(P_d(G),A)^G$ is the closure of the invariant, locally
compact, finite propagation operators on this module; it is $\K(L_d)\otimes(A\rtimes_rG)$, filtered by $\K(L_d)\odot C_c(G,A)_r$, where $C_c(G,A)_r$ is
spanned by the $au_\gamma$ with $|\gamma|\le r$ (Appendix~\ref{app:bc-index}). For separable $A$ every class has a
representative on this module (Lemma~\ref{lem:ample-rep}).

\newcommand{\assemblyindexlemmastatement}{%
Let $B$ be a separable $G$-$C^*$-algebra. Then $\varinjlim_d\mu_d$ is the Baum--Connes assembly map with
coefficients in $B$, defined through Kasparov's descent and the basic projection. The same holds in degree $1$, with
$\exp(-2\pi iP_F)$ in place of $I(F)$ (Section~\ref{sec:degree-one}). In particular $G$ satisfies the Baum--Connes conjecture
with coefficients in $B$ if and only if $\varinjlim_d\mu_d$ is an isomorphism for $B$.
}
\begin{lemma}\label{lem:roe-coeff}
\assemblyindexlemmastatement
\end{lemma}

For $B=\C$ this is Roe's theorem \cite{Roe02}; the proof with coefficients is given in Appendix~\ref{app:bc-index}.

\subsection{Polynomial representatives in both degrees}\label{sec:degree-one}
For an even normalised cycle, put $e=\left[\begin{smallmatrix}1&0\\0&0\end{smallmatrix}\right]$. There is a self-adjoint noncommutative polynomial $q_\eps$ such that
\begin{equation}\label{eq:even-poly}
 q_\eps(e)=e,\qquad
 \|q_\eps(I(F))-\theta(I(F))\|<\eps/4,
 \quad \theta(a)=aa^*(1+(a-a^*)(a^*-a))^{-1}.
\end{equation}
Indeed, for an idempotent $a$ with $a^2=a$ and $\|a\|\le6$, the element $\theta(a)$ is a projection, and the positive element in the inverse has spectrum in $[1,145]$. Approximate its inverse uniformly by a real polynomial $p$ with $p(1)=1$, and take the self-adjoint part of $aa^*p(1+(a-a^*)(a^*-a))$. The interpolation at $1$ can be imposed by subtracting the scalar error of an initial approximation. For every such idempotent $a$, the resulting $q_\eps(a)$ is an $\eps$-almost projection, homotopic through such almost projections to $\theta(a)$. Exact interpolation makes the scalar part of $q_\eps(I(F))$ exactly $e$. Its non-scalar part has compact coefficients and propagation at most $30\deg(q_\eps)d$.

An odd normalised cycle has self-adjoint $T$ and $F$. Set $P_F=(F+1)/2$, so $0\le P_F\le1$. Its index unitary is $\exp(-2\pi iP_F)$, whose difference from $1$ lies in the Roe algebra. Choose a polynomial
\begin{equation}\label{eq:odd-poly}
 w_\eps(x)=1+x(1-x)k_\eps(x),\qquad
 \sup_{0\le x\le1}|w_\eps(x)-e^{-2\pi ix}|<\eps/3.
\end{equation}
This follows by polynomial approximation of the continuous, endpoint-extended function $(e^{-2\pi ix}-1)/(x(1-x))$. Both defects of $w_\eps(P_F)$ have norm less than $\eps$, and the straight line to the exponential consists of $\eps$-almost unitaries. Its scalar part is $1$, its non-scalar coefficients are compact, and its propagation is at most $6\deg(w_\eps)d$.

We denote these polynomial index representatives by $J^0_\eps(F)=q_\eps(I(F))$ and $J^1_\eps(F)=w_\eps(P_F)$; in degree zero the scalar projection $e$ is subtracted in $K$-theory. They are functorial under algebraic $*$-homomorphisms. Their defects are self-adjoint in both parities.

For $X=\sum_g\rho_g\otimes X_g$ with finite support, let $\|X\|_1=\sum_g\|X_g\|$. A finite-rank fibre projection $Q$ can compress the non-scalar part of $J^*_\eps(F)$ with arbitrarily small $\|\cdot\|_1$ error. This follows by approximating each of its finitely many compact coefficients. A common $Q$ works for a compact path of such coefficients. We use these finite-dimensional compressions for $B=\C$, and norm approximation for coefficient injectivity below.

\subsection{Coverings and transport of cycles}\label{sec:transport}
Let $\pi_m\colon G_m\to G$ be surjective with isometry radius $R_m$ and kernel $K_m$. Nontrivial elements
of $K_m$ have length greater than $2R_m$, since $\pi_m$ is injective on $B(e,R_m)$.

\begin{lemma}\label{lem:rips-covering}
If $4d<R_m$, then $K_m$ acts freely on $P_d(G_m)$, and the push-forward of measures $P_d(G_m)\to P_d(G)$ is a
$G_m$-equivariant covering whose fibres are the $K_m$-orbits.
\end{lemma}

\begin{proof}
If $k\in K_m$ fixes $\nu$, it permutes $\supp\nu$; for $v\in\supp\nu$ the element $v^{-1}kv\in K_m$ has length at most
$d<R_m$, so $k=e$. Every finite set of diameter at most $d$ in $G$ lifts isometrically. If two measures have the
same image, translate by an element of $K_m$ so that their supports share a point $v$; both supports lie in
$B(v,d)$, on which $\pi_m$ is injective, so the measures agree. To see the covering property, take a measure $\nu$ and a vertex $v$ of positive weight. In a sufficiently small neighbourhood that weight remains positive, and all supports lie in $B(v,d)$. The unique lifts in $B(\tilde v,d)$ give a continuous inverse on this neighbourhood. Its other inverse images are its disjoint kernel translates, by the fibre description.
\end{proof}

If $4d<R_m$, then push-forward is injective on supports of diameter at most $d$, so $Z_d(G_m)$ is the preimage of $Z_d(G)$
and $D_d(G_m)$ is the lift of $D_d(G)$ with supports in $B(e,d)$. The choices made for $G_m$ alone and for the covering
therefore agree. Then $L_d$ is the same space for $G_m$ and for $G$, and we define the
\emph{pushdown}
\[
\Psi_m=\pi_m\otimes\mathrm{id}\colon\C G_m\odot\B(L_d)\to\C G\odot\B(L_d).
\]
It is a $*$-homomorphism. On elements of propagation less than $R_m/2$ it is the pushdown of
the local lifting construction of Section~\ref{sec:onl}, with inverse the lifting $\Phi_m$; on group rings it is $\pi_m$. On elements of
propagation at most $R_m$ it only relabels coefficients, since $\pi_m$ is injective on $B(e,R_m)$; so it preserves
$\|\cdot\|_1$ and commutes with compressions $1\otimes P$. The upward
control of Theorem~\ref{thm:asym-isometric} applies to $\Phi_m$ with $Y=G$, $X_m=G_m$ and coefficients
$\K(L_d)$. The function $\phi_m(\nu)=\nu(e)$ on $P_d(G_m)$ is the lift of $\phi$ to the sheet over $W_d$ through $e$, and
$\sum_{g\in G_m}g\phi_m=1$. So normalised cycles for $G_m$ alone and for the covering are the same.

\begin{lemma}[Transport of cycles]\label{lem:transport}
Let $16d<R_m$ and let $F_m=\sum_{g\in G_m}g\sqrt{\phi_m}T\sqrt{\phi_m}g^{-1}$ be a normalised cycle on $P_d(G_m)$.
Then $\bar F=\Psi_m(F_m)$ is a normalised cycle on $P_d(G)$, and $\Psi_m(p(F_m))=p(\bar F)$ for every noncommutative
polynomial $p$, in particular for $I$, $q_\eps\circ I$ and $w_\eps(\frac12(\cdot+1))$. Conversely, every normalised cycle $\bar F$ on
$P_d(G)$ is $\Psi_m(F_m)$ for a normalised cycle $F_m$ on $P_d(G_m)$. If $p(F_m)$ has propagation at most $R_m/2$, then
$\Phi_m(p(\bar F))=p(F_m)$.
\end{lemma}

\begin{proof}
Write points of $Z_d$ as measures and $\supp$ for their supports. For $H=G_m$ or $G$, $x\in H$, and $a\ge d$, put
\[
 \mathcal V_a^H(x)=\{\nu\in P_d(H):\supp\nu\subseteq B_H(x,a)\},
\]
where $B_H(x,a)$ is a word-metric ball. These are compact support regions, and
$\mathcal V_a^H(x)\subseteq\operatorname{int}\mathcal V_{a+d}^H(x)$: if $\nu(v)>0$, then every sufficiently close measure $\mu$ has $\mu(v)>0$, so its support lies in $B_H(v,d)\subseteq B_H(x,a+d)$.

We will use the support of the normalising cut-offs to localise products. Each summand of a normalised operator has both support variables in
$gW_d=\mathcal V_{2d}^H(g)$. If this region meets $\mathcal V_a^H(x)$, a vertex in the support of a measure in their intersection gives $d_H(x,g)\le a+2d$.
Thus $gW_d\subseteq\mathcal V_{a+4d}^H(x)$. Starting with a function supported in $\mathcal V_{2d}^H(x)$, one factor of the operator therefore reaches only $\mathcal V_{6d}^H(x)$, and two factors reach only $\mathcal V_{10d}^H(x)$. The same bounds hold for adjoints.

Since $16d<R_m$, push-forward identifies $\mathcal V_{12d}^{G_m}(x)$ with $\mathcal V_{12d}^G(\pi_mx)$. Its interior is an open sheet containing the compact region $\mathcal V_{10d}^{G_m}(x)$. Its $K_m$-translates are disjoint: an intersection would give $d_{G_m}(x,kx)\le24d<2R_m$ for some nonidentity $k\in K_m$, contradicting the kernel-length bound for $x^{-1}kx\in K_m$.

We first express the pushed-down operator in normalised form. Let $\chi$ be the projection onto the sheet over $W_d$ through $e$, and $\bar T$ the operator on
$\ell^2(Z_d(G))\otimes H_0$ corresponding to $\chi T\chi$; then $\|\bar T\|\le1$. Since $\sqrt{\phi_m}=\sqrt{\phi_m}\chi$,
$F_m=\sum_gg\sqrt{\phi_m}\chi T\chi\sqrt{\phi_m}g^{-1}$. Let $z\in D_d(G_m)$, so $\supp z\subseteq B(e,d)$. The term with index $g$
has nonzero entries at $(z,z')$ only if $z,z'\in gW_d(G_m)$, that is, $\supp z,\supp z'\subseteq B(g,2d)$. This forces
$|g|\le3d$ and $\supp z'\subseteq B(e,5d)$. All these supports lie in $B(e,5d)$, where the push-forward is bijective and
$\pi_m$ is injective on $\{g:|g|\le3d\}$. So the entries of $F_m$ at $(z,z')$ equal those of
$\sum_{\bar g}\bar g\sqrt\phi\,\bar T\sqrt\phi\,\bar g^{-1}$ at the images. By equivariance this holds at all pairs, and
\[
\Psi_m(F_m)=\sum_{\bar g\in G}\bar g\sqrt\phi\,\bar T\sqrt\phi\,\bar g^{-1},
\]
the coefficients being relabelled because $F_m$ has propagation $6d\le R_m$.

We next verify the cycle conditions. The map $\Psi_m$ is a $*$-homomorphism on all of $\C G_m\odot\B(L_d)$, since $\pi_m$ is a homomorphism. So
$\Psi_m(p(F_m))=p(\bar F)$ for every polynomial $p$. In particular $\bar F^*\bar F-1$ and $\bar F\bar F^*-1$ have the coefficients of
$F_m^*F_m-1$ and $F_m F_m^*-1$, and these are compact. For commutators let $f\in C_c(P_d(G))$, and write
$f=\sum_{\bar g}f\cdot\bar g\phi$, a finite sum. Fix $\bar g$ and a lift $g$. The operator $[\bar F,f\cdot\bar g\phi]$ has nonzero
entries only at pairs of points with supports in $B(\bar g,10d)$: one point lies in $\bar gW_d$, so it is in some $hD_d$ with
$d(h,\bar g)\le3d$, and the other is within propagation $6d$ of it. By the normal form above, which is local, it is the conjugate of $[F_m,\tilde f]$ by the
identification of the measures supported in $B(g,10d)$ with those supported in $B(\bar g,10d)$, where $\tilde f$ is the lift of
$f\cdot\bar g\phi$ to this region. That commutator is compact because $F_m$ is a cycle. By density, $[\bar F,f]$ is compact
for all $f\in C_0(P_d(G))$.

Conversely, given $\bar F=\sum_{\bar g}\bar g\sqrt\phi\,\bar T\sqrt\phi\,\bar g^{-1}$, let $T$ be the lift of $\chi\bar T\chi$ to the sheet and
$F_m$ the corresponding normalised operator. Then $\Psi_m(F_m)=\bar F$ by the normal form computation. The
coefficients of $F_m^*F_m-1$ and $F_mF_m^*-1$ are those of $\bar F^*\bar F-1$ and $\bar F\bar F^*-1$, because these operators have
propagation at most $12d\le R_m$ and $\Psi_m$ relabels coefficients there. The commutator condition follows by the
same local identification, run in the other direction. Finally, if $p(F_m)$ has propagation at most $R_m/2$, then
$\Phi_m$ inverts $\Psi_m$ on it, so $\Phi_m(p(\bar F))=p(F_m)$.
\end{proof}

\subsection{Transport of classes}
\begin{lemma}\label{lem:transport-classes}
Let $16d<R_m$ and let $B$ be a separable $G$-$C^*$-algebra, on which $G_m$ acts through $\pi_m$. There is a homomorphism
\[
\Psi_{m*}\colon KK^{G_m}_*(C_0(P_d(G_m)),B)\to KK^G_*(C_0(P_d(G)),B)
\]
with $\Psi_{m*}[F_m]=[\Psi_m(F_m)]$ for every normalised cycle $F_m$, odd or even. For $d\le d'$ with $16d'<R_m$ it commutes with the
maps induced by $P_d\subseteq P_{d'}$.
\end{lemma}

Conceptually, $\Psi_{m*}$ is the free-quotient map for the free and proper action of $K_m=\ker\pi_m$ on
$\tilde X=P_d(G_m)$, whose quotient $G$-space is $X=P_d(G)$ (Lemma~\ref{lem:rips-covering}). It is the $K$-homology analogue of Segal's
isomorphism $K_{G/N}(X/N)\cong K_G(X)$ for compact $G$ and a normal subgroup $N$ acting freely \cite[Proposition~2.1]{Seg68}. For trivial
$G$ the corresponding lifting map appears in \cite[Lemma~4.4]{WY12}. The proof, in Appendix~\ref{app:bc-quotient}, uses Meyer's theory of
exactly equivariant Kasparov triples \cite{Mey00}, which works for all locally compact groups.

\begin{remark}\label{rem:transport-coeff}
The cycle-transport lemma also holds on the standard Hilbert $B$-modules when the action factors through the quotient. The cut-off compressions are adjointable module contractions, and the local matrix comparisons respect the action on $B$. Compact coefficients and compact local commutators are preserved. In particular both polynomial indices lift and descend exactly whenever their propagation is below the local multiplicativity radius. The class-transport proof in Appendix~\ref{app:bc-quotient} already includes these coefficients and arbitrary homotopies; transport therefore defines a map on $KK$-classes.
\end{remark}

\subsection{Quantitative surjectivity}\label{sec:qs}
We now formulate the finite data needed for the pushdown argument. We do this
with actual paths; it also specifies the meaning of size-dependent controls
used in the introduction.

\begin{definition}\label{def:finite-input}
An \emph{input of propagation $r$, size at most $n$, and tolerance $\eps$}
is one of the following:
\begin{enumerate}[label=(\roman*)]
\item a self-adjoint $p\in M_k(\C G_r)$, $k\le n$, with
$\|p^2-p\|\le\eps$, together with a scalar projection $e_l$ of rank
$0\le l\le k$; it represents $[p]-[e_l]$;
\item $u\in M_k(\C G_r)$, $k\le n$, with both unitary defects at most
$\eps$; it represents $[u]$.
\end{enumerate}
The \emph{size of a finite propagation homotopy} in a stabilised Roe algebra
is the dimension of a common finite-dimensional fibre corner containing its
non-scalar parts, including all matrix amplifications and stabilising blocks.
On this corner it is a path in $M_N(\C G)$; its scalar background is recorded
separately. Thus $N$ records the full active dimension of the stabilised homotopy.
\end{definition}

\begin{definition}\label{def:qs}
The assembly map of $G$ is \emph{quantitatively surjective with size-dependent
controls} if there is $\Lambda\ge1$ such that the following holds. Given
positive integers $d,r,n$ and $0<\eps<1/(128\Lambda)$, there are integers
$d'\ge d$, $t\ge r$, and $N\ge n$, independent of the input, such that every
input from Definition~\ref{def:finite-input} has a normalised cycle $F$ on
$P_{d'}(G)$ and a finite compression $c$ of a polynomial index $J^*_{\xi}(F)$,
where $0<\xi\le\eps$ and the coefficient $\ell^1$ compression error is at most
$\eps$, satisfying the following conditions.

\begin{enumerate}[label=(\roman*)]
\item The polynomial index has propagation at most $t$.
\item There is a path of
propagation at most $t$, size at most $N$, and defects at most $\Lambda\eps$
from the input to $c$, after scalar stabilisation. In degree zero, if $e_Q$
is the active scalar projection of the compressed index, the endpoints are
$p\oplus e_Q$ and $c\oplus e_l$, with the same extra identity blocks and
the necessary zero padding. In degree one the endpoints are $u$ and $c$,
stabilised by identity blocks.
\end{enumerate}
The output bounds $d',t,N$ may depend on
$G,d,r,n,\eps$, while $\Lambda$ is fixed for $G$.
\end{definition}

The input corner is fixed in advance. Any other corner of the same finite
dimension can be moved to it by a unitary path in the fibre; this costs no
group propagation and uses only a finite-dimensional enlargement. Thus the
definition describes finite-size representatives in the Roe-algebra picture
as well. Restoring the compressed endpoint uses a straight segment with
defects at most $4\eps$, after decreasing the allowed compression error if
necessary. Consequently the paths give the expected equalities in
quantitative $K$-theory with enlarged output parameters. These explicit paths provide the parameters used in the covering argument.
The integer $d$ serves only as a lower bound for the Rips scale $d'$.

\begin{proposition}\label{prop:qs-ordinary}
For a finitely generated group $G$, ordinary assembly surjectivity in either
degree is equivalent to size-dependent quantitative surjectivity in that
degree as in Definition~\ref{def:qs}.
\end{proposition}
\begin{proof}
The reverse implication follows by choosing a finite group-ring
representative of any ordinary $K$-class and then making $\eps$ small enough
that the path and endpoint segment have defects less than $1/4$.

Conversely, suppose $\mu_G$ is surjective in the given degree, put
$\Lambda=20$, and fix $d,r,n$ and $0<\eps<1/(128\Lambda)$.
Lemma~\ref{lem:app-uniform}, with $\Gamma=G$, gives $d',t,N$ independent of
the input, together with the cycles, compressions and paths of
Definition~\ref{def:qs}, with $\xi=\eps/100$, compression error $\eps/100$,
and defects at most $20\eps=\Lambda\eps$.
\end{proof}

This equivalence lets us express ordinary surjectivity through simultaneous
bounds on propagation and total matrix size. We can then compare these
bounds with the covering radii.

\begin{definition}[Surjectivity with asymptotic controls]\label{def:asym-surj}
Let $\{\pi_m\colon G_m\to G\}$ be an asymptotically faithful family and let
$(K,h)$ be a pushdown control on self-adjoint finite matrix defects
(Definition~\ref{def:pushdown}). Put $\eps_l=2^{-l-10}$. The approximants have
\emph{surjectivity with asymptotic controls} if there is $A\ge1$ such that,
for every fixed positive integer tuple $(d,r,n,l)$ and arbitrarily large $m$,
there are integers $k_m\ge d$, $t_m\ge r$, and $N_m\ge n$ with
\begin{equation}\label{eq:asym-witness-control}
 16k_m<R_m,\qquad 2t_m\le h(m,N_m),
\end{equation}
having this property: every input over $G_m$ of propagation $r$, size at most
$n$, and tolerance $\eps_l$ has a normalised cycle on $P_{k_m}(G_m)$ and a
path as in Definition~\ref{def:qs}, of propagation at most $t_m$, size at
most $N_m$, and defects at most $A\eps_l$. The index is
$J^*_{\eps_l/100}(F_m)$, its propagation is at most $t_m$, its coefficient
$\ell^1$ compression error is less than $\eps_l/100$, and all path defects
belong to the class on which the pushdown bound holds. All these bounds
are independent of the input.
\end{definition}

\begin{theorem}\label{thm:qs-framework}
Let $\{G_m\to G\}$ have operator norm localisation with asymptotic controls. If the approximants
have surjectivity with asymptotic controls relative to $(K,h)$, then $\mu_G$
is quantitatively surjective with size-dependent controls in both degrees.
\end{theorem}
\begin{proof}
Fix $d,r,n$ and a sufficiently small input tolerance $\eps$. Choose $l$ so
that $2\eps\le\eps_l<4\eps$, and take one sufficiently large $m$ supplied
by Definition~\ref{def:asym-surj}, also with $2r\le S_m$ and $c_m^{-1}<2$.
These choices do not depend on the input.

We first lift the input using (ii) of Figure~\ref{fig:directions}. Local multiplicativity gives
\[
 \Phi_m(p)^2-\Phi_m(p)=\Phi_m(p^2-p),
 \qquad
 \Phi_m(u)^*\Phi_m(u)-1=\Phi_m(u^*u-1),
\]
and the analogous second unitary identity. Upward norm control, Theorem~\ref{thm:asym-isometric}, bounds
these defects by $2\eps\le\eps_l$.

We now apply the surjectivity hypothesis to obtain $F_m$, its compressed polynomial
index, and the stated actual path.

To push this path down as in (iii), apply the quotient homomorphism coefficientwise to the
path. Its finite support and finite matrix size make this a continuous
path. Algebraic multiplicativity identifies its defects with the
pushdowns of the original self-adjoint defects. They have propagation at
most $2t_m$, so their norms are at most
$KA\eps_l<4KA\eps$. Its initial endpoint is the original input, with
the prescribed scalar stabilisations.

We finally transport the cycle using (i). By Lemma~\ref{lem:transport}, $\bar F_m=\Psi_m(F_m)$ is a normalised cycle
on $P_{k_m}(G)$ and its polynomial index is the pushdown of that of $F_m$.
Coefficientwise quotienting does not increase $\|\cdot\|_1$ and commutes
with the fibre compression. The other endpoint is therefore a compression
of $J^*_{\eps_l/100}(\bar F_m)$ with error less than
$\eps_l/100<\eps$. This is exactly the data in Definition~\ref{def:qs}.
Taking, for example, $\Lambda=\max\{4KA,2^{12}\}$ ensures the required
small-tolerance range as well. The argument treats both degrees.
\end{proof}

% The downward norm estimate controls the defects of the finite homotopy,
% while the inequality $16k_m<R_m$ allows us to transport the topological
% cycle by the local transport lemma.

\section{Lacunary hyperbolic groups}\label{sec:lacunary}
We now apply the preceding framework to lacunary hyperbolic groups. The
upward norm estimate holds for every such group. Injectivity follows
from it and hyperbolic assembly, through the weak equivalence of
Section~\ref{sec:weak-equivalence}; surjectivity will additionally use
rapid decay and fast-growing injectivity radii.

\subsection{Operator norm localisation for the approximants}
Throughout this subsection, $S$ denotes the finite marking and we put
$D_S=\max\{2,2|S|\}$.

\begin{theorem}\label{thm:onl-lacunary}
Let $G=\varinjlim G_m$ be lacunary hyperbolic. Then the approximants have ONL with asymptotic controls in the sense of Definition~\ref{def:asym-onl},
with constant $c=1/(2D_S)$ and controls $g_m(R)\le432\delta_mR+4\delta_m+26$.
\end{theorem}

We collect some facts first. Hyperbolic groups have finite asymptotic dimension; this was shown by
Roe \cite{Roe05}, confirming an idea of Gromov \cite{Gro93}. We need the cover from his proof with explicit
constants, and we run his argument on the vertex set. Let $\Gamma$ be $\delta$-hyperbolic in the sense of
Definition~\ref{def:four-point}, with $\delta\in\N$, $\delta\ge1$. Write $|x|=d(x,e)$ and
$(x|y)=\frac12(|x|+|y|-d(x,y))$. Then $(x|z)\ge\min\{(x|y),(y|z)\}-\delta$ for all $x,y,z\in\Gamma$, which is Roe's
hypothesis. Put $s=8\delta$ and $M=|B(s)|$, so $M\le2D_S^s$ when $D_S\ge2$.

\begin{lemma}[Roe's cover]\label{lem:roe-cover}
Let $d\ge2$ be an even integer. Then $\Gamma$ is covered by $2M$ families of subsets such that, within each
family, distinct members are at distance at least $d/2$, and every member has diameter at most $4d+4\delta$.
\end{lemma}

\begin{proof}
We construct the cover on each annulus as in \cite[proof of Theorem~2.1]{Roe05}. Let $A_k=\{x\mid kd\le|x|<(k+1)d\}$. Let
$\{x_i\}$ be a maximal subset of the sphere $\{|x|=kd\}$ with pairwise distances greater than $d$, and put
$U_i=\{x\in A_k\mid(x|x_i)\ge(k-\frac12)d-\delta\}$. For $x\in A_k$, let $x'$ be the vertex at distance $kd$
from $e$ on a geodesic from $e$ to $x$. Then $(x|x')=kd$, and by maximality $d(x',x_i)\le d$ for some $i$.
So $(x'|x_i)\ge(k-\frac12)d$ and $(x|x_i)\ge(k-\frac12)d-\delta$, and the $U_i$ cover $A_k$. If $x\in U_i$ then $d(x,x_i)=|x|+|x_i|-2(x|x_i)\le2d+2\delta$, so $\diam U_i\le4d+4\delta$.

We next colour the sets in each annulus. For $k=0$ there is a single set, so let $k\ge1$. Call $U_i\ne U_{i'}$ \emph{adjacent} if $d(U_i,U_{i'})<d/2$.
In that case pick $y\in U_i$ and $z\in U_{i'}$ with $d(y,z)<d/2$. Since $|y|,|z|\ge kd$, we have $(y|z)>(k-\frac14)d$, and
$(x_i|y),(x_{i'}|z)\ge(k-\frac12)d-\delta$, so $(x_i|x_{i'})\ge(k-\frac12)d-3\delta$. Let $x''$ be the vertex at distance $(k-\frac12)d$
from $e$ on a geodesic from $e$ to $x_i$; this is an integer because $d$ is even. Then $(x_i|x'')=(k-\frac12)d$, so
$(x_{i'}|x'')\ge(k-\frac12)d-4\delta$ and $d(x_{i'},x'')\le\frac d2+8\delta$. Every vertex $z'$ with $d(z',x'')\le\frac d2+s$ is within $d/2$ of a
vertex of $B(x'',s)$, on a geodesic from $x''$ to $z'$. A closed ball of radius $d/2$ contains at most one of the $x_{i'}$,
whose pairwise distances exceed $d$. Hence $U_i$ has fewer than $M$ neighbours. Colouring the countably many $U_i$
greedily, in some enumeration, with the colours $\{1,\dots,M\}$ gives adjacent sets distinct colours. So within $A_k$, sets
of the same colour are at distance at least $d/2$.

To combine the covers of the annuli, use the colours $\{1,\dots,M\}$ on the annuli with $k$ even, and a disjoint set of $M$ colours
on those with $k$ odd. Annuli $A_k$ and $A_{k+2}$ are at distance at least $d$. So each colour class consists
of sets at pairwise distance at least $d/2$.
\end{proof}

So $\Gamma$ has the metric sparsification property of Chen--Tessera--Wang--Yu \cite[Definition~3.1, Remark~3.2]{CTWY08}
with constant $1/(2M)$ and function $m\mapsto8m+4\delta$. Indeed, given $m$ and a finite measure $\mu$, apply
Lemma~\ref{lem:roe-cover} with $d=2m$; one of the $2M$ families carries at least $\mu(\Gamma)/(2M)$.

\begin{lemma}[{\cite[Proposition~4.1 and its proof]{CTWY08}}]\label{lem:ms-onl}
If $X$ has the metric sparsification property with constant $c$ and function $f$, then for every $\eps>0$,
$X$ has ONL with constant $\sqrt c(1-\eps)$ and control $k\mapsto f(3m)+2m$, where $m=2k+1$.\qed
\end{lemma}

In the proof of \cite[Proposition~4.1]{CTWY08} the localising vector is supported in an $m$-neighbourhood of a
piece of diameter at most $f(3m)$, for any $m>2k$. This gives the control stated.

\begin{proof}[Proof of Theorem~\ref{thm:onl-lacunary}]
We may assume $\delta_m\in\N$ and $\delta_m\ge1$. By Lemma~\ref{lem:roe-cover}, the remark after it and Lemma~\ref{lem:ms-onl},
$G_m$ has ONL with constant $c_0=1/(2\sqrt M)$ and control
\[
f(k)=8(6k+3)+4\delta_m+2(2k+1)=52k+4\delta_m+26,
\]
where $M\le2D_S^s$ and $s=8\delta_m$. Let $n=4\delta_m\ge4$. Then $c_0^{-1/n}\le(2\sqrt2)^{1/n}D_S^{s/(2n)}\le2D_S$, so
Lemma~\ref{lem:onl-constant} gives ONL with constant $1/(2D_S)$ and control
\[
g_m(R)=4(n-1)R+f(2nR)=(108n-4)R+4\delta_m+26\le432\delta_mR+4\delta_m+26,
\]
and, on self-adjoint operators, control $2(n-1)R+f(nR)\le216\delta_mR+4\delta_m+26$. Any $R\le(R_m-4\delta_m-26)/(432\delta_m+1)$ satisfies
$R+g_m(R)\le R_m$, and this bound tends to infinity because $\delta_m=o(R_m)$.
\end{proof}

We record the following results on hyperbolic groups from the literature for convenience.
\begin{theorem}\label{thm:hyperbolic-input}
If $\Gamma$ is $\delta$-hyperbolic in the convention of \eqref{eq:fourpoint}, then
\begin{enumerate}[label=(\roman*)]
\item $P_d(\Gamma)$ is a model for $\underline E\Gamma$ whenever $d\ge d_\delta:=32\delta+20$ \cite[Theorem 1]{MS02};
\item $\mu_{\Gamma,B}$ is an isomorphism for every separable $\Gamma$-$C^*$-algebra $B$ \cite{Laf12}.
\end{enumerate}
\end{theorem}

\subsection{Weak equivalence of the asymptotic lift}\label{sec:weak-equivalence}
This subsection makes arrow (ii) of Figure~\ref{fig:directions} precise. Let
$\{\pi_m\colon G_m\to G\}$ be an asymptotically faithful covering family of marked groups,
with isometry radii $R_m$ and kernels $K_m$, and let $(B,\beta)$ be a separable
$G$-$C^*$-algebra, on which $G_m$ acts through $\pi_m$. For
$a=\sum_{|g|\le r}b_gu_g\in C_c(G,B)$ and $2r<R_m$ put
\[
 \Phi_m(a)=\sum_{|g|\le r}b_gu_{\tilde g}\in C_c(G_m,B),
\]
where $\tilde g$ is the unique lift of $g$ in $B_{G_m}(e,r)$. We use the same notation
after tensoring with $M_n(\C)$ or with $\K(\mathcal H')$ for a separable Hilbert space
$\mathcal H'$, such as $\K(L_d)$. If $|g|+|h|<R_m$ then $\tilde g\tilde h$ is the lift of
$gh$, since $\pi_m$ is injective on $B_{G_m}(e,R_m)$. Hence $\Phi_m$ preserves adjoints,
units, and every product of elements whose propagations sum to less than $R_m$.
Let
\[
 \mathcal Q_B=\prod_mB\rtimes_rG_m\Big/\bigoplus_mB\rtimes_rG_m,\qquad
 \Phi_\infty(a)=[(\Phi_m(a))_m],
\]
with arbitrary coordinates for the finitely many $m$ with $2r\ge R_m$. This is a
$*$-homomorphism $C_c(G,B)\to\mathcal Q_B$, and
$\|\Phi_\infty(a)\|=\limsup_m\|\Phi_m(a)\|$. For $B=\C$ it is the representation
$\lambda_\infty$ of Section~\ref{sec:onl}.

\begin{theorem}[Weak equivalence of the asymptotic lift]\label{thm:weak-equivalence}
In this setting the following hold, for every matrix size and every Hilbert fibre.
\begin{enumerate}[label=(\arabic*)]
\item For every $a$, $\|a\|_{B\rtimes_rG}\le\liminf_m\|\Phi_m(a)\|$. Hence the identity of
$C_c(G,B)$ extends to a surjection from the closure of $\Phi_\infty(C_c(G,B))$ onto
$B\rtimes_rG$. For $B=\C$: $\lambda_\infty$ weakly contains $\lambda$.
\item Suppose that the approximants have ONL with asymptotic controls, and let
$c_m\to1$ and $S_m\to\infty$ be as in Lemma~\ref{lem:cm-to-one}. Then
\[
 c_m\|\Phi_m(a)\|\le\|a\|\qquad(\prop(a)\le S_m).
\]
Hence $\|\Phi_m(a)\|\to\|a\|$, $\lambda_\infty$ and $\lambda$ are weakly equivalent, and
$\Phi_\infty$ extends to an injective $*$-homomorphism $B\rtimes_rG\to\mathcal Q_B$.
\item Under the hypothesis of \textup{(2)}, let $A=\K(\mathcal H')\otimes(B\rtimes_rG)$ and
$A_m=\K(\mathcal H')\otimes(B\rtimes_rG_m)$, filtered by propagation. Let $r\ge1$ and
$0<\eps<1/48$, and let $m$ satisfy $4r\le S_m$ and $c_m^{-1}\le2$. Then $\Phi_m$ maps
$(r,\eps)$-quasi-projections and quasi-unitaries to $(r,2\eps)$-quasi-projections and
quasi-unitaries, paths of $(r,\eps)$-quasi-projections to paths of
$(r,2\eps)$-quasi-projections, and paths of $(2r,3\eps)$-quasi-unitaries to paths of
$(2r,6\eps)$-quasi-unitaries. It therefore induces homomorphisms
\[
 \Phi_{m*}\colon K_*^{r,\eps}(A)\longrightarrow K_*^{r,2\eps}(A_m),\qquad *=0,1,
\]
compatible with the structure maps.
\end{enumerate}
By Theorem~\ref{thm:onl-lacunary}, \textup{(2)} and \textup{(3)} hold for every lacunary
presentation, with every separable coefficient algebra.
\end{theorem}

\begin{proof}
Represent $B$ faithfully and nondegenerately on a Hilbert space $H$, and use the regular
representation of $B\rtimes_rG$ on $\ell^2(G)\otimes H$ in its right-regular form,
$(bu_g\zeta)(h)=\beta_h(b)\zeta(hg)$; it is unitarily equivalent to the usual one by
$\zeta\mapsto\zeta\circ\mathrm{inv}$. An element $a=\sum b_gu_g$ supported in $B(e,r)$
then has matrix entries $a_{h,hg}=\beta_h(b_g)$, so it is an operator of propagation at
most $r$ on $\ell^2(G)\otimes H$ in the sense of Section~\ref{sec:onl}. The same holds over
$G_m$. The entries of $\Phi_m(a)$ are $\Phi_m(a)_{\tilde h,\tilde h\tilde g}=\beta_{\pi_m(\tilde h)}(b_g)$,
because the action factors through $\pi_m$. These are the entries of $a$ at the image pair.
So $\Phi_m(a)$ is the local lift of Definition~\ref{def:lifting}, with deck group $K_m$, and
the same holds after tensoring with $M_n(\C)$ or $\K(\mathcal H')$.

(1) The first half of the proof of Theorem~\ref{thm:asym-isometric} uses only that
$\pi_m$ is isometric on large balls: given $\eta>0$, a finitely supported $\xi$ with
$\|a\xi\|\ge(1-\eta)\|a\|\|\xi\|$ lifts to one sheet, isometrically and compatibly with
$a$, as soon as $R_m$ exceeds twice the radius of a ball containing $\supp\xi$ and
$\supp a\xi$. So $\|\Phi_m(a)\|\ge(1-\eta)\|a\|$ for all large $m$. In particular
$\|a\|\le\|\Phi_\infty(a)\|$, which is the stated surjection.

(2) The second half of the proof of Theorem~\ref{thm:asym-isometric} applies to the
local lift of any operator of propagation at most $S_m$, with any Hilbert fibre.
Together with (1) and $c_m\to1$ this gives $\|\Phi_m(a)\|\to\|a\|$. So $\Phi_\infty$ is
isometric on $C_c(G,B)$ and extends to an isometric, hence injective, $*$-homomorphism.

(3) Extend $\Phi_m$ to unitisations by $1\mapsto1$. If $p$ is an
$(r,\eps)$-quasi-projection then
$\Phi_m(p)^2-\Phi_m(p)=\Phi_m(p^2-p)$, of propagation at most $2r\le S_m$, so its norm is at
most $c_m^{-1}\eps\le2\eps$. The two defects of a quasi-unitary are treated in the same
way. Along a path the estimate holds pointwise, and $\Phi_m$ is bounded on each
propagation subspace below $S_m$, by (2), so continuity is preserved. For paths of
$(2r,3\eps)$-quasi-unitaries the defects have propagation at most $4r\le S_m$, giving
$(2r,6\eps)$-quasi-unitaries; $6\eps<1/8$. Identity and zero blocks are preserved, so the
stabilisations and equivalences of Definition~\ref{def:qk} are respected. Additivity and
compatibility with the structure maps are clear.
\end{proof}

We use (3) to transport a vanishing class to individual approximants. Such a
vanishing is witnessed by a path of fixed propagation, which lifts to paths over
$A_m$ for all large $m$, with the stated defect bounds. This gives the finite
propagation input to hyperbolic injectivity in the following proof.

\subsection{Injectivity, with coefficients}
\begin{theorem}\label{thm:main-inj}
For every lacunary hyperbolic group $G$ and every separable $G$-$C^*$-algebra $B$, the map $\mu_{G,B}$ is injective in both degrees.
\end{theorem}

\begin{proof}[Proof of Theorem~\ref{thm:main-inj}]
Let $x\in K_*^{\mathrm{top}}(G;B)$ have zero assembly. We follow arrows (ii) and (i) of
Figure~\ref{fig:directions}.

Choose $d$ and a normalised cycle $F$ on
$P_d(G)$ representing $x$ (Lemma~\ref{lem:ample-rep}), and put $\eps=1/1000$. The
polynomial index $J^*_\eps(F)$ has some propagation $r_0$, and its ordinary $K$-class is
$\mu_d[F]=\mu_{G,B}(x)=0$. Lemma~\ref{lem:app-path}, applied in $\K(L_d)\otimes(B\rtimes_rG)$
with a scalar second endpoint, gives an actual path, of one finite propagation $r\ge r_0$
and with defects at most $\eps'=20\eps$, from the stabilised $J^*_\eps(F)$ to a scalar
representative. In even degree the scalar subtraction is
handled by adding the same scalar blocks; in odd degree the endpoint is the identity.
In the notation of
Theorem~\partref{thm:weak-equivalence}{3}, with $\mathcal H'=L_d$, the path witnesses
$[J^*_\eps(F)]=0$ in $K_*^{r,\eps'}(A)$.

We now lift this path. By Theorem~\ref{thm:onl-lacunary} the approximants have
ONL with asymptotic controls. Put
$d_m'=\max(d,32\delta_m+20)$, and take $m$ so large that $4r\le S_m$,
$c_m^{-1}\le2$ and $16d_m'<R_m$; the last condition holds for large $m$ because
$\delta_m=o(R_m)$. Theorem~\partref{thm:weak-equivalence}{3}, applied pointwise, maps the path to a
path of defects at most $2\eps'<1/4$ over $A_m$, from $\Phi_m(J^*_\eps(F))$ to the same scalar endpoint. By
Lemma~\ref{lem:transport} and Remark~\ref{rem:transport-coeff}, $F=\Psi_m(F_m)$ for a
normalised cycle $F_m$ on $P_d(G_m)$, and $\Phi_m(J^*_\eps(F))=J^*_\eps(F_m)$ because
$r_0\le R_m/2$. Hence $\mu_{G_m,B}[F_m]=0$ in ordinary $K$-theory.

Hyperbolic injectivity and the Rips model in
Theorem~\ref{thm:hyperbolic-input} now show that the image of $[F_m]$ in
$KK_*^{G_m}(C_0(P_{d_m'}(G_m)),B)$ is zero.

To return to $G$, apply the class-transport map of
Lemma~\ref{lem:transport-classes} at scale $d_m'$. Its naturality in the Rips scale
identifies the result with the image of $[\Psi_m(F_m)]=[F]$ on $P_{d_m'}(G)$. That image
is therefore zero, and $x=0$ in the direct limit.

The same argument treats both degrees and all the stated coefficients.
\end{proof}

\subsection{Rapid decay and pushdown}\label{sec:assembly}
\begin{definition}\label{def:rd}
A group $G$ has \emph{ball-form rapid decay} with constants $C\ge1$ and $s\ge0$ if
\begin{equation}\label{eq:rd-ball}
 \|\lambda_G(a)\|\le C(1+r)^s\|a\|_2
 \quad\text{for }a\in\C G\text{ supported in }B_G(e,r).
\end{equation}
This is the form of property \RD{} convenient for our purposes. It implies a Sobolev RD estimate with every exponent strictly greater than $s$. Indeed, decompose a function into dyadic annuli, apply the ball estimate to each annular part, and use Cauchy--Schwarz with the summable weights $2^{-j(s'-s)}$ for $s'>s$. The converse implication follows immediately from support in a ball.
\end{definition}

For $b=\sum_gb_g\lambda_g\in M_n(\C H)$ write
\[
 \|b\|_{2,\mathrm{HS}}^2=\sum_g\|b_g\|_{\mathrm{HS}}^2.
\]
If $G$ satisfies \eqref{eq:rd-ball}, then
\begin{equation}\label{eq:rd-matrix}
 \|b\|\le C(1+r)^s\|b\|_{2,\mathrm{HS}}
 \quad(b\in M_n(\C G_r)).
\end{equation}
Indeed, pointwise domination on $\ell^2(G,\C^n)$ gives $\|b\|\le\|\lambda_G(g\mapsto\|b_g\|)\|$, and $\|b_g\|\le\|b_g\|_{\mathrm{HS}}$. This estimate is uniform in $n$ when the Hilbert--Schmidt coefficient norm is used.

\begin{proposition}[RD controls finite matrix pushdown]\label{prop:rd-pushdown}
Let $\pi\colon H\to G$ be a marked quotient with isometry radius $R$, and suppose $G$ satisfies \eqref{eq:rd-ball}. If $b=b^*\in M_n(\C H_r)$ and $kr\le R$, then
\begin{equation}\label{eq:power-pushdown}
 \|\pi(b)\|\le C^{1/k}(1+kr)^{s/k}n^{1/(2k)}\|b\|.
\end{equation}
Consequently, with $K=C e^{2s+1/2}$, one has
\begin{equation}\label{eq:log-pushdown}
 \|\pi(b)\|\le K\|b\|
 \quad\text{if }r\ge1\text{ and }r\lceil\log(enr)\rceil\le R.
\end{equation}
\end{proposition}

\begin{proof}
The quotient is an algebraic $*$-homomorphism, so $\pi(b)^k=\pi(b^k)$. On the support of $b^k$, contained in $B_H(e,kr)$, it is injective. Hence
\[
 \|\pi(b^k)\|_{2,\mathrm{HS}}=\|b^k\|_{2,\mathrm{HS}}
 \le\sqrt n\,\|b^k\|,
\]
where the last inequality uses the canonical trace $\tau\otimes\Tr$ on $M_n(C_r^*H)$. Both $b$ and $\pi(b)$ are self-adjoint. Applying \eqref{eq:rd-matrix} to the power and taking $k$th roots proves \eqref{eq:power-pushdown}.

For \eqref{eq:log-pushdown}, put $L=\log(enr)\ge1$ and $k=\lceil L\rceil$. Then $C^{1/k}\le C$, $n^{1/(2k)}\le e^{1/2}$, and
\[
 \log(1+kr)\le\log2+\log(L+1)+\log r\le2L.
\]
This gives $(1+kr)^{s/k}\le e^{2s}$.
\end{proof}

We apply this proposition to $p^2-p$ in even degree and to $u^*u-1$ and $uu^*-1$ in odd degree. Scalar blocks are unchanged by pushdown, and the active finite corner has exactly the size counted in Lemma~\ref{lem:finite-witness}.

In the sense of Definition~\ref{def:pushdown}, this gives a pushdown control
\begin{equation}\label{eq:rd-h}
 h(m,n)=\max\bigl(\{0\}\cup\{r\in\N:r\le R_m/2,
          \ r\lceil\log(enr)\rceil\le R_m\}\bigr),
 \qquad K=Ce^{2s+1/2}.
\end{equation}
For each fixed $n$, this radius tends to infinity.

\subsection{Finite-size solutions}\label{sec:witnesses}
The topological preimage in a hyperbolic group can be placed on one fixed Rips complex. We need the following uniform finite-size form of ordinary surjectivity.

\begin{lemma}[Finite witnesses for fixed inputs]\label{lem:finite-witness}
Fix a marked $\delta$-hyperbolic group $\Gamma$, positive integers $d,r,n$, and $0<\eps<1/128$. Set $d'=\max(d,d_\delta)$. There are integers $t,N$, with
\begin{equation}\label{eq:witness-floor}
 t\ge\max\{1,r,16d'\},\qquad N\ge n,
\end{equation}
that work for every input of either of the following kinds:
\begin{enumerate}[label=(\roman*)]
\item $p=p^*\in M_k(\C\Gamma_r)$, $k\le n$, $\|p^2-p\|\le\eps$, together with a scalar projection $e_l$ of rank $0\le l\le k$;
\item $u\in M_k(\C\Gamma_r)$, $k\le n$, with both unitary defects at most $\eps$.
\end{enumerate}
For each input there is a normalised cycle $F$ of the corresponding parity on $P_{d'}(\Gamma)$ and a finite-rank fibre projection $Q$ such that the following hold.

The polynomial index $J^*_{\eps/100}(F)$ has propagation at most $t$ and an $\eps/100$-compression in coefficient $\ell^1$-norm to $Q$. Write $c$ for its active finite matrix and $e_Q$ for its active scalar projection in even degree. There is a path of propagation at most $t$ and size at most $N$, with defects at most $20\eps$, joining the input to this compressed index after scalar stabilisation. Explicitly, in even degree the path joins stabilisations of $p\oplus e_Q$ and $c\oplus e_l$, with the same additional identity blocks and the necessary zero padding. In odd degree it joins stabilisations of $u$ and $c$ by identity blocks.
\end{lemma}

\begin{proof}
By Theorem~\ref{thm:hyperbolic-input}, $P_{d'}(\Gamma)$ is a model for $\underline E\Gamma$ and $\mu_\Gamma$ is surjective, so every input class has a preimage in $K^\Gamma_*(P_{d'}(\Gamma))$. Apply Lemma~\ref{lem:app-uniform} to each parity with $d_0=d_\delta$, take the larger of the two pairs $(t,N)$, and increase $t$ to satisfy \eqref{eq:witness-floor}.
\end{proof}

\subsection{A rate function}\label{sec:rate}
\begin{lemma}[Finitely many bounded-hyperbolicity markings]\label{lem:finite-markings}
For a fixed marking $S$ and integer $\delta$, the set $\mathcal G_S(\delta)$ of $\delta$-hyperbolic marked groups is finite.
\end{lemma}

\begin{proof}
By Theorem~\ref{thm:hyperbolic-input}, $P_{d_\delta}(\Gamma)$ is contractible. A Cayley-graph loop therefore bounds a finite simplicial disc in this Rips complex, after subdivision. Replace each edge of that disc by a geodesic word of length at most $d_\delta$. The boundary of every triangle then gives a trivial word of length at most $3d_\delta$. It follows that the set of all trivial words of length at most $3d_\delta$ normally generates the kernel of the marked map $F(S)\to\Gamma$. Thus this finite set of short words determines the marked group. There are only finitely many possible subsets of the finite set of words of length at most $3d_\delta$.
\end{proof}

\begin{definition}[The witness rate]\label{def:rate-function}
Set $\eps_l=2^{-l-10}$ for $l\ge1$. For each $\Gamma\in\mathcal G_S(\delta)$ and each tuple $(d,r,n,l)$, choose integers
\[
 t=t(\Gamma,\delta;d,r,n,l),\qquad N=N(\Gamma,\delta;d,r,n,l)
\]
as in Lemma~\ref{lem:finite-witness} for $\eps=\eps_l$. To make a definite choice, take the first permissible pair in the ordering by $t+N$, then by $t$. Let $F^0_S(\delta,j)$ be the maximum of $16\max(d_\delta,j)$ and all the integers
\begin{equation}\label{eq:rate-size}
 2t\left\lceil\log(eN\,2t)\right\rceil
\end{equation}
for $\Gamma\in\mathcal G_S(\delta)$ and $1\le d,r,n,l\le j$. Finally define
\begin{equation}\label{eq:rate-envelope}
 F_S(\delta,j)=
 \max_{1\le\delta'\le\delta,\ 1\le j'\le j}F^0_S(\delta',j').
\end{equation}
An empty maximum over marked groups contributes nothing.
\end{definition}

Lemmas~\ref{lem:finite-witness} and~\ref{lem:finite-markings} prove that this function is finite. The envelope makes it nondecreasing in both variables. The factor $2$ in \eqref{eq:rate-size} reserves propagation for the self-adjoint defects of the path. In particular, condition \eqref{eq:fir} eventually dominates the exact matrix-size and propagation requirement for every fixed input tuple, at every stage in the specified infinite subsequence.

\begin{definition}[Fast-growing injectivity radii]\label{def:fir}
Let $G=\varinjlim G_m$ be a lacunary presentation with marking $S$
(Definition~\ref{def:lacunary}), with four-point hyperbolicity constants
$\delta_m\in\N$ of the groups $G_m$ and isometry radii $R_m$ of the limit maps
$G_m\to G$. The presentation, together with the data $(\delta_m,R_m)_m$, has
\emph{fast-growing injectivity radii} if it satisfies \eqref{eq:fir}:
\[
 R_m\ge mF_S(\delta_m,m)\qquad\text{for infinitely many }m.
\]
A lacunary hyperbolic group $G$ \emph{has fast-growing injectivity radii} if some
lacunary presentation of $G$, with some choice of this data, does. When we say this
of a group, such a presentation and its data are fixed and are part of the
hypothesis.
\end{definition}

\begin{remark}\label{rem:fir}
Condition \eqref{eq:fir} compares integer numbers and can be imposed
recursively when constructing a presentation. Along the subsequence in \eqref{eq:fir} one has $R_m\ge16m(32\delta_m+20)$, so
$R_m\to\infty$ and $\delta_m=o(R_m)$ there: that subsequence is lacunary as a
consequence of \eqref{eq:fir}. The sequence $(G_m)_m$ and the constants
$\delta_m,R_m$ are therefore retained as part of the data in the definition
and in the results below.
\end{remark}

\subsection{The main assembly theorem}
\begin{theorem}[Ordinary assembly for fast-growing injectivity radii]\label{thm:main-rate}
For each finite marking $S$ there is a finite, nondecreasing function
\[
 F_S\colon\N\times\N\longrightarrow\N
\]
with the following property. Suppose $G=\varinjlim G_m$ is a lacunary presentation with this marking which has fast-growing injectivity radii in the sense of Definition~\ref{def:fir}, and suppose $G$ has \RD.
Then $\mu_{G,\C}$ is an isomorphism in both degrees. The function is the one in Definition~\ref{def:rate-function}. The map is also quantitatively surjective with size-dependent controls as in Definition~\ref{def:qs}.
\end{theorem}

\begin{proof}[Proof of Theorem~\ref{thm:main-rate}]
Injectivity is Theorem~\ref{thm:main-inj}. For surjectivity we check the hypotheses of
Theorem~\ref{thm:qs-framework}. The approximants have ONL with asymptotic controls by
Theorem~\ref{thm:onl-lacunary}. Proposition~\ref{prop:rd-pushdown}, for some ball RD
constants $(C,s)$ of $G$, gives the pushdown control \eqref{eq:rd-h} on self-adjoint
finite matrix defects, with $N$ the active matrix size (arrow (iii) of
Figure~\ref{fig:directions}). It remains to show that the approximants have
surjectivity with asymptotic controls relative to $(K,h)$, with $A=20$.

Fix a positive integer tuple $(d,r,n,l)$ and pass to the infinite subsequence in
\eqref{eq:fir}. For $m$ in it with $m\ge\max(2,d,r,n,l)$, the group $G_m$ lies in
$\mathcal G_S(\delta_m)$. Put $k_m=\max(d,d_{\delta_m})$, and let $t_m,N_m$ be the
integers chosen in Definition~\ref{def:rate-function} for $\Gamma=G_m$, $\delta=\delta_m$
and this tuple. Lemma~\ref{lem:finite-witness} supplies, for every input over $G_m$ of
propagation $r$, size at most $n$ and tolerance $\eps_l$, a normalised cycle $F_m$ on
$P_{k_m}(G_m)$, the compression of $J^*_{\eps_l/100}(F_m)$ with error less than
$\eps_l/100$, and a path of propagation at most $t_m$, size at most $N_m$ and defects at
most $20\eps_l$. Its defects are self-adjoint finite matrices of size $N_m$. These bounds
do not depend on the input.

Since the tuple is bounded by $m$, Definition~\ref{def:rate-function} and \eqref{eq:fir}
give
\begin{equation}\label{eq:enough-radius}
 16k_m\le F_S(\delta_m,m)\le R_m/m<R_m,\qquad
 2t_m\lceil\log(eN_m\,2t_m)\rceil\le R_m.
\end{equation}
As $\log(eN_m\,2t_m)\ge\log(2e)>1$, the ceiling is at least $2$, so $2t_m\le R_m/2$.
Hence $2t_m\le h(m,N_m)$, which is \eqref{eq:asym-witness-control}. This holds for
arbitrarily large $m$, so the approximants have surjectivity with asymptotic controls.

Theorem~\ref{thm:qs-framework} now shows that $\mu_G$ is quantitatively surjective with
size-dependent controls in both degrees, and Proposition~\ref{prop:qs-ordinary} gives
ordinary surjectivity.

The tolerance $\eps_l$ and the input parameters are fixed before $m$ is chosen;
the maximum in $F_S(\delta_m,m)$ covers every tuple with entries at most $m$.
\end{proof}

\subsection{Examples}\label{sec:examples}
We use graphical small cancellation theory in Gruber's formulation
\cite{Gru15}. A graphical presentation is given by a labelled graph
$\Theta$, and $G(\Theta)$ is the quotient of the free group on the alphabet
by the words read on closed paths in $\Theta$. Classical presentations
arise from disjoint unions of cycle graphs. A subfamily of components
still satisfies the same small cancellation condition, since deleting
components introduces no new pieces.

\begin{theorem}[{Gruber \cite[Theorem~2.16, Lemma~5.5]{Gru15}}]\label{thm:gruber}
Let $\Theta=\Gamma\sqcup\bigsqcup_n\Theta_n$ satisfy $Gr(7)$, with
$\Gamma$ finite and the $\Theta_n$ finite, connected, and pairwise
non-isomorphic as labelled graphs. Then $G(\Gamma)$ is hyperbolic, and
the isometry radii of $G(\Gamma)\to G(\Gamma\sqcup\Theta_n)$ tend to
infinity as $n\to\infty$.
\end{theorem}

The theorem in the reference uses injectivity on balls. It gives the
stated conclusion after the fixed change of radius explained in
Section~\ref{sec:afc}. Classical $C'(1/6)$ and strong graphical
$C'(\lambda)$, $\lambda\le1/6$, satisfy the required graphical condition
\cite[Lemma~2.3 and the subsequent remarks]{Gru15}.

\begin{construction}[Small cancellation groups with fast-growing injectivity radii]\label{constr:lacunary}
Fix an infinite $Gr(7)$-labelled family as above. Choose its first retained
component, call the resulting finite graph $\Gamma_1$, and put
$G_1=G(\Gamma_1)$ and $L_0=0$. Having chosen $\Gamma_m$ and $G_m$, choose
an integer hyperbolicity bound $\delta_m\ge1$ and put
\[
 L_m=\max\{L_{m-1},mF_S(\delta_m,m)\}.
\]
Choose the next unused component so far out in the family that its
addition gives a quotient $G_m\to G_{m+1}$ which is isometric on balls
of radius $L_m+1$. Let $G$ be the group of all retained components.
\end{construction}

The choice is possible by Theorem~\ref{thm:gruber}. It uses only the
already chosen finite prefix and the finite number $F_S(\delta_m,m)$.

\begin{proposition}\label{prop:construction}
The output of Construction~\ref{constr:lacunary} is lacunary hyperbolic
and has isometry radii $R_m\ge mF_S(\delta_m,m)$ for its limit maps.
\end{proposition}
\begin{proof}
Every later quotient is isometric on balls of radius $L_m+1$.
Lemma~\ref{lem:inj-limit} therefore makes $L_m$ an isometry radius for
$G_m\to G$. Since $F_S(\delta,m)\ge16(32\delta+20)$, these radii tend
to infinity and $\delta_m/R_m\to0$.
\end{proof}

Classical examples are obtained from an infinite $C'(1/6)$ family of
distinct cyclically reduced relators, retaining a subfamily by the
construction. In particular, one may start with a stronger
$C'(1/20)$ family. We record both analytic routes available for them.

\begin{theorem}[{Sledd \cite{Sle21}; Bodart--Osin \cite[Theorem~4.9]{BO26}}]\label{thm:classical-inputs}
Every finitely generated classical $C'(1/6)$ group has Assouad--Nagata
dimension at most $2$ and satisfies
\begin{equation}\label{eq:classical-rd}
 \|\lambda_G(f)\|\le(10r+1)^3\|f\|_2
 \qquad(\supp f\subseteq B_G(e,r)).
\end{equation}
The defining presentation may be infinite.
\end{theorem}

\begin{theorem}\label{thm:c16}
Let $G$ have a classical $C'(1/6)$ presentation and a lacunary
presentation with fast-growing injectivity radii (Definition~\ref{def:fir}). Then its ordinary reduced
assembly map is an isomorphism in both degrees and is quantitatively
surjective with size-dependent controls. Such groups exist by
Construction~\ref{constr:lacunary}.
\end{theorem}
\begin{proof}
Equation~\eqref{eq:classical-rd} supplies \RD, so
Theorem~\ref{thm:main-rate} applies to the approximating sequence in the
hypothesis. Proposition~\ref{prop:qs-ordinary} gives the stated
quantitative formulation.

There is also a direct ONL pushdown proof. Sledd's theorem gives finite
asymptotic dimension with linear control, hence metric sparsification
and ONL with linear control \cite{CTWY08}. Rescale the ONL constant to
$1/2$ and write the control as $ar+b$. Lemma~\ref{lem:onl-pushdown}
gives $K=2$ and $h(m)$ growing linearly with $R_m$. For every fixed input
tuple, the finite witnesses used to define $F_S$ satisfy
$t_m\le R_m/(2m)$ on the subsequence in \eqref{eq:fir}. Thus $2t_m\le h(m)$ for
large $m$. Lemma~\ref{lem:finite-witness} and
Theorem~\ref{thm:qs-framework} give the same conclusion in both degrees.
\end{proof}

\begin{remark}
The Haagerup and Higson--Kasparov theorems \cite{AO15,HK01} prove
the Baum--Connes conjecture with coefficients for every classical $C'(1/6)$ group. Proposition~\ref{prop:qs-ordinary}
also supplies the existential size-dependent formulation in that
generality. The example above illustrates the two directions of norm
control in the covering argument.
The RD statement in \eqref{eq:classical-rd} extends the earlier
$C'(1/10)$ result \cite{AD12}.
\end{remark}

\subsection{Graphical groups containing expanders}\label{sec:construction}
For the expander examples, Appendix~\ref{app:rd} supplies rapid decay. It
needs the defining graphs to be sparse: each girth must be large compared
with the previous one \eqref{rd-eq:separation}. Condition \eqref{eq:fir} is a
requirement of the same kind. Both are lower bounds on the next girth in terms
of the graphs already chosen, so one subsequence meets both. This is the
content of Theorem~\ref{thm:main-expanders}.

\begin{lemma}[Finite prefixes and limit radii]\label{lem:graphical-radii}
Let a strong graphical $C'(1/20)$ presentation be given by finite connected graphs with increasing girths $\sigma_i$ and $\diam\Theta_i\le A\sigma_i$. Write $G_m$ for the group defined by its first $m$ graphs and $G$ for the full group. Then:
\begin{enumerate}[label=(\roman*)]
\item $G_m$ has $2A\sigma_m$-slim geodesic triangles;
\item it is $\delta_m$-hyperbolic in the four-point convention for
\begin{equation}\label{eq:delta-prefix}
 \delta_m=\lceil16A\sigma_m+16\rceil;
\end{equation}
\item a nonidentity element in $\ker(G_m\to G)$ has word length at least $\sigma_{m+1}$, and hence $R_m=\lfloor\sigma_{m+1}/8\rfloor$ is an isometry radius once this number is positive.
\end{enumerate}
\end{lemma}

\begin{proof}
Assertion (i) follows from \cite[Lemma~5.8]{Gru15}, using the maximum diameter of a defining component. The standard comparison from $\eta$-slim triangles to the four-point convention permits, for example, $\delta=4\eta+4$ on the vertex set \cite[III.H.1]{BH99}; the larger value in \eqref{eq:delta-prefix} therefore suffices.

For (iii), represent a nontrivial kernel element by a shortest word in $G_m$. It becomes trivial using finitely many relators of the full presentation. A minimal graphical diagram must use a face from some component of index greater than $m$; otherwise the word was already trivial in $G_m$. That face has perimeter at least $\sigma_{m+1}$. By \cite[Lemma 5.9]{Gru15}, each face perimeter is at most the length of the whole boundary word. This proves the lower bound directly for the map to the limit. The isometry-radius assertion follows from $4R_m<\sigma_{m+1}$ and the convention in Section~\ref{sec:afc}.
\end{proof}

\begin{theorem}[Expander examples]\label{thm:main-expanders}
Let $(\Theta_i)_i$ be an expander family of finite connected graphs of degree
at most $D$, with $\girth\Theta_i\to\infty$ and
$\diam\Theta_i\le A\girth\Theta_i$. Then there are a subsequence
$(\Theta_{i_m})_m$, with girths $\sigma_m$, and a labelling of it by a finite
alphabet $S$ with the following properties. Let $G_m$ be the graphical group of
the first $m$ graphs and $G=\varinjlim G_m$ the graphical group of the whole
subsequence.
\begin{enumerate}[label=(\roman*)]
\item The labelled subsequence satisfies hypotheses (H1)--(H4) of
Appendix~\ref{app:rd} with $\lambda=1/20$ and $\varkappa=1$. Each $\Theta_{i_m}$
embeds isometrically in the Cayley graph of $G$.
\item $G$ and all $G_m$ have ball-form \RD{} with exponent $13/2$ and one
common constant.
\item $G$ is lacunary hyperbolic and has fast-growing injectivity radii
(Definition~\ref{def:fir}). The lacunary presentation is $G=\varinjlim G_m$,
with $\delta_m=\lceil16A\sigma_m+16\rceil$ and
$R_m=\lfloor\sigma_{m+1}/8\rfloor$, and \eqref{eq:fir} holds for every $m$.
\item $\mu_{G,\C}$ is an isomorphism in degrees $0$ and $1$.
\end{enumerate}
The labelling is Osajda's \cite[Theorem~2.7]{Osa20} with one further finite
label coordinate (Theorem~\ref{rd-thm:construction}).
\end{theorem}

Large-girth expander families with this diameter bound are among those used in
\cite[\S3]{Osa20}. For these groups the assembly map is an isomorphism with trivial coefficients
and fails to be surjective for some coefficient algebra
\cite[Corollaries 3.3--3.4]{Osa20}, while Theorem~\ref{thm:main-inj} gives
injectivity for all separable coefficients.

\begin{proof}[Proof of Theorem~\ref{thm:main-expanders}]
By Theorem~\ref{rd-thm:construction}, after passing to a
subsequence there is a labelling by a finite alphabet $S$ satisfying
(H1)--(H3) with $\lambda=1/20$ and $\varkappa=1$, and every further subsequence
satisfies them with the same alphabet. From now on $S$ is fixed. Put
$c=A\log(2|S|)$.

We now choose the subsequence recursively, starting with any first graph. Suppose graphs of
girths $\sigma_1<\dots<\sigma_m$ have been chosen, and put
$\delta_m=\lceil16A\sigma_m+16\rceil$. Choose the next graph with girth
\begin{equation}\label{eq:joint-choice}
 \sigma_{m+1}\ge
 \max\left\{100e^{c\sigma_m},\ 16m\,F_S(\delta_m,m)+16\right\}.
\end{equation}
This is possible: the girths tend to infinity, and the right side depends
only on $S$, $A$, $m$ and $\sigma_m$. In particular $F_S$ is the function of
Definition~\ref{def:rate-function} for the marking $S$. It is defined on all
hyperbolic groups with this marking, so it is fixed before we choose the
subsequence.

The first bound in \eqref{eq:joint-choice} is (H4). So
(H1)--(H4) hold, and Theorem~\ref{rd-thm:main} gives \textup{(ii)}. The
embedding in \textup{(i)} is \cite[Lemma~2.4]{ACGH}.

To check the assembly rate, apply Lemma~\ref{lem:graphical-radii}. It shows that $G_m$ is
$\delta_m$-hyperbolic and $R_m=\lfloor\sigma_{m+1}/8\rfloor$ is an isometry
radius for $G_m\to G$. The second bound in \eqref{eq:joint-choice} gives
$R_m\ge mF_S(\delta_m,m)$ for every $m$. Since
$F_S(\delta,m)\ge16(32\delta+20)$, also $R_m\to\infty$ and $\delta_m=o(R_m)$.
So $G=\varinjlim G_m$ is a lacunary presentation with fast-growing injectivity
radii, which is \textup{(iii)}, and Theorem~\ref{thm:main-rate} gives
\textup{(iv)}.
\end{proof}

\begin{remark}\label{rem:irrational-traces}
We record the choices that combine Theorem~\ref{thm:main-expanders} with the
construction of \cite{OAI26}. First choose $0<\lambda<1/24$ with
$12(L_0+1)\lambda\le1/16$, where $L_0$ is the universal constant in
\cite[equation~(4)]{OAI26}. Take the regular expander family and Osajda labelling
provided by \cite[Proposition~3.1]{OAI26}, and apply the finite-coordinate
refinement in Theorem~\ref{rd-thm:construction} with this value of $\lambda$.
Refinement preserves the bound on matching labelled paths, while the intrinsic
graph metrics and expansion constants remain the same. Fix the resulting finite
alphabet $S$, and then pass to a further subsequence satisfying
\eqref{eq:joint-choice}. This preserves the intrinsic-ball conditions and the
diverging gaps between the integers $\lfloor|V(\Theta_i)|/2\rfloor$ required in
\cite[Proposition~3.1]{OAI26}. The proof of Theorem~\ref{thm:main-expanders}, with
this smaller $\lambda$, gives rapid decay and fast-growing injectivity radii for
the resulting group $G$. Hence $G$ satisfies the ordinary reduced Baum--Connes
conjecture.

After fixing these data, choose a prime $p>2$ large enough to encode all sets of
outgoing labels from $S^{\pm1}$, and then choose the codes as in
\cite[Section~4.1]{OAI26}. The construction in \cite[Definition~4.3]{OAI26} gives
a quotient $V$ of the lamp module $\mathbb F_p^{(G)}$. Its semidirect product
$V\rtimes G$ has a projection of irrational canonical trace in its reduced
group $C^*$-algebra, so its ordinary assembly map is not surjective. Thus the
same choice of base group gives an isomorphism for its own assembly map and
the expander geometry used to construct these forbidden traces.
\end{remark}

\appendix
\section{Folklore on the Baum--Connes assembly map}\label{app:bc}
This appendix collects the standard facts about the assembly map used in
Section~\ref{sec:bc}, in the precise forms needed there: one standard module
for all classes, homotopies that change the representation, action and
operator together, separable coefficients, and finite propagation paths with
stated bounds. The notation is that of Section~\ref{sec:bc}: the Rips complex
$P_d(\Gamma)$ with its free dense subset $Z_d$ and fundamental domain $D_d$,
$L_d=\ell^2(D_d)\otimes H_0$, the cut-off function $\phi(\nu)=\nu(e)$, and
normalised cycles.

\subsection{Exactly equivariant representatives}\label{app:bc-meyer}
\begin{lemma}\label{lem:app-meyer}
Let $\Gamma$ be a finitely generated group, $d>0$, and $B$ a separable graded
$\Gamma$-$C^*$-algebra. Then $A=C_0(P_d(\Gamma))$ is a proper
$\Gamma$-algebra in the sense of \cite[\S8]{Mey00}, and the following hold
for triples for $A,B$.
\begin{itemize}
\item every class is represented by an essential triple $(E,\varphi,F_0)$ with $F_0$ an exactly equivariant self-adjoint contraction;
\item if two such triples $T_0,T_1$ define the same class, there is an essential
homotopy $(\mathcal E,\varphi,\mathcal F_0)$ with $\mathcal F_0$ exactly equivariant. Its endpoints are, up to equivariant unitary equivalence,
locally compact perturbations of $T_0$ and $T_1$, joined to them by straight lines of exactly equivariant triples.
\end{itemize}
\end{lemma}
\begin{proof}
Let $a\in C_c(P_d(\Gamma))$ and $b\in A$. The translates of $\operatorname{supp}a$ have bounded overlap $M$, so for finite $F\subseteq \Gamma$
\[
\textstyle\bigl\|\sum_{g\notin F}c_{ab}(g)^*c_{ab}(g)\bigr\|=\bigl\|\sum_{g\notin F}b^*|\alpha_g(a)|^2b\bigr\|\le M\|a\|^2\sup\{|b(x)|^2:x\in\bigcup_{g\notin F}g\operatorname{supp}a\},
\]
which tends to $0$ as $F$ increases, since these translates leave every compact set. So the coefficient functions are square-integrable. By
\cite[Proposition~8.6]{Mey00}, $A$ has Meyer's property AE: every essential triple has an exactly $\Gamma$-equivariant locally compact
perturbation. By \cite[Lemma~3.3]{Mey00}, every triple is homotopic to its essential part, whose operator is a connection, hence a locally compact
perturbation of the original one on that part. A self-adjoint contraction can be chosen \cite[proof of Proposition~3.4]{Mey00}.
\end{proof}

\subsection{Representatives on the standard module}\label{app:bc-rep}
\begin{restatedlemma}[\nameref*{lem:ample-rep}]{lem:ample-rep}
\standardmodulelemmastatement
\end{restatedlemma}
\begin{proof}
Put $X=P_d(G)$ and $A=C_0(X)$. We treat even degree; the odd case is the same with $B\otimes\mathrm{Cl}_1$. We use the graded versions of Meyer's results. For absorption of representations, treat the grading as an action of $\Z/2$ and apply Gabe--Szab\'o to $G\times\Z/2$. The almost intertwining isometries constructed below are even.

By Lemma~\ref{lem:app-meyer}, the class is represented by an essential triple $(E,\varphi,F_0)$ with $F_0$ an exactly equivariant self-adjoint
contraction. Since $A$ is proper (Lemma~\ref{lem:app-meyer}), $E$ is square-integrable \cite[Proposition~8.4, as in the proof of Proposition~8.6]{Mey00}, so
$E\oplus\mathcal H_B\cong\mathcal H_B$ equivariantly \cite[Theorem~8.5]{Mey00}. Here $\mathcal H_B\cong\ell^2(G)\otimes\ell^2(D_d)\otimes H_0\otimes B$ is the standard module, because $G$ acts
freely on $Z_d$. Add the degenerate cycle $(\mathcal H_B,M,1\otimes\sigma_0\otimes1)$, with $\sigma_0$ an odd self-adjoint unitary on $H_0$. Transporting by the unitary,
the class is $[(\mathcal H_B,\psi,F_1)]$ for an equivariant non-degenerate representation $\psi$ and an exactly equivariant self-adjoint contraction $F_1$.

We next construct almost intertwining isometries. For $a\in A$ let $\omega_a$ be its modulus of continuity for the $G$-invariant metric of
\cite[\S6.1]{OOY15}. We claim that for every $\delta>0$ there is an even equivariant isometry $W\in\mathcal L(\mathcal H_B)$ with
$\|W\psi(a)-M(a)W\|\le\omega_a(\delta)$ for all $a$. Choose $e_1,\dots,e_n\in C_c(X)$, $e_j\ge0$, with $\sum_je_j\le1$, $\sum_je_j=1$ on $\operatorname{supp}\phi$, and supports of
diameter at most $\delta$. Choose $z_j\in Z_d\cap\operatorname{supp}e_j$ ($Z_d$ is dense), mutually orthogonal graded subspaces $H_j\subseteq H_0$ with infinite-dimensional
even and odd parts, and even unitaries $S_j\colon\mathcal H_B\to H_j\otimes B$ of the form $s_j\otimes\mathrm{id}_B$, with $s_j\colon\ell^2(Z_d)\otimes H_0\to H_j$ a graded unitary. Put
\[
W_0\xi=\textstyle\sum_j\delta_{z_j}\otimes S_j\psi(\sqrt{e_j})\xi,\qquad W=\sum_{g\in G}\gamma_gW_0\gamma_g^{-1}\psi(\sqrt{g\phi}).
\]
The summands of $W$ have orthogonal ranges, in the spans of $\delta_{gz_j}\otimes H_j\otimes B$, because $G$ acts freely on $Z_d$ and the $H_j$ are orthogonal. For finite
$F\subseteq G$ the partial sums $W_F$ satisfy $W_F^*W_F=\psi(\sum_{g\in F}g\phi)$, so $W_F\xi$ converges for every $\xi$, because $\psi$ is non-degenerate; the adjoints converge in the same way, as
row operators. So $W$ converges strictly, and $W^*W=\sum_g\psi(\sqrt{g\phi})\psi(g\textstyle\sum_je_j)\psi(\sqrt{g\phi})=\sum_g\psi(g\phi)=1$. Re-indexing the sum shows $W\gamma_h=\gamma_hW$. Next,
\[
W\psi(a)-M(a)W=\textstyle\sum_g\gamma_g\bigl(W_0\psi(g^{-1}a)-M(g^{-1}a)W_0\bigr)\psi(\sqrt\phi)\gamma_g^{-1},
\]
and $(W_0\psi(b)-M(b)W_0)\eta=\sum_j\delta_{z_j}\otimes S_j\psi\bigl(\sqrt{e_j}(b-b(z_j))\bigr)\eta$ has inner product at most $\omega_b(\delta)^2\langle\eta,\eta\rangle$. As
$\omega_{g^{-1}a}=\omega_a$, $\langle\gamma_g\zeta,\gamma_g\zeta\rangle=\beta_g(\langle\zeta,\zeta\rangle)$ with $\beta_g$ positive, the summands have orthogonal ranges, and $\sum_g\langle\psi(\sqrt{g\phi})\xi,\psi(\sqrt{g\phi})\xi\rangle=\langle\xi,\xi\rangle$, the claim follows.

To apply the absorption theorem, let $B'=\K(\mathcal H_B)$ with the action $\operatorname{Ad}\gamma$. It is separable and strongly stable, since
$\mathcal H_B\otimes\ell^2\cong\mathcal H_B$ through $H_0$. Regard $\psi$ and $M$ as equivariant representations $A\to\mathcal M(B')=\mathcal L(\mathcal H_B)$, i.e.\ cocycle representations with
trivial cocycle in the sense of Gabe--Szab\'o \cite[Definition~2.2]{GS25}. For a contraction $b\in B'$ put $c=Wb$. Then
$c^*\gamma_gc\gamma_g^{-1}=b^*\gamma_gb\gamma_g^{-1}$ exactly, and $\|b^*\psi(a)b-c^*M(a)c\|\le\omega_a(\delta)$. Since $\sup_{a\in\mathcal F}\omega_a(\delta)\to0$ for every compact $\mathcal F\subseteq A$, $\psi$ is weakly contained in $M$
\cite[Definition~4.3]{GS25}, and by \cite[Corollary~4.12]{GS25} $M^\infty$ is asymptotically unitarily equivalent to $\psi^\infty\oplus M^\infty$. Forming
infinite repeats with even isometries in $1\otimes\mathcal B(H_0)\otimes1$, which are invariant and commute with $M$, gives $M^\infty=M$. Asymptotic unitary
equivalence is symmetric (replace $u_s$ by $u_s^*$) and transitive (multiply the paths), and it is preserved by Cuntz sums and by
conjugation with invariant unitaries, because $B'$ is an ideal in $\mathcal M(B')$. Since $\psi^\infty$ is equivalent to $\psi\oplus\psi^\infty$ by an invariant unitary, it follows that $\psi\oplus M$ is
asymptotically unitarily equivalent to $M$ \cite[Definition~4.1]{GS25}. So there is a norm-continuous path $(u_s)_{s\ge0}$ of unitaries
$u_s\colon\mathcal H_B\oplus\mathcal H_B\to\mathcal H_B$ such that $u_s(\psi\oplus M)(a)u_s^*-M(a)$, $u_s(\gamma_g\oplus\gamma_g)u_s^*-\gamma_g$ and $u_s(\epsilon\oplus\epsilon)u_s^*-\epsilon$ are compact
and tend to $0$ in norm as $s\to\infty$, where $\epsilon$ is the grading. For the action comparison, compactness means that the difference is compact after composition on the right with $\gamma_g^{-1}$; the difference itself is semilinear if the action on $B$ is nontrivial.

We use this path to construct a homotopy of cycles. Fix $t\ge0$ and put $\rho_s=\operatorname{Ad}u_s(\psi\oplus M)$, $\gamma^s=\operatorname{Ad}u_s(\gamma\oplus\gamma)$ and $\epsilon^s=\operatorname{Ad}u_s(\epsilon\oplus\epsilon)$ for
$s\in[t,\infty]$, with the value at $s=\infty$ given by $(M,\gamma,\epsilon)$. These are norm-continuous up to $s=\infty$ and differ from their values at $s=t$ by compacts.
Let $F_2(s)$ be the $\epsilon^s$-odd part of $u_t(F_1\oplus1\otimes\sigma_0\otimes1)u_t^*$. Then $(C([t,\infty],\mathcal H_B),(\rho_s),(F_2(s)))$, with the
pointwise action and grading, is a cycle for $A,B[t,\infty]$: each condition holds at $s=t$, and the corrections are compact and
norm-continuous in $s$. At $s=t$ it is unitarily equivalent to $(\mathcal H_B,\psi,F_1)\oplus(\mathcal H_B,M,1\otimes\sigma_0\otimes1)$. At $s=\infty$ it is $(\mathcal H_B,M,F')$, with
$F'=F_2(\infty)$ self-adjoint and odd, and with $(\gamma_gF'\gamma_g^{-1}-F')M(f)$ compact.

Finally, we average the resulting operator. We may assume $\|F'\|\le1$ by functional calculus. Roe's average $\sum_gg\sqrt\phi F'\sqrt\phi g^{-1}$ is exactly equivariant, and
its difference from $F'$ is
$\sum_g\sqrt{g\phi}[F',\sqrt{g\phi}]+\sum_g\sqrt{g\phi}(gF'g^{-1}-F')\sqrt{g\phi}$. Multiplied by $f\in C_c$, on the left or, using the adjoint form, on the right, only finitely
many terms survive, each compact. So it is a locally compact perturbation, and a contraction.\qedhere

\end{proof}

\subsection{Assembly as an index map}\label{app:bc-index}
\begin{restatedlemma}{lem:roe-coeff}
\assemblyindexlemmastatement
\end{restatedlemma}
\begin{proof}
For $B=\C$ this is Roe's theorem \cite{Roe02}, applied to the proper cocompact $G$-space $X=P_d(G)$. Roe works on the right
side of his diagram with the Kasparov cycle given by an exactly equivariant $F$ on $L^2(X)$. He forms the Hilbert
$C^*_r(G)$-module $L^2(X)_G$ (his Lemmas~2.1--2.2) and identifies $C^*_G(X)$ with $\K(L^2(X)_G)$ (Lemma~2.3). He identifies the range of
the basic projection in the descended module with $L^2(X)_G$ (Lemma~3.4) and checks that the descended operator
agrees with $F$ there, modulo compacts. So Kasparov's assembly of $[F]$ is the class of $F$ on $L^2(X)_G$. Its index is
$[I(F)]-[e]$ in degree $0$ and $[\exp(-2\pi iP_F)]$ in degree $1$, the standard boundary-map formulas
\cite[Definitions~2.8.1, 2.8.4 and~2.8.5]{WY20}.

For general $B$ we write the argument out with coefficients; for $B=\C$ it is Roe's. Put $\mathcal H=\mathcal H_B$, with the action
$\gamma_g=g\otimes1\otimes\beta_g$ and the multiplication representation $M$. Write $\mathcal H_c$ for its finitely supported vectors, and $c=\sqrt\phi$, so
that $\sum_g(gc)^2=1$. Represent the class by a self-adjoint, exactly equivariant $F$ of finite propagation (Lemma~\ref{lem:ample-rep} and
averaging). We use $\langle\gamma_g\xi,\gamma_g\eta\rangle_B=\beta_g\langle\xi,\eta\rangle_B$, $\gamma_gM(f)\gamma_g^{-1}=M(gf)$ and $u_gbu_g^*=\beta_g(b)$ throughout.

Following Roe's Lemmas~2.1--2.2, we first construct the module. On $\mathcal H_c$ put
\[
\langle\xi,\eta\rangle=\textstyle\sum_g\langle\xi,\gamma_g\eta\rangle_B\,u_g\in B\rtimes_rG,\qquad\xi\cdot(bu_h)=\gamma_h^{-1}(\xi b).
\]
The sum is finite, and the inner product is $B\rtimes_rG$-linear and Hermitian. It is positive: in the regular representation of $B\rtimes_rG$ on
$\ell^2(G)\otimes B$, given by $(b\zeta)(h)=\beta_{h^{-1}}(b)\zeta(h)$ and $(u_g\zeta)(h)=\zeta(g^{-1}h)$, one computes for finitely supported $\zeta$
\[
\textstyle\langle\zeta,\langle\xi,\xi\rangle\zeta\rangle=\bigl\langle\sum_h\gamma_{h^{-1}}(\xi)\zeta(h),\sum_h\gamma_{h^{-1}}(\xi)\zeta(h)\bigr\rangle_B\ge0.
\]
Let $(\mathcal H)_G$ be the completion. The same computation shows that $U\colon\xi\otimes\zeta\mapsto\sum_h\gamma_{h^{-1}}(\xi)\zeta(h)$ is a unitary
$(\mathcal H)_G\otimes_{B\rtimes_rG}(\ell^2(G)\otimes B)\to\mathcal H$; it is onto because $U(\xi\otimes\delta_eb)=\xi b$. Since the regular representation is faithful,
$\alpha(T)=U(T\otimes1)U^*$ is an isometric $*$-homomorphism $\mathcal L((\mathcal H)_G)\to\mathcal L(\mathcal H)$. Moreover
$\langle\delta_e\otimes v\otimes b,\gamma_g(\delta_e\otimes w\otimes b')\rangle_B=0$ for $g\ne e$, so $v\otimes bx\mapsto(\delta_e\otimes v\otimes b)\cdot x$ is a unitary
$L_d\otimes(B\rtimes_rG)\to(\mathcal H)_G$, and $\K((\mathcal H)_G)=\K(L_d)\otimes(B\rtimes_rG)$.

We now identify its compact operators with the Roe algebra, as in Roe's Lemma~2.3. For $\xi,\eta\in\mathcal H_c$, evaluating both sides on $U(\zeta'\otimes\delta_eb)=\zeta'b$ gives
\[
\textstyle\alpha(\theta_{\xi,\eta})=\sum_g\gamma_g\theta^{\mathcal H}_{\xi,\eta}\gamma_g^{-1},\qquad\theta^{\mathcal H}_{\xi,\eta}\zeta=\xi\langle\eta,\zeta\rangle_B,
\]
an invariant, locally compact operator of finite propagation. Conversely, let $T$ be invariant and locally compact, of finite propagation.
Then $T=\sum_gM(g\phi)T=\sum_g\gamma_g\bigl(M(\phi)T\bigr)\gamma_g^{-1}$, where $M(\phi)T$ is compact and supported in a fixed compact set. Approximate it by finite
sums of $\theta^{\mathcal H}_{\xi,\eta}$ with $\xi,\eta\in\mathcal H_c$ supported near that set. The translates of a compact set have bounded overlap, so
$K\mapsto\sum_g\gamma_gK\gamma_g^{-1}$ is bounded on such $K$, and $\alpha(\K((\mathcal H)_G))$ is dense in, hence equal to, $C^*(P_d(G),B)^G$. Under the unitary constructed above
the filtration by propagation corresponds to $\K(L_d)\odot C_c(G,B)_r$, with our word-propagation convention.

To compare this index with assembly, we follow Roe's Lemmas~3.1--3.4. Kasparov's descended module is $\mathcal H\rtimes_rG=\mathcal H\otimes_B(B\rtimes_rG)$,
with $fu_g(\eta\otimes x)=M(f)\gamma_g\eta\otimes u_gx$ for $f\in C_0(X)$, and $\tilde F=F\otimes1$ \cite{Kas88}. The basic projection is
$p=\sum_gc\,(gc)\,u_g$, and $p=p^*=p^2$ because $\sum_g(gc)^2=1$, as in Roe's Lemma~3.3. For $\xi\in\mathcal H_c$ put
\[
\textstyle I(\xi)=\sum_gM(c)\gamma_g\xi\otimes u_g,
\]
a finite sum. Using $\langle\gamma_g\xi,M(c^2)\gamma_h\eta\rangle_B=\beta_g\langle\xi,M(g^{-1}c^2)\gamma_{g^{-1}h}\eta\rangle_B$ and $u_g^*bu_h=\beta_{g^{-1}}(b)u_{g^{-1}h}$,
\[
\textstyle\langle I(\xi),I(\eta)\rangle=\sum_{g,k}\langle\xi,M(g^{-1}c^2)\gamma_k\eta\rangle_Bu_k=\langle\xi,\eta\rangle,
\]
and $I(\xi\cdot bu_h)=I(\xi)\cdot bu_h$. So $I$ extends to an isometry of Hilbert $B\rtimes_rG$-modules. Next, $pI(\xi)=\sum_lM(c\sum_k(kc)^2)\gamma_l\xi\otimes u_l=I(\xi)$, and
for $\eta\in\mathcal H_c$
\[
\textstyle p(\eta\otimes bu_h)=\sum_kM\bigl(c\,(kc)\bigr)\gamma_k(\eta b)\otimes u_{kh}=I\bigl(\gamma_h^{-1}(M(c)\eta b)\bigr).
\]
So $I$ is a unitary from $(\mathcal H)_G$ onto $p(\mathcal H\rtimes_rG)$.

It remains to compare the operators under this unitary. Since $F$ is equivariant, $\tilde FI(\xi)-I(F\xi)=\sum_gK\gamma_g\xi\otimes u_g$ with $K=[F,M(c)]$. This $K$ is compact, and
$K=KM(c')$ for any $c'\in C_c(X)$ equal to $1$ on the $\prop(F)$-neighbourhood of $\operatorname{supp}c$. So $\tilde FI-IF=(K\otimes1)J$ with
$J(\xi)=\sum_gM(c')\gamma_g\xi\otimes u_g$. The inner-product computation above gives $\langle J\xi,J\xi\rangle=\langle\xi,M(h')\xi\rangle$ with $h'=\sum_gg^{-1}(c'^2)$ bounded and
invariant, so $J$ is bounded. And $K\otimes1$ is compact, because $B$ acts on $B\rtimes_rG$ by compact operators. Hence $\tilde FI-IF$ is compact.

Kasparov's assembly of $[F]$ is $[p]\otimes j^G_r[F]$, the class of $p\tilde Fp$ on $p(\mathcal H\rtimes_rG)$ \cite[\S3]{Roe02}. By the unitary $I$ and the compact operator comparison, it is the
class of $F$ on $(\mathcal H)_G$, and the identification of the Roe algebra gives its index in $K_*(\K((\mathcal H)_G))=K_*(C^*(P_d(G),B)^G)$. Every class in $KK^G_*(C_0(P_d(G)),B)$ is represented on
$\ell^2(Z_d)\otimes H_0\otimes B$ by Lemma~\ref{lem:ample-rep}, and by an exactly equivariant normalised operator after averaging. Roe works with
$L^2(X,\mu)$ for a $G$-invariant measure finite on compact sets; we use instead the finitely supported vectors of $\ell^2(Z_d)$, as in
\cite[Remark~2.4]{Roe02}. Roe writes out
one parity; the other follows with Clifford gradings, as he indicates.
\end{proof}

\subsection{The free-quotient map}\label{app:bc-quotient}
\begin{proof}[Proof of Lemma~\ref{lem:transport-classes}]
Write $A=C_0(\tilde X)$. All triples below are for $A,B$, or for $A,B[0,1]$ when they are homotopies; inner products are $B$-valued. Since
$G_m$ acts on $B$ through $\pi_m$, each $k\in K_m$ acts on these modules by a $B$-linear map.

Lemma~\ref{lem:app-meyer}, applied to $G_m$ acting on $\tilde X=P_d(G_m)$, gives
exactly equivariant representatives of every class, and of homotopies between two representatives of one class, as stated there.

We first normalise these representatives. For $F_0$ exactly equivariant put
\[
N(F_0)=\sum_{g\in G_m}g\,\varphi(\sqrt{\phi_m})\,F_0\,\varphi(\sqrt{\phi_m})\,g^{-1}=\sum_g\varphi(\sqrt{g\phi_m})F_0\varphi(\sqrt{g\phi_m}).
\]
It is exactly equivariant, and a contraction by the isometry $\xi\mapsto(\varphi(\sqrt{g\phi_m})\xi)_g$. Moreover
$N(F_0)-F_0=\sum_g\varphi(\sqrt{g\phi_m})[F_0,\varphi(\sqrt{g\phi_m})]=\sum_g[\varphi(\sqrt{g\phi_m}),F_0]\varphi(\sqrt{g\phi_m})$. Multiplied on the left by $\varphi(f)$, $f\in C_c(\tilde X)$, the first
form is a finite sum of compact operators; multiplied on the right, the second form is. So $N(F_0)$ is a locally compact perturbation of
$F_0$. It has bounded support propagation: $\varphi(h)N(F_0)\varphi(h')=0$ unless some translate of $\operatorname{supp}\phi_m$ meets both supports. $N$ is natural for
equivariant unitaries and is applied fibrewise to homotopies.

We now form the free quotient. For an essential triple $(E,\varphi,F_0)$ as above, put $F=N(F_0)$. On the vectors $\xi=\varphi(f)\xi$, $f\in C_c(\tilde X)$,
define $\langle\xi,\eta\rangle^-=\sum_{k\in K_m}\langle\xi,k\eta\rangle$, a finite sum. Take a partition of unity $(\psi_U)$ on $X$ subordinate to open sets $U$ whose
preimages are disjoint unions $\bigsqcup_kk\tilde U$, and let $\tilde\psi_U$ be the lift of $\psi_U$ to $\tilde U$. Then $\sum_{U,k}k\tilde\psi_U=1$, and
\[
\textstyle\sum_k\langle\xi,k\xi\rangle=\sum_U\bigl\langle\zeta_U,\zeta_U\bigr\rangle\ge0,\qquad\zeta_U=\sum_{l\in K_m}\varphi(\sqrt{\tilde\psi_U})\,l\xi,
\]
again finite sums. Let $\bar E$ be the completion modulo null vectors, a countably generated Hilbert $B$-module, with class map $\xi\mapsto[\xi]$; then
$[k\xi]=[\xi]$ for $k\in K_m$. So $\bar g[\xi]=[g\xi]$, for any lift $g$ of $\bar g$, is a well-defined action of $G$, and $f[\xi]=[\varphi(f\circ p)\xi]$ defines $\bar\varphi$, through the
extension of $\varphi$ to bounded functions.

To describe this quotient locally, let $\tilde U\subseteq\tilde X$ be open with pairwise disjoint $K_m$-translates, and $U=p(\tilde U)$. Then $\xi\mapsto[\xi]$ is isometric on
$\varphi(C_c(\tilde U))E$, because the cross terms vanish. We implement the local identification by cut-off operators. For $\tilde u\in C_c(\tilde U)$ real, define
\[
L_{\tilde u}\colon\bar E\to E,\qquad L_{\tilde u}[\xi]=\varphi(\tilde u)\textstyle\sum_{k\in K_m}k\xi,
\]
a finite sum. Since $\langle\zeta,L_{\tilde u}[\xi]\rangle=\langle[\varphi(\tilde u)\zeta],[\xi]\rangle$, it is adjointable, with $L_{\tilde u}^*\zeta=[\varphi(\tilde u)\zeta]$, and $\|L_{\tilde u}\|\le\|\tilde u\|_\infty$ by the isometry. Let $S\in\mathcal L(E)$
satisfy $S=\varphi(\tilde f_1)S\varphi(\tilde f_2)$ with $\tilde f_i\in C_c(\tilde U)$, and let $\tilde u=1$ on $\operatorname{supp}\tilde f_1\cup\operatorname{supp}\tilde f_2$. Put $\bar S=L_{\tilde u}^*SL_{\tilde u}$, so
$\bar S[\xi]=[S\sum_kk\xi]$; it does not depend on $\tilde u$, $\bar S^*=\overline{S^*}$, and $\bar S$ is compact if $S$ is, since $L^*_{\tilde u}\theta_{x,y}L_{\tilde u}=\theta_{L^*_{\tilde u}x,L^*_{\tilde u}y}$.

We next descend the normalised operator. Let $\tilde W$ be the interior of the sheet over $W_d$ through $e$; its $K_m$-translates are disjoint because $4d<R_m$. It contains
$\operatorname{supp}\phi_m$, including its boundary. Indeed, such a measure $\nu$ is supported in $B(e,d)$. Choose any $v\in\supp\nu$ with $\nu(v)>0$. In a sufficiently small neighbourhood, $\mu(v)>0$, so $\supp\mu\subseteq B(v,d)\subseteq B(e,2d)$. Thus this whole compact support lies in the interior of the sheet. Each term
$T_g=g\varphi(\sqrt{\phi_m})F_0\varphi(\sqrt{\phi_m})g^{-1}$ of $F$ is supported in $g\tilde W$. Set $\bar T_{\bar g}[\xi]=[T_g\sum_kk\xi]$, the operator $\bar S$ above for $S=T_g$ and $\tilde U=g\tilde W$. It is
independent of the lift $g$, because $T_{k'g}=k'T_gk'^{-1}$ and $[k'\zeta]=[\zeta]$, and $\bar F=\sum_{\bar g\in G}\bar T_{\bar g}$. This is an exactly $G$-equivariant contraction, by the isometry $[\xi]\mapsto(\bar\varphi(\sqrt{\bar g\phi})[\xi])_{\bar g}$, and
$\bar F[\xi]=[F\xi]$ on compactly supported $\xi$.

To verify the cycle conditions, decompose $f\in C_c(X)$ as the finite sum
$f=\sum_{x\in G}f_x$, where $f_x=f\cdot x\phi$. Each $f_x$ is supported in
$\mathcal V_{2d}^G(x)$, using the word-ball support regions from the proof of Lemma~\ref{lem:transport}.
Choose a lift $\tilde x$ of $x$ and lift $f_x$ to $\tilde f_x$ supported in
$\mathcal V_{2d}^{G_m}(\tilde x)$. The support calculation in that proof applies to the summands $T_g$ on $E$ as well: one factor reaches
$\mathcal V_{6d}^{G_m}(\tilde x)$, and two factors reach $\mathcal V_{10d}^{G_m}(\tilde x)$.
Consequently both support variables of $(F^2-1)\varphi(\tilde f_x)$ and
$[F,\varphi(\tilde f_x)]$ lie in $\mathcal V_{10d}^{G_m}(\tilde x)$.

This compact region is contained in the open sheet
$\tilde U=\operatorname{int}\mathcal V_{12d}^{G_m}(\tilde x)$, whose $K_m$-translates are disjoint since $24d<2R_m$.
Choose a real $\tilde u\in C_c(\tilde U)$ equal to $1$ on $\mathcal V_{10d}^{G_m}(\tilde x)$.
All intermediate supports in the products defining the squared defect lie in this sheet. The local identification therefore gives
\[
\begin{aligned}
(\bar F^2-1)\bar\varphi(f_x)
 &=L_{\tilde u}^*\bigl((F^2-1)\varphi(\tilde f_x)\bigr)L_{\tilde u},\\
[\bar F,\bar\varphi(f_x)]
 &=L_{\tilde u}^*[F,\varphi(\tilde f_x)]L_{\tilde u}.
\end{aligned}
\]
Both operators are compact. Summing over $x$ and using density proves the cycle conditions for every $f\in C_0(X)$, so
$(\bar E,\bar\varphi,\bar F)$ is a triple for $C_0(X),B$. The construction is natural for equivariant unitaries and commutes with direct sums and
with evaluation at $t\in[0,1]$. It sends degenerate triples to degenerate triples, because $N(F_0)=F_0$ when $F_0$ commutes with $\varphi$.

On the standard module $E=\ell^2(Z_d(G_m))\otimes H_0\otimes B$, one gets $\bar E=\ell^2(Z_d(G))\otimes H_0\otimes B$, since $Z_d(G_m)$ is the preimage of $Z_d(G)$ and
$K_m$ acts freely. Comparing matrix entries as in Lemma~\ref{lem:transport} gives $\bar F=\Psi_m(F)$. For a normalised cycle $F_m=\sum g\sqrt\phi T\sqrt\phi g^{-1}$, the
segment $(1-s)F_m+sN(F_m)$ consists of normalised cycles, with $T_s=(1-s)T+sF_m$, and its push-down by Lemma~\ref{lem:transport} joins
$\Psi_m(F_m)$ to $\Psi_m(N(F_m))=\overline{N(F_m)}$.

We can now define $\Psi_{m*}[(E,\varphi,F_0)]=[(\bar E,\bar\varphi,\overline{N(F_0)})]$. By Lemma~\ref{lem:app-meyer} and the naturality and evaluation properties above, representatives of
one class give homotopic results, so $\Psi_{m*}$ is well defined; it is additive, and by the last paragraph $\Psi_{m*}[F_m]=[\Psi_m(F_m)]$. In odd
degree, replace $B$ by $B\otimes\mathrm{Cl}_1$. Meyer's results hold for graded $\sigma$-unital coefficients \cite[\S\S3, 8]{Mey00}, and the dictionary to
self-adjoint $F$ with $F^2-1$ locally compact is the usual one. For naturality in $d$, restriction along $C_0(P_{d'})\to C_0(P_d)$ preserves
essentiality and exact equivariance, $\phi(\nu)=\nu(e)$ does not depend on $d$, and $N$ and the free quotient commute with restriction.
\end{proof}

\subsection{Finite propagation witnesses}\label{app:bc-witness}
Equalities in ordinary $K$-theory of finite propagation elements are
witnessed by actual finite propagation paths, with bounds uniform over
compact sets of inputs. Both statements below are used with explicit
parameters in Sections~\ref{sec:bc} and~\ref{sec:lacunary}.

\begin{lemma}[Finite propagation paths]\label{lem:app-path}
Let $\Gamma$ be a finitely generated group, $B$ a separable
$\Gamma$-$C^*$-algebra, $\mathcal H'$ a separable Hilbert space, and
$A=\K(\mathcal H')\otimes(B\rtimes_r\Gamma)$, filtered by
$\K(\mathcal H')\odot C_c(\Gamma,B)_r$. Let $0<\eps<1/128$. All elements
below lie in matrix algebras over the unitisation, have finite propagation,
and have norm less than $2$.
\begin{enumerate}[label=(\arabic*)]
\item Let $x,y$ be self-adjoint with $\|x^2-x\|,\|y^2-y\|\le\eps$, and let
$e,e'$ be scalar projections with $[x]-[e]=[y]-[e']$ in $K_0(A)$, the classes
being those of the spectral projections. Then there is a path of finite
propagation, with all defects at most $20\eps$, joining $x\oplus e'$ and
$y\oplus e$, after adding the same scalar identity blocks to both and the
necessary zero padding.
\item Let $x,y$ have both unitary defects at most $\eps$, with $[x]=[y]$ in
$K_1(A)$, the classes being those of the polar unitaries. Then there is a
path as in \textup{(1)} joining $x$ and $y$, after stabilising both by
identity blocks.
\item If $B=\C$ and the non-scalar parts of $x$ and $y$ lie in a common
finite-dimensional fibre corner, the path in \textup{(1)} or \textup{(2)}
can be chosen with its non-scalar part in one finite-dimensional fibre
corner, that is, as a path in some $M_N(\C\Gamma)$ together with its scalar
background.
\item If $\|x'-x\|<\eps/100$, with $x'$ self-adjoint when $x$ is as in
\textup{(1)}, then the straight segment from $x'$ to $x$ has defects less
than $2\eps$, and $x'$ and $x$ have the same ordinary $K$-class.
\end{enumerate}
\end{lemma}

\begin{proof}
We first prove \textup{(1)} and \textup{(2)} by choosing exact representatives of the given classes. In even degree, put
$p_x=\chi_{(1/2,\infty)}(x)$ and define $p_y$ in the same way. The defect bound places the spectrum of $x$ within distance $2\eps$ of $\{0,1\}$, so
$\|x-p_x\|\le2\eps$; the same estimate holds for $y$. In odd degree, the two defect bounds imply that $x$ is invertible, and its polar unitary
$v_x=x(x^*x)^{-1/2}$ satisfies $\|x-v_x\|\le2\eps$. Again, the same holds for $y$; these are the functional calculus estimates of
\cite[Remark~1.4]{OOY15}.

The equality in \textup{(1)} gives a homotopy of projections between $p_x\oplus e'$ and $p_y\oplus e$ after adding a common projection $q$ and zero padding.
We can replace this extra projection by scalar blocks as follows. Add $1-q$ to both endpoints and observe that
\[
 \begin{pmatrix}
 q+\sin^2\theta(1-q)&\sin\theta\cos\theta(1-q)\\
 \sin\theta\cos\theta(1-q)&\cos^2\theta(1-q)
 \end{pmatrix},\qquad 0\le\theta\le\pi/2,
\]
is a path of projections from $q\oplus(1-q)$ to $1\oplus0$. Appending these paths at the two ends therefore gives the required scalar stabilisations.
In odd degree, equality in $K_1(A)$ gives a homotopy between $v_x$ and $v_y$ after adding identity blocks. Thus, in either degree, we obtain a
norm-continuous path $\gamma$ of projections or unitaries in one matrix algebra over the unitisation. If $X,Y$ denote the corresponding stabilisations
of the original almost representatives, then
\[
 \|X-\gamma(0)\|\le2\eps,\qquad \|Y-\gamma(1)\|\le2\eps.
\]

We now construct a finite propagation approximation to $\gamma$. Put $\delta=\eps/100$. By uniform continuity, choose a subdivision
$0=t_0<\cdots<t_m=1$ such that the oscillation of $\gamma$ on each subinterval is less than $\delta/2$. At each $t_j$, use density to choose a finite
propagation element $a_j$, with the same scalar part as $\gamma(t_j)$, such that $\|a_j-\gamma(t_j)\|<\delta/2$. In even degree we take its self-adjoint
part, which preserves this estimate. Linear interpolation between the $a_j$ gives a path $a(t)$ with
$\|a(t)-\gamma(t)\|<\delta$ for every $t$. There are only finitely many $a_j$, so the union of their group supports supplies one propagation bound
for the whole interpolated path.

The next step is to join $X$ to $a_0$ and $a_m$ to $Y$ by straight segments. Each point of the first segment is within $2\eps$ of $\gamma(0)$,
and the analogous statement holds at the other endpoint. Together with the interpolated path, these segments have a common finite propagation bound,
obtained by also including the supports of $X$ and $Y$. To estimate their defects, expansion of the products gives, for a projection $P$, a
self-adjoint $h$, and a unitary $U$,
\[
 \begin{aligned}
 \|(P+h)^2-(P+h)\|&\le3\|h\|+\|h\|^2,\\
 \max\{\|(U+h)^*(U+h)-1\|,\|(U+h)(U+h)^*-1\|\}
 &\le2\|h\|+\|h\|^2,
 \end{aligned}
\]
where the second estimate holds for arbitrary $h$. Every point of our path is within $2\eps$ of the relevant exact representative. Therefore its
norm is at most $1+2\eps<2$, and its defects are at most $6\eps+4\eps^2<20\eps$. This proves \textup{(1)} and \textup{(2)}.

To obtain \textup{(3)}, we choose the non-scalar parts in a common finite-dimensional fibre corner. The finitely many approximating values have
only finitely many compact coefficients. A finite-rank fibre projection $Q$ can approximate all these coefficients simultaneously by compression,
and we can require its range to contain the given endpoint corner. We leave the scalar parts fixed and compress the non-scalar parts by $Q$.
Choosing the initial approximations and the compression errors sufficiently small retains the bound $\|a_j-\gamma(t_j)\|<\delta/2$ used above.
Compression preserves self-adjointness and does not enlarge group supports. Moreover, the non-scalar parts of $X$ and $Y$ already lie in this corner.
Consequently the non-scalar parts of both the interpolated path and the endpoint segments lie in that corner, with the defect estimates just proved.

Finally, we prove \textup{(4)} by putting $h=x'-x$ and considering $x+sh$, $0\le s\le1$. Its group support is contained in the union of those of
$x$ and $x'$. Since $\|x\|<2$, expansion of the projection defect, or of
the two unitary defects, bounds each of them by
\[
 \eps+5\|h\|+\|h\|^2
 <\eps+\frac{5\eps}{100}+\frac{\eps^2}{10000}<2\eps.
\]
In the self-adjoint case, the spectra along this segment therefore stay uniformly away from $1/2$, so their spectral projections form a continuous
path. In the unitary case, both defects are less than $1$, so the elements remain invertible and their polar unitaries form a continuous path.
Thus the two endpoints represent the same ordinary $K$-class.
\end{proof}

\begin{lemma}[Uniform finite witnesses]\label{lem:app-uniform}
Let $\Gamma$ be a finitely generated group, let $d,r,n$ be positive integers
and $0<\eps<1/128$, and consider the inputs over $\Gamma$ of one parity, of
propagation $r$, size at most $n$, and tolerance $\eps$, as in
Definition~\ref{def:finite-input}. Suppose that the classes they represent
lie in the image of $\mu_\Gamma$. Then there are integers $d'\ge d$, $t\ge r$
and $N\ge n$, independent of the input, with the following property.

For each input there is a normalised cycle $F$ of the corresponding parity
on $P_{d'}(\Gamma)$ and a finite-rank fibre projection $Q$ such that the
polynomial index $J^*_{\eps/100}(F)$ has propagation at most $t$ and an
$\eps/100$-compression in coefficient $\ell^1$-norm to $Q$. Write $c$ for
its active finite matrix and $e_Q$ for its active scalar projection in even
degree. There is a path of propagation at most $t$ and size at most $N$, with
defects at most $20\eps$, joining the input to this compressed index after
scalar stabilisation. Explicitly, in even degree the path joins
stabilisations of $p\oplus e_Q$ and $c\oplus e_l$, with the same additional
identity blocks and the necessary zero padding. In odd degree it joins
stabilisations of $u$ and $c$ by identity blocks.

If all these classes have preimages in $K^\Gamma_*(P_{d_0}(\Gamma))$ for one
fixed $d_0$, one can take $d'=\max(d,d_0)$.
\end{lemma}

\begin{proof}
We first construct a witness for one input. Its class has a preimage at some Rips scale $d_1\ge d$, and Lemma~\ref{lem:ample-rep} represents this
preimage by a normalised cycle $F$ on $P_{d_1}(\Gamma)$. Put $\eta=\eps/100$ and $J=J^*_\eta(F)$. By Section~\ref{sec:degree-one}, $J$ has finite
propagation, compact non-scalar coefficients, and defects at most $\eta$, and its ordinary $K$-class is the assembly image of $[F]$.

We next choose the compression. Write $s$ for the scalar background of $J$, namely $e$ in even degree and $1$ in odd degree. There are only finitely
many compact coefficients in $J-s$, so we can choose a finite-rank fibre projection $Q$ such that the reconstructed compression
\[
 \widetilde c=s+Q(J-s)Q
 \quad\text{satisfies}\quad
 \|J-\widetilde c\|_1<\eta.
\]
Here matrix amplifications of $Q$ are understood, so $Q$ commutes with $s$. The coefficient $\ell^1$-norm bounds the operator norm. Consequently,
the product expansions in the preceding proof show that every defect along the straight segment from $J$ to $\widetilde c$ is bounded by
\[
 \eta+5\eta+\eta^2<7\eta<\eps.
\]
Thus functional calculus along this segment identifies their ordinary $K$-classes. If $c$ is the restriction of $\widetilde c$ to the active corner,
then in even degree its class is $[c]-[e_Q]$, where $e_Q$ is the restriction of $e$ to that corner; in odd degree it is $[c]$, with the identity on
the complementary corner. The input and this compressed index therefore represent the same class.

We may now apply Lemma~\partref{lem:app-path}{1}--(3). In even degree it gives a path from $p\oplus e_Q$ to $c\oplus e_l$ after the stated scalar
stabilisations; in odd degree it gives a path from $u$ to $c$ after adding identity blocks. In both cases the path has defects at most $20\eps$
and has its non-scalar part in a finite-dimensional fibre corner. Place the input in its fixed canonical corner and include that corner in all subsequent
enlargements. We have therefore obtained a witness of finite size and finite propagation. Increasing its propagation bound, we can also bound
the propagation of the uncompressed polynomial index $J$ and ensure that the bound is at least $r$.

We now make these choices uniform over the inputs. For each $k\le n$, the space $M_k(\C\Gamma_r)$ is finite-dimensional, because the word ball of
radius $r$ in $\Gamma$ is finite. The self-adjoint matrices satisfying $\|p^2-p\|\le\eps$ form a closed bounded subset of this space: the spectral
estimate in the preceding proof gives $\|p\|\le1+2\eps$. Similarly, the matrices satisfying both unitary defect bounds form a closed bounded
subset, since $\|u\|^2=\|u^*u\|\le1+\eps$. These input sets are therefore compact. In even degree we consider them separately for each rank
$0\le l\le k$ of the scalar projection.

For each of these finitely many compact sets, choose a finite cover by operator-norm balls of radius $\eps/100$, centred at inputs. Collect all
the centres as $x_1,\ldots,x_M$, retaining their sizes and scalar ranks. For each $x_j$, choose a topological preimage at a scale $d_j\ge d$, and put
$d'=\max_j d_j$. Naturality of assembly under the inclusions of Rips complexes sends these preimages to classes on $P_{d'}(\Gamma)$ with the same
assembly images. Choose normalised representatives there and apply the construction above to each centre. We obtain finitely many witnesses,
all at the same Rips scale, with propagation bounds $t_j$ and size bounds $N_j$.

To extend such a witness to an input $x$ in the ball centred at $x_j$, prepend the straight segment from $x$ to $x_j$, with the same scalar
stabilising blocks. By Lemma~\partref{lem:app-path}{4}, this segment has defects less than $2\eps$ and preserves the ordinary $K$-class. It has
propagation at most $r$ and stays in the fixed input corner. Therefore the same normalised cycle and compression give a witness for $x$, and
concatenation retains the defect bound $20\eps$. Taking
\[
 t=\max_j t_j,\qquad N=\max\{n,N_1,\ldots,N_M\}
\]
now gives bounds independent of the input. The dimensions $N_j$ include all matrix amplifications, scalar stabilisations, and finite-dimensional
corner enlargements used in the corresponding paths.

Under the final hypothesis, each centre has a preimage at the prescribed scale $d_0$. We can therefore start with $d_j=\max(d,d_0)$ for every
$j$, and the same argument gives the asserted choice of $d'$.
\end{proof}
\section{Rapid decay and the refined graphical construction}\label{app:rd}
\subsection{Statement of the rapid decay theorem}
\label{rd-sec:statement}

All graph metrics below are metrics on vertices, with edges of length one. Paths are reduced unless stated otherwise. For a finitely generated group $H$, write $|g|$ for word length and $\lambda_H(f)\xi(x)=\sum_a f(a)\xi(a^{-1}x)$ for left convolution. Property \RD{} means that its operator norm is bounded by $C(r+1)^\alpha\norm f_2$ for functions supported in the ball of radius $r$, with fixed $C,\alpha$. In this appendix the identity of a group is written $1$. A labelling by a finite alphabet $S$ assigns an element of $S$ to each edge with a chosen orientation, and the inverse label to the reverse orientation. It is reduced if the outgoing labels at each vertex are distinct. A piece is a reduced labelled path admitting two distinct occurrences in the disjoint union of the defining graphs. Distinct occurrences include those related by graph automorphisms.

Let $\Theta=\bigsqcup_{i\ge1}\Theta_i$, where $\Theta_i$ are finite connected graphs with finite girths $\sigma_i$. Write
\begin{equation}\label{rd-eq:groups}
 G_m=\langle S\mid\text{labels of closed paths in }\Theta_1,\ldots,\Theta_m\rangle,
 \qquad G_\infty=G(\Theta).
\end{equation}
The generating alphabet is the same at every stage. We impose the following hypotheses, for fixed constants $D\ge2$, $A>0$, $\varkappa\ge1$ and $0<\lambda\le1/20$.
\begin{enumerate}[label=(H\arabic*)]
\item Degrees are at most $D$, $\sigma_i$ strictly increase to infinity, and $\diam\Theta_i\le A\sigma_i$.
\item The labelling is reduced, and every piece occurring in $\Theta_i$ has length strictly less than $\lambda\sigma_i$. In particular, a piece occurring in both $\Theta_i$ and $\Theta_k$ has length less than $\lambda\min(\sigma_i,\sigma_k)$.
\item Every component of the off-diagonal fibre product of two defining components has fewer than $\varkappa\min(\sigma_i,\sigma_k)$ vertices. Here a fibre-product vertex is a pair of vertices, edges match equal oriented labels, and in a self-product the diagonal component is omitted.
\item With $c=A\log(2|S|)$,
\begin{equation}\label{rd-eq:separation}
 \sigma_{i+1}\ge100\exp(c\sigma_i)\qquad(i\ge1).
\end{equation}
\end{enumerate}
Hypothesis (H2) is the strong, girth-relative convention used in Osajda's construction. It is stronger than the general cycle-relative graphical $Gr'(\lambda)$ condition. We may discard a finite initial segment when arranging the hypotheses.

\begin{remark*}[The fibre product in (H3)]
The fibre product in (H3) is the standard pullback of two $S$-labelled graphs in the sense of Stallings \cite{Sta83}. Regard $\Theta_i$ and $\Theta_k$ as labelled graphs over the rose $R_S$, with one vertex and one loop for each $s\in S$. The fibre product is $\Theta_i\times_{R_S}\Theta_k$: its vertices are pairs $(u,v)$, and there is an edge $(u,v)\to(u',v')$ labelled $s$ exactly when $u\to u'$ and $v\to v'$ are both edges labelled $s$. A path in the fibre product is therefore a word that can be read simultaneously in both graphs, from the two projected starting points. In a self-product the diagonal $\{(v,v)\}$ is a copy of $\Theta_i$, the trivial overlap of a component with itself. Since the labelling is reduced, a component containing one diagonal vertex lies entirely on the diagonal, so the remaining components are genuine self-overlaps.

Each off-diagonal component is thus a maximal common labelled subgraph of the two defining graphs: the graph-theoretic form of a piece. In the Cayley graph, the intersection of two distinct copies of defining components is the image of one such component, and both projections of that component are embeddings (Lemma~\ref{rd-lem:intersection}). Bounding the components therefore bounds the overlaps of relator copies. The rapid decay estimate uses the number of vertices in these overlaps.

Hypothesis (H2) bounds the length of reduced paths in a component by $\lambda\min(\sigma_i,\sigma_k)$, and hence its diameter, and Lemma~\ref{rd-lem:intersection} shows that the components are trees. Since the labelling is reduced, the fibre product has degree at most $D$. A priori, however, a tree of diameter less than $\lambda\sigma$ and degree at most $D$ can have of order $(D-1)^{\lambda\sigma/2}$ vertices, exponential in the girth. The pullback of two labelled graphs can in general be considerably more complicated than a path: small cancellation constrains the lengths of common subwords, while the vertex count also depends on branching. Condition (H3) asks for linear growth, and $\varkappa$ is only the implied constant. Theorem~\ref{rd-thm:construction} realises it with $\varkappa=1$ by refining the labelling.
\end{remark*}

\begin{theorem}\label{rd-thm:main}
Under (H1)--(H4), there is a constant $C$ depending only on the displayed structural constants and the alphabet such that, for every $m\in\mathbb N\cup\{\infty\}$, every integer $r\ge0$, and every finitely supported complex function $f$ on $G_m$ with $\supp f\subseteq B_{G_m}(1,r)$,
\begin{equation}\label{rd-eq:main}
 \norm{\lambda_{G_m}(f)}\le C(r+1)^{13/2}\norm{f}_2.
\end{equation}
Consequently $G_\infty$ has property \RD. It is also lacunary hyperbolic.
\end{theorem}

The proof establishes the estimate directly in each group $G_m$, including $G_\infty$, with the same constants.

\begin{theorem}\label{rd-thm:construction}
Let $(\Theta_i)_i$ be finite connected graphs of degree at most $D$, with
$\girth\Theta_i\to\infty$ and $\diam\Theta_i\le A\girth\Theta_i$, and let
$0<\lambda\le1/20$. After passing to a subsequence, there is a labelling of
$(\Theta_i)_i$ by a finite alphabet $S$ with the following properties.
\begin{enumerate}[label=(\alph*)]
\item Hypotheses (H1)--(H3) hold, with this $\lambda$ and with $\varkappa=1$.
\item Every further subsequence satisfies (H1)--(H3) with the same alphabet and
the same constants.
\item Every further subsequence whose girths satisfy \eqref{rd-eq:separation}
satisfies (H1)--(H4). Such subsequences exist, and Theorem~\ref{rd-thm:main}
applies to their graphical groups.
\end{enumerate}
The labelling is the product of an Osajda labelling \cite[Theorem~2.7]{Osa20}
and one further label coordinate with $q$ letters, where $q$ depends only on
$D$ and $A$. If $(\Theta_i)_i$ is an expander family, the groups in
\textup{(c)} contain isometrically embedded expanders.
\end{theorem}

Small cancellation alone bounds the diameter of an overlap. The further
coordinate is what gives the cardinality bound (H3). It is chosen by the
Lov\'asz local lemma in Section~\ref{rd-sec:labelling}.

\subsection{Copies, intersections, and reduced diagrams}
\label{rd-sec:geometry}

Fix one of the groups in \eqref{rd-eq:groups}, denoted temporarily by $H$, and its Cayley graph $X$. Only defining components present in $H$ are considered. By \cite[Lemma 2.4]{ACGH}, each component embeds isometrically and convexly in $X$, with a unique label-preserving map after specifying a vertex and its image. Throughout, convexity means that \emph{every} geodesic between two points of the copy stays in it.

Fix a base vertex $o_i$ of each component $\Theta_i$. The \emph{standard copy} $Y^{(i)}$ is the image of the label-preserving map $\Theta_i\to X$ sending $o_i$ to the identity $1\in H$. For a vertex $v$ of $\Theta_i$, let $y_v^{(i)}\in H$ be its image under this map. Equivalently, $y_v^{(i)}$ is the element of $H$ spelled by the labels along any path from $o_i$ to $v$ in $\Theta_i$; this does not depend on the path, since the label of a closed path in $\Theta_i$ is a relator. In particular $y_{o_i}^{(i)}=1$. Put
\begin{equation}\label{rd-eq:displacement}
 z_{uv}^{(i)}=(y_u^{(i)})^{-1}y_v^{(i)},
\end{equation}
the element spelled along any path from $u$ to $v$. These elements depend on $H$ as well as on $i$, $u$ and $v$. We suppress this, since each argument below takes place in one fixed group $H$.

For $\omega\in H$, left multiplication by $\omega$ is a label-preserving automorphism of $X$. Hence $\omega Y^{(i)}$ is again a copy of $\Theta_i$, with vertices $(\omega Y^{(i)})_v=\omega y_v^{(i)}$. Conversely, a copy $D$ of $\Theta_i$ with $D_{o_i}=\omega$ equals $\omega Y^{(i)}$, by uniqueness of the label-preserving map once the image of $o_i$ is specified. Reducedness implies that an automorphism preserving all labels and fixing one vertex fixes the whole connected graph. A nontrivial label-preserving automorphism would therefore give two distinct occurrences of a path of length at least $\lambda\sigma_i$, contrary to (H2). Such automorphisms are trivial. Different components cannot be label-isomorphic for the same reason. In particular, $\omega\mapsto\omega Y^{(i)}$ is a bijection from $H$ onto the set of copies of $\Theta_i$, and we write $D=\omega Y^{(i)}$, $D_v=\omega y_v^{(i)}$.

\begin{lemma}\label{rd-lem:intersection}
For distinct copies $D,E$ of components of girths $\sigma_D,\sigma_E$, their nonempty intersection is a tree $T(D,E)$, and
\begin{equation}\label{rd-eq:intersection}
 \diam T(D,E)<\lambda\min(\sigma_D,\sigma_E),
 \qquad |V(T(D,E))|<\varkappa\min(\sigma_D,\sigma_E).
\end{equation}
It is the image of one off-diagonal fibre-product component, and both projections of that component are embeddings. A geodesic interval of length at least $\lambda\sigma_D$ in $D$ identifies $D$ among copies containing intervals of length at least their own girth times $\lambda$ starting at the same point along the same oriented geodesic.
\end{lemma}
\begin{proof}
Let $\iota_D\colon\Theta_i\to X$ and $\iota_E\colon\Theta_j\to X$ be the embeddings with images $D$ and $E$. If two vertices lie in $D\cap E$,
then every ambient geodesic between them lies in both copies, by convexity. Thus the nonempty intersection is itself convex and connected.

We first bound the lengths of reduced paths in this intersection. Such a path lifts through $\iota_D$ and $\iota_E$ to paths with the same labels
in the two defining components. If $i\ne j$, these are distinct occurrences. If $i=j$ and the lifted paths agreed at any position, the two
component maps would send the same vertex to the same point of $X$. Uniqueness of the label-preserving map would then give $D=E$, contrary to
the hypothesis. Hence the two lifts are distinct in this case as well. The path is therefore a piece, and (H2) bounds its length strictly by
$\lambda\min(\sigma_D,\sigma_E)$. If the intersection contained a cycle, repeated traversal would give reduced paths of arbitrarily large
length, contradicting this bound. Consequently the intersection is a tree $T(D,E)$. Applying the same bound to the geodesic between any two
of its vertices gives the diameter estimate in \eqref{rd-eq:intersection}.

We now identify this tree with a fibre-product component. Choose $z\in T(D,E)$, and let $u\in\Theta_i$ and $v\in\Theta_j$ be its two preimages.
Let $C$ be the component of the fibre product containing $(u,v)$. Along any path in $C$, the two projected paths have the same labels and their
images in $X$ start at the same point. Their images therefore agree at every step. This defines a map
\[
 C\longrightarrow T(D,E),\qquad (u',v')\longmapsto\iota_D(u')=\iota_E(v').
\]
Conversely, every path in $T(D,E)$ starting at $z$ lifts to such a pair of paths. Since the tree is connected, the map is onto the entire
intersection. If two vertices of $C$ have the same first coordinate, equality of their images and injectivity of $\iota_E$ force their second
coordinates to agree. Interchanging the two components proves the same assertion for the other projection; edges are determined in the same
way. Thus both projections are embeddings, and the displayed map identifies $C$ with $T(D,E)$. When $i=j$, a diagonal vertex of $C$ would again
identify the two component maps, so $C$ is off-diagonal. Hypothesis (H3) now gives
$|V(T(D,E))|=|V(C)|<\varkappa\min(\sigma_D,\sigma_E)$.

Finally, suppose that two candidate copies contain intervals starting at the same point along a fixed oriented geodesic, each of length at
least $\lambda$ times the girth of its component. Their common initial interval has length at least
$\lambda\min(\sigma_D,\sigma_E)$. The diameter estimate just proved rules this out for distinct copies. Therefore the interval determines
the copy as asserted.
\end{proof}

We use \emph{graphically reduced} diagrams: no interior edge originates in the defining graph, in the terminology of \cite[Section 2.3 and Lemma 2.4]{GScancel}. Such diagrams exist for all null words. Faces are labelled by simple cycles. In a graphically reduced diagram, adjacent faces sharing an interior edge map to distinct copies: if they mapped to the same copy, injectivity of that copy would identify the two lifts of the shared edge, so the edge would originate in the defining graph. Every interior arc is therefore a piece with the bounds of (H2).

\subsubsection*{Strebel's classification}
We recall the combinatorial input from Arzhantseva--Cashen--Gruber--Hume \cite{ACGH}. Their standing hypothesis is that the defining graph is $Gr'(1/6)$-labelled: every piece contained in a simple cycle $c$ has length less than $|c|/6$ \cite[Definition 2.3]{ACGH}. Hypothesis (H2) implies this, because a simple cycle in $\Theta_i$ has length at least $\sigma_i$ and $\lambda\le1/20$. The embedding statement used above is \cite[Lemma~2.4]{ACGH}.

We use the terminology of \cite[Definitions 3.1--3.3, 3.11]{ACGH}. A diagram is a finite, simply connected, planar $2$-dimensional CW complex; it is \emph{simple} if it is homeomorphic to a disc, and a diagram \emph{over $\mathcal R$} if its edges are labelled by $S$ and every $2$-cell reads a relator. An \emph{arc} is a maximal path of length at least one whose interior vertices have valence two; it is \emph{interior} if its interior lies in the interior of the diagram, and \emph{exterior} if it lies in the boundary. A \emph{face} is the image of a closed $2$-cell. For a face $\Pi$, the interior degree $i(\Pi)$ and exterior degree $e(\Pi)$ are the numbers of interior and exterior arcs in its boundary. A \emph{$(3,7)$-diagram} is one in which every interior vertex has valence at least three and every face with $e(\Pi)=0$ has $i(\Pi)\ge7$. A \emph{combinatorial geodesic $n$-gon} $(D,(\gamma_j)_j)$ is a $(3,7)$-diagram whose boundary is a concatenation of immersed subpaths $\gamma_0,\dots,\gamma_{n-1}$, the \emph{sides}, such that every boundary face whose exterior part is a single arc contained in one side has interior degree at least $4$. A valence-two vertex on more than one side is \emph{distinguished}. The polygon is \emph{degenerate} if merging two consecutive sides still gives a combinatorial geodesic $(n-1)$-gon. Bigons and triangles are $2$-gons and $3$-gons.

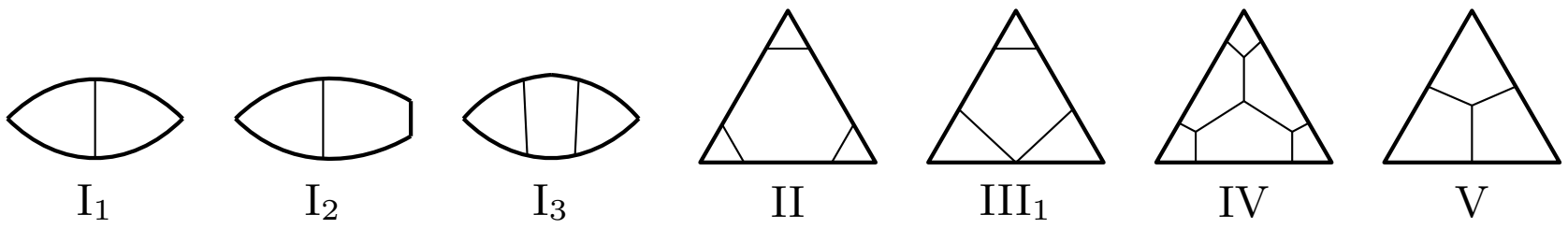
\begin{figure}[ht]
\centering
\begin{tikzpicture}[scale=0.62,font=\footnotesize,line join=round]
\begin{scope}
 \draw[thick] (0,0) .. controls (0.6,0.6) and (1.4,0.6) .. (2,0) coordinate[pos=0.5] (a);
 \draw[thick] (0,0) .. controls (0.6,-0.6) and (1.4,-0.6) .. (2,0) coordinate[pos=0.5] (b);
 \draw (a)--(b); \node at (1,-0.95) {$\mathrm I_1$};
\end{scope}
\begin{scope}[xshift=2.6cm]
 \draw[thick] (0,0) .. controls (0.6,0.6) and (1.4,0.55) .. (2,0.2) coordinate[pos=0.5] (a);
 \draw[thick] (0,0) .. controls (0.6,-0.6) and (1.4,-0.55) .. (2,-0.2) coordinate[pos=0.5] (b);
 \draw[thick] (2,0.2)--(2,-0.2); \draw (a)--(b); \node at (1,-0.95) {$\mathrm I_2$};
\end{scope}
\begin{scope}[xshift=5.2cm]
 \draw[thick] (0,0) .. controls (0.3,0.35) and (0.65,0.47) .. (1,0.5) coordinate[pos=0.7] (a1);
 \draw[thick] (1,0.5) .. controls (1.35,0.47) and (1.7,0.35) .. (2,0) coordinate[pos=0.3] (a2);
 \draw[thick] (0,0) .. controls (0.6,-0.6) and (1.4,-0.6) .. (2,0) coordinate[pos=0.37] (b1) coordinate[pos=0.63] (b2);
 \draw (a1)--(b1); \draw (a2)--(b2); \node at (1,-0.95) {$\mathrm I_3$};
\end{scope}
\begin{scope}[xshift=7.9cm,yshift=-0.5cm]
 \draw[thick] (0,0)--(2,0)--(1,1.732)--cycle;
 \draw (0.75,1.299)--(1.25,1.299); \draw (0.25,0.433)--(0.5,0); \draw (1.75,0.433)--(1.5,0);
 \node at (1,-0.45) {$\mathrm{II}$};
\end{scope}
\begin{scope}[xshift=10.5cm,yshift=-0.5cm]
 \draw[thick] (0,0)--(2,0)--(1,1.732)--cycle;
 \draw (0.75,1.299)--(1.25,1.299); \draw (0.35,0.606)--(1,0)--(1.65,0.606);
 \node at (1,-0.45) {$\mathrm{III}_1$};
\end{scope}
\begin{scope}[xshift=13.1cm,yshift=-0.5cm]
 \draw[thick] (0,0)--(2,0)--(1,1.732)--cycle;
 \draw (0.8,1.386)--(1,1.2)--(1.2,1.386); \draw (1,1.2)--(1,0.7);
 \draw (1,0.7)--(0.45,0.35)--(0.26,0.45); \draw (0.45,0.35)--(0.45,0);
 \draw (1,0.7)--(1.55,0.35)--(1.74,0.45); \draw (1.55,0.35)--(1.55,0);
 \node at (1,-0.45) {$\mathrm{IV}$};
\end{scope}
\begin{scope}[xshift=15.7cm,yshift=-0.5cm]
 \draw[thick] (0,0)--(2,0)--(1,1.732)--cycle;
 \draw (1,0.65)--(0.5,0.866); \draw (1,0.65)--(1.5,0.866); \draw (1,0.65)--(1,0);
 \node at (1,-0.45) {$\mathrm V$};
\end{scope}
\end{tikzpicture}
\caption{Strebel's shapes of combinatorial geodesic bigons and triangles, redrawn after \cite[Figure 5]{ACGH}. Thick lines are sides, thin lines interior arcs. Each shape stands for the family obtained by face combination with $\mathrm I_1$ ladders at non-degenerate corners.}
\label{rd-fig:strebel}
\end{figure}

We recall the following statements from \cite{ACGH}.
\begin{enumerate}[label=(S\arabic*),leftmargin=3em]
\item \cite[Lemma 3.4]{ACGH}, after Gruber. If a word represents the identity in $G(\Theta)$, it is the boundary word of an $S$-labelled diagram over $\mathcal R$ in which every interior arc is a piece.
\item \cite[Proposition 3.5]{ACGH}. If $P$ is an $n$-gon in $X$ whose sides are geodesics, there is an $S$-labelled diagram $D$ over $\mathcal R$ with sides reading the same words, which becomes a combinatorial geodesic $n$-gon after forgetting interior vertices of valence two. Its interior arcs are pieces.
\item \cite[Lemma 3.7]{ACGH}, Strebel's curvature formula \cite[p.~253]{Str90}. For a simple diagram without vertices of degree two,
\[
6=2\sum_v(3-d(v))+\sum_{e(\Pi)=0}(6-i(\Pi))+\sum_{e(\Pi)=1}(4-i(\Pi))+\sum_{e(\Pi)\ge2}(6-2e(\Pi)-i(\Pi)).
\]
Consequently a $(3,7)$-diagram with more than one face has at least two faces with $e(\Pi)=1$ and $i(\Pi)\le3$ (Greendlinger's lemma).
\item \cite[Lemma 3.8]{ACGH}. The sides of a combinatorial geodesic polygon are embedded, and such a polygon has at least two sides.
\item \cite[Lemma 3.10]{ACGH}. For a face $\Pi$ and a side $\alpha$, $\Pi\cap\alpha$ is empty or connected. If consecutive faces $\Pi_1,\dots,\Pi_k$ of a chain all meet $\alpha$, their union meets $\alpha$ in a connected set.
\item \cite[\S3.2]{ACGH}, combination. A \emph{vertex combination} wedges two polygons at distinguished vertices that each lie in a single face and in no other maximal cell. An \emph{edge combination} performs a vertex combination and blows the wedge point up to an interior edge. A \emph{face combination} performs an edge combination and then forgets the new edge, merging its two incident faces into one face.
\item \cite[Theorem 3.13]{ACGH}, Strebel's classification \cite[Theorem 43]{Str90}. Let $D$ be a simple diagram that is not a single face. If $D$ is a combinatorial geodesic bigon, it has shape $\mathrm I_1$. If $D$ is a combinatorial geodesic triangle, it has one of the shapes $\mathrm I_2$, $\mathrm I_3$, $\mathrm{II}$, $\mathrm{III}_1$, $\mathrm{IV}$, $\mathrm V$ of Figure~\ref{rd-fig:strebel}. Shapes $\mathrm I_2$ and $\mathrm I_3$ degenerate to bigons. Each shape stands for the infinite family obtained by performing face combination with a shape-$\mathrm I_1$ bigon, at a non-degenerate distinguished vertex, arbitrarily many times. (Strebel's further shape $\mathrm{III}_2$ belongs to a different definition of combinatorial geodesic polygon and does not occur here \cite[footnote to Theorem 3.13]{ACGH}.)
\end{enumerate}

In words, as drawn in \cite[Figure 5]{ACGH}:
\begin{itemize}
\item $\mathrm I_1$ (bigon): a \emph{ladder}, a row of faces in which consecutive faces share one interior arc (a rung) and every face meets both sides. The two corners lie in the two end faces.
\item $\mathrm I_2$: a ladder whose end face also carries a short third side; that end face contains two corners and meets all three sides.
\item $\mathrm I_3$: a ladder in which the third corner lies on the exterior of a middle face; that face meets all three sides.
\item $\mathrm{II}$: one central face meeting all three sides in arcs, and three corner faces, each meeting two sides and sharing one interior arc with the central face.
\item $\mathrm{III}_1$: one central face meeting two sides in arcs and the third side in a single vertex; one corner face opposite that vertex; and two faces at the remaining corners, each meeting two sides.
\item $\mathrm{IV}$: an interior vertex $t$ of valence three; three faces around $t$, each meeting \emph{exactly one} side and having four interior arcs; and three corner faces.
\item $\mathrm V$: an interior vertex $t$ of valence three with three interior arcs from $t$ to the three sides; three faces around $t$, each meeting two sides and containing one corner.
\end{itemize}

The three descriptions in the next lemma record what sits at the centre of the triangle: a vertex common to all three sides (i); one face meeting all three sides, which covers a single face and the shapes $\mathrm{II}$, $\mathrm{III}_1$, $\mathrm I_2$, $\mathrm I_3$ (ii); or three faces around an interior vertex, each meeting two sides, which is shape $\mathrm V$ (iii). The only other central configuration in Strebel's list is $\mathrm{IV}$, three faces around an interior vertex each meeting one side. It is excluded because $\lambda<1/8$ forbids one-sided faces with at most four interior arcs. Description (i) arises from the decomposition of a singular filling diagram into disc components and tree segments.

\begin{figure}[ht]
\centering
\begin{tikzpicture}[scale=0.95,font=\footnotesize]
\begin{scope}
 \coordinate (o) at (0,0); \coordinate (a) at (0,1.6);
 \coordinate (v) at (0.9,0.8); \coordinate (x) at (4.2,0.8);
 \draw[thick] (o)--(v)--(a);
 \draw[thick] (v) .. controls (1.8,1.2) and (3.3,1.2) .. (x)
   coordinate[pos=0.2] (u1) coordinate[pos=0.4] (u2) coordinate[pos=0.6] (u3) coordinate[pos=0.8] (u4);
 \draw[thick] (v) .. controls (1.8,0.4) and (3.3,0.4) .. (x)
   coordinate[pos=0.2] (l1) coordinate[pos=0.4] (l2) coordinate[pos=0.6] (l3) coordinate[pos=0.8] (l4);
 \foreach \i in {1,2,3,4} \draw (u\i)--(l\i);
 \fill (v) circle (1.4pt);
 \node[below left] at (o) {$1$}; \node[above left] at (a) {$a$};
 \node[right] at (x) {$x$}; \node[left] at (v) {$v\,$};
 \node at (2.1,-0.45) {(i)};
\end{scope}
\begin{scope}[xshift=5.4cm]
 \coordinate (o) at (0,0); \coordinate (a) at (0,1.6); \coordinate (x) at (4.2,0.8);
 \fill[gray!20] (0,0.35)--(0,1.25)--(0.9,1.429)--(2.6,1.105)--(2.6,0.495)--(0.9,0.171)--cycle;
 \draw[thick] (o)--(a)--(x)--cycle;
 \draw (0,1.25)--(0.9,1.429); \draw (0,0.35)--(0.9,0.171); \draw (2.6,1.105)--(2.6,0.495);
 \draw (0,0.17)--(0.45,0.086); \draw (0,1.43)--(0.45,1.514);
 \foreach \X in {3.1,3.6,3.95} \draw (\X,{1.6-0.1905*\X})--(\X,{0.1905*\X});
 \node at (1.3,0.8) {$Z$};
 \node[below left] at (o) {$1$}; \node[above left] at (a) {$a$}; \node[right] at (x) {$x$};
 \node[left] at (0,0.8) {$\gamma_a$}; \node[above] at (2.1,1.2) {$\beta$}; \node[below] at (2.1,0.4) {$\gamma_x$};
 \node at (2.1,-0.45) {(ii)};
\end{scope}
\begin{scope}[xshift=10.6cm]
 \coordinate (o) at (0,0); \coordinate (a) at (0,1.6); \coordinate (x) at (4.2,0.8);
 \coordinate (t) at (1.1,0.8);
 \fill[gray!25] (0,0.8)--(0,1.3)--(0.6,1.486)--(1.6,1.295)--(t)--cycle;
 \fill[gray!12] (0,0.8)--(t)--(1.6,0.305)--(0.6,0.114)--(0,0.3)--cycle;
 \fill[gray!35] (t)--(1.6,1.295)--(3.0,1.029)--(3.0,0.571)--(1.6,0.305)--cycle;
 \draw[thick] (o)--(a)--(x)--cycle;
 \draw (t)--(0,0.8); \draw (t)--(1.6,1.295); \draw (t)--(1.6,0.305);
 \draw (0,0.3)--(0.6,0.114); \draw (0,1.3)--(0.6,1.486);
 \foreach \X in {3.0,3.5,3.9} \draw (\X,{1.6-0.1905*\X})--(\X,{0.1905*\X});
 \fill (t) circle (1.4pt); \node[right] at (t) {$t$};
 \node at (0.45,1.15) {$Z_{12}$}; \node at (0.45,0.47) {$Z_{31}$}; \node at (2.25,0.8) {$Z_{23}$};
 \node[below left] at (o) {$1$}; \node[above left] at (a) {$a$}; \node[right] at (x) {$x$};
 \node at (2.1,-0.45) {(iii)};
\end{scope}
\end{tikzpicture}
\caption{The three descriptions of Lemma~\ref{rd-lem:diagrams}, drawn in the style of
\cite[Figures~5--6]{ACGH}; the rungs are corner chains (ladders of two-sided faces).
The triangle has vertices $1,a,x=ab$ and sides $\gamma_a$, $\beta=a\gamma_b$, $\gamma_x$; the side
$\gamma_a$ is short, $|a|\le r$. (i) A vertex $v$ on all three sides. (ii) A central face
$Z$ meeting all three sides (shapes II, III$_1$ and the degenerate I$_2$, I$_3$). (iii) A
central vertex $t$ with faces $Z_{12},Z_{31}$ meeting the short side and $Z_{23}$ meeting
$\beta$ and $\gamma_x$ (shape V). In the trilinear form \eqref{rd-eq:T}, $f$ is evaluated on
the short side, $g$ on $\beta$, and $h$ on $\gamma_x$.}
\label{rd-fig:shapes}
\end{figure}
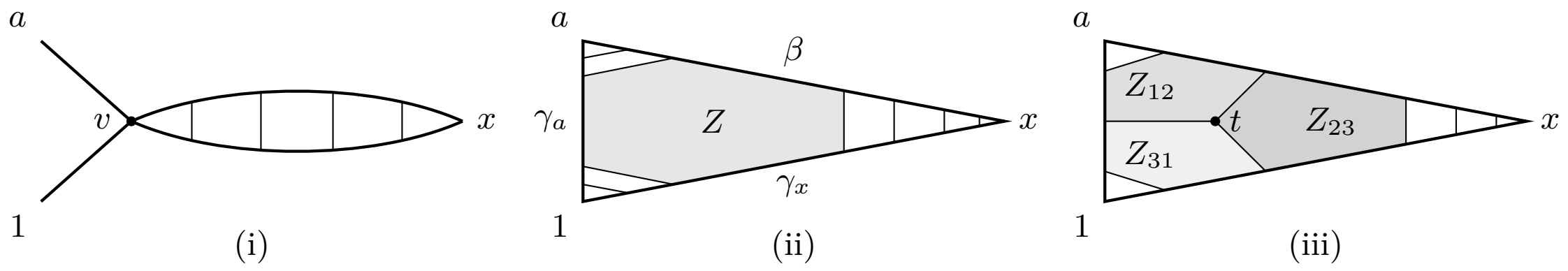

\begin{lemma}\label{rd-lem:diagrams}
A geodesic triangle admits one of the following descriptions in a graphically reduced
diagram (Figure~\ref{rd-fig:shapes}).
\begin{enumerate}[label=(\roman*)]
\item A vertex is on all three sides.
\item There is a central face $Z$ meeting all three sides, with vertex intersections allowed. At each of the three corners, the two incident sides are joined to $Z$ by a chain of two-sided faces, joined by interior arcs, cut vertices, or common tree segments.
\item There is a central vertex $t$ incident to three faces $Z_{12},Z_{23},Z_{31}$, where $Z_{uv}$ meets sides $u,v$. Each central face has at most three interior arcs and has one exterior arc on each of those two sides. The central faces share three interior arcs, from $t$ to the three sides. Corner chains attach outside this configuration.
\end{enumerate}
Every face in a corner chain has two exterior geodesic arcs and at most two interior arcs. A central face in (ii) has at most four interior arcs. No face has its entire exterior contained in just one side. In (iii), the copies of central faces sharing an interior arc are distinct.
\end{lemma}
\begin{proof}
Choose a graphically reduced diagram $D$ for the boundary word of the triangle, as supplied by \cite[Lemma~2.4]{GScancel}. Its interior arcs are pieces. By (H2), an interior arc on the boundary of a face $\Pi$ from $\Theta_i$ has length less than $\lambda\sigma_i\le\lambda p$, where $p\ge\sigma_i$ is the perimeter of $\Pi$.
We apply the degree argument in the proof of \cite[Proposition~3.5]{ACGH} to this same diagram. Every interior face has at least seven interior arcs. For a boundary face whose exterior is one arc contained in a side, that arc has length at most $p/2$, so the remaining boundary requires at least four interior arcs. After suppressing interior vertices of valence two, $D$ is therefore a combinatorial geodesic triangle. We retain its graphical reduction throughout the argument.

We first exclude faces whose exterior lies on one side. Let $\Pi$ be a face whose exterior part is a single arc $\alpha$ contained in one side. The side is a geodesic, and the rest of $\partial\Pi$ is another path between the endpoints of $\alpha$, so $|\alpha|\le p/2$. The rest of $\partial\Pi$ consists of $i(\Pi)$ interior arcs, so $p/2\le p-|\alpha|<i(\Pi)\lambda p$, that is, $i(\Pi)>1/(2\lambda)\ge10$. In every shape of (S7), every face has at most four interior arcs. This persists under the face combinations of (S7): they take place at distinguished vertices, which lie on two sides, so they never touch a one-sided face. Hence no face is one-sided, and in particular shape $\mathrm{IV}$, whose three faces around $t$ are one-sided with four interior arcs, does not occur. A face cannot have its whole exterior in one side with $e(\Pi)\ge2$ either: by (S5) its intersection with that side is connected, and planarity then leaves no room for an interior arc between two of its exterior arcs. This proves the penultimate assertion.

Suppose now that $D$ is simple. If it is a single face, that face meets all three sides and gives (ii) with empty corner chains. Otherwise (S7) applies, and we go through the shapes.
\begin{itemize}
\item $\mathrm{II}$: the central face is $Z$, and the three corner faces are the corner chains. This is (ii).
\item $\mathrm{III}_1$: the central face is $Z$, meeting the third side in a vertex, which (ii) allows. The corner face and the two faces at the other corners are the corner chains.
\item $\mathrm I_2$, $\mathrm I_3$: these degenerate to bigons. Take for $Z$ the face containing the extra corner. It meets all three original sides, possibly at that corner alone on one side. The rest of the ladder forms the corner chain on the far side. This is (ii).
\item $\mathrm V$: the vertex $t$ and its three faces give (iii).
\item $\mathrm{IV}$: excluded above.
\end{itemize}
Face combination with $\mathrm I_1$ ladders at the corners, as in (S7), merges the face at a corner with the end face of a ladder and so extends that corner chain. Ladder faces meet both sides of their corner and have at most two interior arcs (two rungs). Counting arcs in the figure, the central face in (ii) has at most three interior arcs in $\mathrm{II}$ and $\mathrm{III}_1$, at most three in $\mathrm I_2$ and $\mathrm I_3$ after combination, and none for a single face; this is within the stated bound of four. Each central face in (iii) has two interior arcs from $t$, plus at most one more from a combination at its corner, so at most three.

For completeness, singular diagrams reduce to this description as follows. Separate their maximal disc components and the edges not in any disc. Their incidence structure is a tree. An original geodesic side meets any disc in a connected boundary segment: an excursion would return through the same cut vertex, contradicting that the side is geodesic. A disc cannot have just one nonempty side, by (S4); its inherited sides are geodesic bigon or triangle sides. The median of the three corner paths in the incidence tree is either a point common to the three sides, or a unique central disc. If a side meets that disc only in a vertex, the vertex is common to all three sides. Otherwise the disc has three sides, and the classification just given applies. Each leg outside it consists of common tree segments and bigon discs in series. A bigon disc is either one face or an I$_1$ ladder, so these legs extend the corner chains. Finally, distinctness of adjacent central copies follows from graphical reduction, as explained above.
\end{proof}

\begin{lemma}\label{rd-lem:arc}
Suppose a face of perimeter $p$ has exactly two exterior geodesic arcs and at most $k$ interior arcs. Each exterior arc has length at least $p/2-k\lambda\sigma$, where $\sigma$ is the girth of its component, and hence at least $(1/2-k\lambda)p$. The boundary of any face from $\Theta_i$ has ambient diameter at most $A\sigma_i$.
\end{lemma}
\begin{proof}
The other exterior arc has length at most $p/2$, and the total interior length is less than $k\lambda\sigma$. Subtract these lengths from $p$. The diameter assertion follows because the entire face boundary maps into one isometric component copy.
\end{proof}

\subsection{A uniform estimate for one copy}
\label{rd-sec:compression}

For a component copy $Y$ of girth $\sigma$, let
\begin{equation}\label{rd-eq:kernel}
 K_\varphi(u,v)=\varphi(y_u^{-1}y_v)
 \quad\text{on }\ell^2(V(Y)).
\end{equation}
This is a compression of a right-translation operator; the displayed matrix fixes the convention. Only values of $\varphi$ on $Y^{-1}Y$ matter.

\begin{theorem}[Uniform compression estimate]\label{rd-thm:compression}
Assume reduced labelling, convex isometric component copies, trivial label automorphisms, and the girth-relative piece bound with $\lambda<1/6$. Then, for every $r\ge0$,
\begin{equation}\label{rd-eq:compression}
 \norm{K_\varphi}\le12(r+1)^{3/2}\norm{\varphi}_2
 \qquad(\supp\varphi\subseteq B_H(1,r)).
\end{equation}
The constant is uniform over all the groups and component copies satisfying the stated hypotheses.
\end{theorem}

To prove this, we use the fact that the matrix $K_\varphi$ sees $\varphi$ only through displacements $y_u^{-1}y_v$ between vertices of the single copy $Y$. Split these displacements at the scale $L=\lceil\lambda\sigma\rceil$, where the small cancellation hypothesis starts to act.
For long displacements, the piece bound gives rigidity: a displacement of length at least $L$ occurs in at most one matrix entry. The Hilbert--Schmidt norm then controls the long part, with constant one \eqref{rd-eq:HS-long}.

For short displacements, we use the tree structure of balls in $Y$. Below length $\sigma/2$ every reduced path in $Y$ is the unique geodesic between its endpoints, so short displacements are spelled by unique free words. Lifting $\varphi$ to the free group $F(S)$ and applying Haagerup's inequality there gives the bound $H(r)=(r+1)^{3/2}$ \eqref{rd-eq:Haagerup}.

To obtain a bound uniform in the girth from this short-displacement comparison, we apply it to powers at scales below $\sigma/2$. While $kr<\sigma/2$, the matrix $K_\varphi^k$ is entrywise dominated by the kernel of the transferred free convolution power $a^{*k}$, whose $\ell^2$ norm is at most $H(r)^{k-1}$. Splitting that kernel into short and long parts gives
\[
 \norm{K_\varphi}^k\le(U+1)H(r)^{k-1}\qquad\eqref{rd-eq:power-bound},
\]
where $U$ is the best constant for short nonnegative symmetric functions, the quantity being estimated. This is used twice. With $k=2$ it gives a crude, girth-dependent bound $U\le2L^{3/2}$ \eqref{rd-eq:U}. Feeding this back with $k\approx2L/r$, the girth-dependent factor appears only under a $k$-th root and is bounded by an absolute constant \eqref{rd-eq:remove-girth}. This is the spectral-radius (power) trick. Finally, a general complex $\varphi$ is reduced to a nonnegative symmetric one at the cost of a factor two, and the long part is added back.

\begin{proof}
We can take $r$ integral. Put $L=\lceil\lambda\sigma\rceil$. Two observations distinguish short and long entries.

We first consider short words. A reduced path in $Y$ of length less than $\sigma/2$ is the unique geodesic between its endpoints in $Y$: two distinct paths of total length less than $\sigma$ would contain a shorter cycle. Convexity makes it the unique ambient geodesic too. If two such paths from $y_u$ to $y_v$ and from $y_{u'}$ to $y_{v'}$ have the same displacement, translate the first by $y_{u'}y_u^{-1}$. Its endpoints are then $y_{u'}$ and $y_{v'}$. The translated path is an ambient geodesic, so convexity places it in $Y$, where uniqueness identifies it with the second path. Left translation preserves labels, hence their freely reduced words agree. Write $w_g$ for this unique word whenever the displacement $g$ is realised by such a path.

For long displacements, suppose that
\[
 y_u^{-1}y_v=y_{u'}^{-1}y_{v'}=g,\qquad |g|\ge L,
\]
then the translate $y_{u'}y_u^{-1}Y$ and $Y$ contain both $y_{u'}$ and $y_{v'}$. Convexity gives a common geodesic of length $|g|\ge\lambda\sigma$. The piece bound forces these copies to coincide, and triviality of label automorphisms forces $u=u'$ and $v=v'$. Thus every long displacement occurs in at most one matrix entry, and
\begin{equation}\label{rd-eq:HS-long}
 \norm{K_\psi}\le\norm{K_\psi}_{\mathrm{HS}}\le\norm{\psi}_2
 \quad\text{if }\psi(g)=0\text{ for }|g|<L.
\end{equation}

To estimate the short part, assume first $L\ge2$. Let $U$ be the maximum of $\norm{K_\psi}$ over nonnegative, inversion-symmetric functions $\psi$, supported on the realised displacements of length at most $L-1$, with $\norm\psi_2=1$. This is a finite-dimensional maximum, and $U\ge1$.

We compare these kernels with convolution on the free group. For such a symmetric function $\varphi$ supported in radius $r\le L-1$, with $\norm\varphi_2=1$, lift it to the free group $F(S)$ by putting $a(w_g)=\varphi(g)$ and zero elsewhere. The short-word observation gives $\norm a_2=1$, symmetry, and support in the free ball of radius $r$. Haagerup's sphere inequality \cite{Haa}, followed by Cauchy--Schwarz over the spheres, gives
\begin{equation}\label{rd-eq:Haagerup}
 \norm{\lambda_{F(S)}(a)}\le(r+1)^{3/2}=:H(r).
\end{equation}

We now apply this comparison to powers. Write $a^{*k}=a*\cdots*a$ for the $k$-fold convolution power in $\C F(S)$, where $(a*b)(w)=\sum_{uv=w}a(u)b(v)$ with products freely reduced. Thus $a^{*k}(w)$ sums $a(w_1)\cdots a(w_k)$ over all $k$-tuples of words whose product reduces to $w$, $\lambda_{F(S)}(a)^k=\lambda_{F(S)}(a^{*k})$, and $\norm{a^{*k}}_2=\norm{\lambda_{F(S)}(a)^{k-1}a}_2\le H(r)^{k-1}$ by \eqref{rd-eq:Haagerup}. If $kr<\sigma/2$, concatenating the paths counted by an entry of $K_\varphi^k$ gives a path in $Y$ of length at most $kr$. Free cancellation of labels is cancellation of actual backtracking edges, by reducedness. The resulting path is the unique geodesic for its endpoints. Therefore the entry is at most the coefficient of $a^{*k}$ at its freely reduced word. Restrict $a^{*k}$ to words realised by paths in $Y$ of length at most $kr$, and transfer those coefficients to $H$, obtaining a nonnegative symmetric function $b_k$. The transfer is injective on these words, so
\begin{equation}\label{rd-eq:power}
 0\le K_\varphi^k\le K_{b_k}\ \text{entrywise},\qquad
 \norm{b_k}_2\le\norm{a^{*k}}_2\le H(r)^{k-1}.
\end{equation}
Split $b_k$ into lengths below $L$ and lengths at least $L$, and use $U$ on the former and \eqref{rd-eq:HS-long} on the latter. Nonnegative entrywise domination implies operator-norm domination. Since $K_\varphi$ is self-adjoint, \eqref{rd-eq:power} yields
\begin{equation}\label{rd-eq:power-bound}
 \norm{K_\varphi}^{k}\le(U+1)H(r)^{k-1}.
\end{equation}

First take $r=L-1$ and $k=2$. This is legitimate because $2(L-1)<2\lambda\sigma<\sigma/2$. Taking the maximum in \eqref{rd-eq:power-bound} gives
\begin{equation}\label{rd-eq:U}
 U^2\le(U+1)L^{3/2},\qquad U\le2L^{3/2}.
\end{equation}
We then use larger powers to obtain a uniform bound. For $1\le r\le L-1$, choose
\[
 k=2\left\lfloor\frac{L-1}{r}\right\rfloor,
 \qquad t=\frac{L}{r+1}.
\]
Then $kr\le2(L-1)<\sigma/2$ and $k\ge (L-1)/r\ge t\ge1$. Equations \eqref{rd-eq:power-bound}--\eqref{rd-eq:U} imply
\begin{equation}\label{rd-eq:remove-girth}
 \norm{K_\varphi}\le H(r)(3t^{3/2})^{1/k}
 \le 3e^{3/(2e)}H(r),
\end{equation}
using $\log t/t\le1/e$. Radius zero is immediate.

Finally, for a general complex short function, put $\psi=(|\varphi|+|\varphi|^\vee)/2$, with $\varphi^\vee(g)=\varphi(g^{-1})$. Then $|K_\varphi|\le2K_\psi$ entrywise and $\norm\psi_2\le\norm\varphi_2$. This costs a factor two. If $r\ge L$, split $\varphi$ into short and long parts and combine this estimate with \eqref{rd-eq:HS-long}. Since $6e^{3/(2e)}+1<12$, the claimed bound follows. If $L=1$, the short part is diagonal and the long part satisfies \eqref{rd-eq:HS-long}, giving an even smaller bound.
\end{proof}

\subsection{Two analytic counting lemmas}
\label{rd-sec:analytic}

We now establish the rapid decay bound by working through the cases of Lemma~\ref{rd-lem:diagrams}. The estimate concerns the trilinear form \eqref{rd-eq:T} below, a sum over geodesic triangles with vertices $1,a,x=ab$ and short side $|a|\le r$. Lemma~\ref{rd-lem:diagrams} sorts every such triangle by what sits at its centre: a common vertex (i), one central face (ii), or a V-configuration of three faces (iii). The rest of the proof works through these cases. Section~\ref{rd-sec:scales} fixes a scale $j$ determined by $r$ and calls a face \emph{small} or \emph{big} according as its component index is below $j$ or not. This refines the three cases into four classes:
\begin{center}\small
\begin{tabular}{lll}
\toprule
Lemma~\ref{rd-lem:diagrams} & class & estimated by\\
\midrule
(i) common vertex & (A) thin & Lemma~\ref{rd-lem:CR}, centre sets $\mathcal N_\Delta$\\
(ii) small central face & (A) thin & Lemma~\ref{rd-lem:CR}, centre sets $\mathcal N_\Delta$\\
(ii) big, long arcs on $\beta,\gamma_x$ & (B) two long arcs & Lemma~\ref{rd-lem:copy}, Theorem~\ref{rd-thm:compression}\\
(ii) big, short arc on $\beta$ or $\gamma_x$ & (C) one short arc & Theorem~\ref{rd-thm:compression}, current copy\\
(iii) faces meeting $\gamma_a$ small & (A) thin & Lemma~\ref{rd-lem:CR}, centre sets $\mathcal N_\Delta$\\
(iii) a big face meeting $\gamma_a$ & (D) vertex centre & Lemma~\ref{rd-lem:CR}, with overlap trees\\
\bottomrule
\end{tabular}
\end{center}
Thin triangles have a centre within $\Delta$ of all three sides, and Lemma~\ref{rd-lem:CR} applies directly. When the central face belongs to a big copy, we estimate the triangle through that copy: Lemma~\ref{rd-lem:copy} for class (B), and a direct use of the compression estimate for class (C), where Lemma~\ref{rd-lem:classes} shows the copy is current. In class (D) the centre sets must also contain the overlap trees $T(D,E)$ of adjacent big copies. The cardinality bound (H3) enters there, and in the recovery of corner points in Section~\ref{rd-sec:corners}. Section~\ref{rd-sec:corners} provides the recovery statements used in classes (B) and (C), Section~\ref{rd-sec:estimates} carries out the four estimates, and Section~\ref{rd-sec:assembly} adds them up.

Fix, for each $a\in H$, one geodesic $\gamma_a$ from $1$ to $a$. The family of paths from $x$ to $y$ is then chosen by translating $\gamma_{x^{-1}y}$. All functions in trilinear estimates are nonnegative and finitely supported. Set
\begin{equation}\label{rd-eq:T}
 T(f,g,h)=\sum_{ab=x}f(a)g(b)h(x).
\end{equation}
Assume $\supp f\subseteq B(1,r)$. Subscripts on $T$ indicate restrictions to a family of triples.

\begin{lemma}[Restricted centre estimate]\label{rd-lem:CR}
Let $C(x,y)$ be symmetric, left-equivariant subsets of $H$ containing $x,y$. Suppose $\diam C(x,y)\le R$ whenever $d(x,y)\le r$, and $|C(x,y)\cap B(x,R)|\le\varsigma$ for all $x,y$. If a family $\mathcal F$ consists of triples whose triangles satisfy
\[
 C(1,a)\cap C(a,x)\cap C(x,1)\ne\varnothing,
\]
then
\begin{equation}\label{rd-eq:CR}
 T_{\mathcal F}(f,g,h)\le\varsigma^{3/2}\norm f_2\norm g_2\norm h_2.
\end{equation}
\end{lemma}
\begin{proof}
This is the positive, restricted form of the argument of \cite[Proposition 1.7]{CR}; here is a direct proof. On indices $u,v\in B(1,R)$ and $w\in H$, define
\[
 F_{uv}=f(u^{-1}v)\ind_{1\in C(u,v)},\quad
 G_{vw}=g(v^{-1}w)\ind_{1\in C(v,w)},\quad
 J_{wu}=h(u^{-1}w)\ind_{1\in C(w,u)}.
\]
For a centre $t$ of $(1,a,x)$, take $u=t^{-1}$, $v=t^{-1}a$, $w=t^{-1}x$. The diameter assumption and the endpoints in $C(1,a)$ imply $u,v\in B(1,R)$. Thus each triple in $\mathcal F$ contributes at least once to $\Tr(FGJ)$. For a fixed value $a=u^{-1}v$, equivariance bounds the number of possible $u$ in $F$ by $|C(1,a)\cap B(1,R)|\le\varsigma$. The same argument applies to $G$ using its bounded index $v$, and to $J$ using $u$ and symmetry. Hence the three Hilbert--Schmidt norms are bounded by $\sqrt\varsigma$ times the respective function norms. The trace inequality gives \eqref{rd-eq:CR}. Finite truncations justify the calculation if an intermediate index set is infinite.
\end{proof}

We next consider triangles whose central face belongs to a big component copy $D$, as in Lemma~\partref{rd-lem:diagrams}{ii} with long exterior arcs, which gives class (B) above. We estimate their contributions by factoring through $D$ and applying the compression estimate. We proceed as follows.
For each triple choose a copy $D$ and a point $a'$ of $D$ on the short side $\gamma_a$, and write $a=a'\ell$ with tail $\ell$. Then $x=a'\ell b$, so the triple is determined by the position of $a'$ in $D$ together with $\ell$ and $b$.

Hypotheses (i) and (ii) of the next lemma say this bookkeeping can be inverted with bounded loss. From $x$ one recovers $D$ up to $U$ choices; from $b$ one recovers the component, the position of $a'$ in $D$, and the tail $\ell$, up to $W$ choices. Section~\ref{rd-sec:corners} supplies these recoveries, and Section~\ref{rd-sec:estimates} obtains $U\lesssim r+1$ and $W\lesssim(r+1)^2$ \eqref{rd-eq:long-UW}.

After Cauchy--Schwarz in $x$, which separates $h$, what remains is a Gram matrix indexed by positions and tails in a copy, summed over all translates $\omega Y^{(i)}$. The change of variables $b\mapsto(i,v,\ell,y_v\ell b)$, encoded by $\Psi$, costs at most $\sqrt W$. Collapsing the tails replaces $f$ by its geodesic prefix transform $\widehat f$ \eqref{rd-eq:fhat}, at a cost of $r+1$ in squared $\ell^2$ norm. Summing over translates turns the Gram matrix into the kernel $K_{\widehat\varphi}$ on a single copy, where $\widehat\varphi=\widetilde{\widehat f}*\widehat f$ is the autocorrelation of $\widehat f$. Theorem~\ref{rd-thm:compression} bounds this kernel uniformly in the girth.

Since $\norm{\widehat\varphi}_2$ involves $\norm{\lambda_H(\widehat f)}$, the resulting bound \eqref{rd-eq:bootstrap-input} contains $M_H(r)$ under a square root. Section~\ref{rd-sec:assembly} then closes the argument with the bootstrap inequality $M\le p\sqrt M+q$ \eqref{rd-eq:final-bootstrap}.

\begin{lemma}\label{rd-lem:copy}
Assign to each triple in a family $\mathcal F$ a component copy $D$ and a point $a'\in D\cap\gamma_a$. Put $\ell=(a')^{-1}a$. Suppose:
\begin{enumerate}[label=(\roman*)]
\item $x$ determines at most $U$ candidate copies containing the assigned $D$;
\item $b$ determines a set $W(b)$ of at most $W$ triples $(i,v,\ell)$ containing the assigned component, position of $a'$ in $D$, and tail.
\end{enumerate}
Define
\begin{equation}\label{rd-eq:fhat}
 \widehat f(y)^2=\sum_\ell f(y\ell)^2\ind_{y=\gamma_{y\ell}(|y\ell|-|\ell|)},
 \qquad \widehat\varphi=\widetilde{\widehat f}*\widehat f.
\end{equation}
Then
\begin{equation}\label{rd-eq:copy-bound}
 T_{\mathcal F}\le U\sqrt W\,
 \left\|\bigoplus_i K^{(i)}_{\widehat\varphi}\right\|^{1/2}\norm g_2\norm h_2,
 \qquad \norm{\widehat f}_2^2\le(r+1)\norm f_2^2.
\end{equation}
Here $\widetilde u(t)=u(t^{-1})$, and the direct sum includes all eligible components.
\end{lemma}
\begin{proof}
Pad the candidate lists to $U$ entries and sum over the candidate index. For one index let $A_D$ be the set of $x$ for which its candidate is $D$; these sets are disjoint. Write $I(y,\ell)$ for the indicator in \eqref{rd-eq:fhat}, and define
\[
 \Lambda_Dg(x)=\sum_{v,\ell}f(D_v\ell)I(D_v,\ell)
 \ind_{(i,v,\ell)\in W(\ell^{-1}D_v^{-1}x)}g(\ell^{-1}D_v^{-1}x).
\]
Positivity, followed by Cauchy--Schwarz and disjointness, bounds the contribution by
\begin{equation}\label{rd-eq:Lambda}
 \norm h_2\left(\sum_D\norm{\Lambda_Dg}_2^2\right)^{1/2}.
\end{equation}
For $D=\omega Y^{(i)}$, expand the square and sum over $\omega$. This gives the Gram matrix
\[
 \mathcal K^{(i)}((v,\ell),(v',\ell'))
 =\sum_\omega f(\omega y_v\ell)I(\omega y_v,\ell)
                    f(\omega y_{v'}\ell')I(\omega y_{v'},\ell').
\]
Let $\mathcal K=\bigoplus_i\mathcal K^{(i)}$, and define
\[
 \Psi\delta_b=\sum_{(i,v,\ell)\in W(b)}
        \delta_{((i,v,\ell),y_v^{(i)}\ell b)}.
\]
Expansion shows $\sum_D\Lambda_D^*\Lambda_D=\Psi^*(\mathcal K\otimes I)\Psi$: the equality of output coordinates is precisely $y_v\ell b=y_{v'}\ell'b'$. Also $\Psi^*\Psi$ is diagonal with entries $|W(b)|$, since the position, tail and product determine $b$. Thus the square sum in \eqref{rd-eq:Lambda} is at most $W\norm{\mathcal K}\norm g_2^2$.

For a vector $\xi(v,\ell)$ put $\eta(v)=\norm{\xi(v,\cdot)}_2$. Cauchy--Schwarz in $\ell$ gives
\begin{align*}
 \langle\mathcal K^{(i)}\xi,\xi\rangle
 &\le\sum_\omega\left(\sum_v\eta(v)\widehat f(\omega y_v)\right)^2\\
 &=\langle K^{(i)}_{\widehat\varphi}\eta,\eta\rangle,
\end{align*}
because $\sum_\omega\widehat f(\omega y_v)\widehat f(\omega y_{v'})
=\widehat\varphi(y_v^{-1}y_{v'})$. Hence $\norm{\mathcal K}\le\sup_i\norm{K^{(i)}_{\widehat\varphi}}$. Summing over the $U$ indices proves the first assertion. Finally, for each $a$ and each integer tail length between $0$ and $|a|$, the prescribed geodesic gives exactly one pair $(y,\ell)$ with $y\ell=a$ and $I(y,\ell)=1$. This proves the norm estimate. All operator identities can first be tested on finitely supported vectors; the resulting bounds justify their extensions.
\end{proof}

For any fixed $H,r$, the number
\begin{equation}\label{rd-eq:M}
 M_H(r)=\sup_{\supp f\subseteq B(1,r),\ f\ne0}
             \frac{\norm{\lambda_H(f)}}{\norm f_2}
\end{equation}
is finite, bounded by $|B(1,r)|^{1/2}$. It can equally be computed using nonnegative $f$. Since $\widehat f$ is supported in $B(1,r)$, Theorem~\ref{rd-thm:compression} gives
\begin{equation}\label{rd-eq:bootstrap-input}
 \left\|\bigoplus_i K^{(i)}_{\widehat\varphi}\right\|
 \le12(2r+1)^{3/2}M_H(r)(r+1)\norm f_2^2.
\end{equation}
Indeed, $\norm{\widehat\varphi}_2\le\norm{\lambda_H(\widehat f)}\norm{\widehat f}_2$. Here $M_H(r)$ is finite by its definition.

\subsection{Corner chains and recovery of endpoints}
\label{rd-sec:corners}

Consider a triangle with vertices $1,a,x=ab$, sides $\gamma_a$, $\beta=a\gamma_b$, $\gamma_x$, and $|a|\le r$. Suppose it has a central face $Z$ as in Lemma~\partref{rd-lem:diagrams}{ii}, with copy $D$.

We now establish the recovery statements used in classes (B) and (C): the central copy and the corner points must be determined, up to few choices, by one endpoint alone. Two facts do this. First, the corner chains next to the short side consist of faces of girth less than $3r$ with short rungs, so the long arcs of $Z$ start at depth at most $\zeta_r$; an endpoint and a starting depth then determine the copy. Second, at a corner the two sides leave $D$ either at the same point or both inside the first face of the corner chain lying in a different copy $E$ (Lemma~\ref{rd-lem:corner}). In the latter case both exit points lie in the overlap tree $T(D,E)$, and $E$ is determined by either exit and its outward direction. So one exit recovers the other within $\{p\}\cup T(D,E)$, which has at most $1+3\varkappa r$ points by (H3) (Lemma~\ref{rd-lem:corner-bounds}).

\begin{lemma}\label{rd-lem:corner}
Fix an original corner $z\in\{1,a,x\}$. Orient its two incident sides from $Z$ towards $z$. Let $p,q$ be their last points in $D$ in these directions. Either $p=q$, or there is a copy $E\ne D$ such that
\begin{equation}\label{rd-eq:corner}
 p,q\in T(D,E),
\end{equation}
and each outward side from its exit remains in $E$ for length greater than $(1/2-3\lambda)\sigma_E$. The copy $E$ is uniquely determined by either exit point and its outward geodesic segment.
\end{lemma}
\begin{proof}
Use the whole corner chain of Lemma~\ref{rd-lem:diagrams}, including external bigon discs. Its initial rung, or initial vertex, lies in $D$. Proceed towards $z$. On a common tree segment the two sides coincide; if their intersections with $D$ end there, then $p=q$. A face whose copy is $D$ has its whole boundary in $D$, so the sweep can continue past its outgoing rung or vertex.

Otherwise, let $F$ be the first face with copy $E\ne D$. Both incoming endpoints are in $D\cap E$. Each exterior arc of $F$ has length at least $(1/2-2\lambda)\sigma_E$, by Lemma~\ref{rd-lem:arc}. Its initial portion in $D$ has length less than $\lambda\min(\sigma_D,\sigma_E)$, by Lemma~\ref{rd-lem:intersection}. Therefore both sides leave $D$ within this same face, and their exits are $p,q$. Convexity prevents a later return. The remaining portions in $E$ have lengths greater than $(1/2-3\lambda)\sigma_E>\lambda\sigma_E$.

If the sweep reaches $z$ without an exit, then $p=q=z$. Finally, two different candidates $E,E'$ containing initial outward intervals of lengths at least $\lambda\sigma_E,\lambda\sigma_{E'}$ would share too long a piece. This proves uniqueness, including when the component girths differ.
\end{proof}

\begin{figure}[ht]
\centering
\begin{tikzpicture}[scale=1.3,font=\footnotesize]
 \coordinate (a) at (0,3); \coordinate (p) at (0,1.75); \coordinate (q) at (1.05,2.5275);
 \fill[gray!20] (0,0.3)--(0,1.4)--(1.5,2.325)--(2.7,1.785)--(2.7,0.3)--cycle;
 \fill[gray!45] (0,1.4)--(0,2.5)--(0.5,2.775)--(1.5,2.325)--cycle;
 \draw[thick] (0,0.3)--(a)--(2.7,1.785);
 \draw (0,2.5)--(0.5,2.775); \draw (0,2.75)--(0.25,2.8875);
 \draw[line width=2.2pt] (p)--(0,1.4)--(1.5,2.325)--(q);
 \fill (p) circle (1.6pt); \fill (q) circle (1.6pt);
 \node[left] at (p) {$p$}; \node[above right] at (q) {$q$}; \node[above] at (a) {$a$};
 \node at (1.4,1.2) {$Z$, copy $D$};
 \node at (0.55,2.2) {$F$};
 \node[right] at (1.3,2.85) {$F$: first face of copy $E\ne D$};
 \node[right] at (0.8,1.75) {$T(D,E)$};
 \draw[->] (0.82,1.77)--(0.55,1.75);
 \node[left] at (0,0.8) {$\gamma_a$}; \node[above right] at (2.3,1.97) {$\beta$};
\end{tikzpicture}
\caption{Lemma~\ref{rd-lem:corner} at the corner $a$. The central face $Z$ lies in the copy
$D$; $F$ is the first face of the corner chain whose copy $E$ differs from $D$. Its
incoming rung lies in $D\cap E$, and both sides leave $D$ inside $F$, at $p\in\gamma_a$ and
$q\in\beta$. The overlap tree $T(D,E)$ (thick) contains $p$ and $q$, and has fewer than
$\varkappa\min(\sigma_D,\sigma_E)$ vertices by (H3). Beyond $p$ and $q$ both sides stay in $E$ for
length greater than $(1/2-3\lambda)\sigma_E$, so either exit with its outward geodesic
determines $E$; the other exit then lies in $\{p\}\cup T(D,E)$.}
\label{rd-fig:corner}
\end{figure}
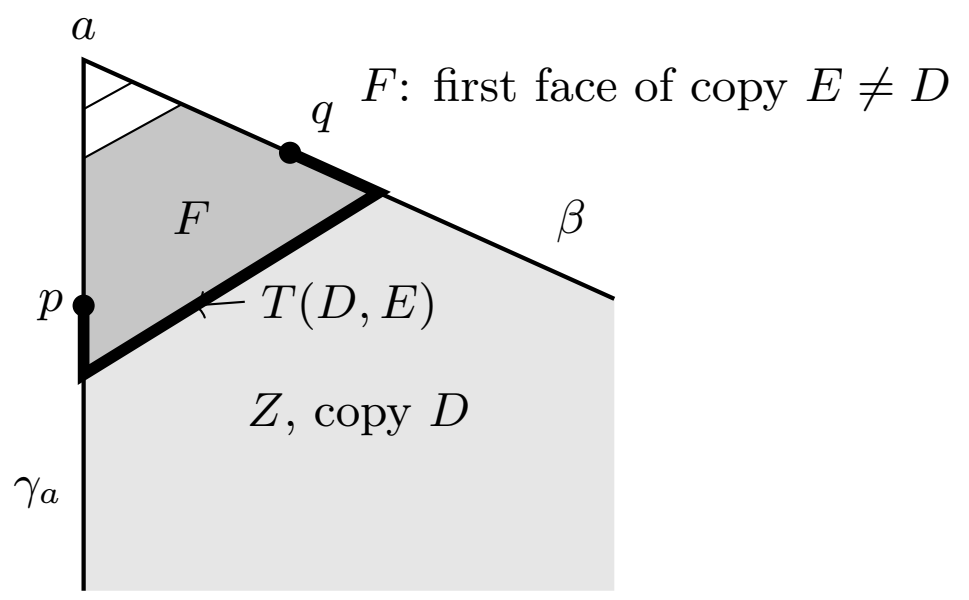

Given an exit $p$ and its outward geodesic $\eta$, denote by $E(p,\eta)$ the unique copy, if any, containing an initial interval of length at least its girth times $\lambda$. Discard a candidate equal to $D$. Then the other exit belongs to
\begin{equation}\label{rd-eq:recovery}
 \{p\}\cup T(D,E(p,\eta)),
\end{equation}
with an empty tree term if no candidate exists. This recovery set is determined by the stated geometric data.

\begin{lemma}\label{rd-lem:corner-bounds}
For $r\ge1$, put $\zeta_r=2r+2$ and $\vartheta_r=1+3\varkappa r$.
\begin{enumerate}[label=(\roman*)]
\item The initial points of $Z\cap\gamma_x$ and $Z\cap\beta$, measured from $1$ and $a$, have depths at most $\zeta_r$.
\item Let $a'$ be the last point of $D\cap\gamma_a$ towards $a$, and $P_a$ the first point of $D\cap\beta$ from $a$. Given either point and its outward geodesic towards $a$, the other has at most $\vartheta_r$ possibilities in \eqref{rd-eq:recovery}.
\item If $\sigma_D\le4r$, the two last points towards $x$ of $D\cap\gamma_x$ and $D\cap\beta$ similarly recover each other with at most $1+4\varkappa r$ possibilities.
\end{enumerate}
\end{lemma}
\begin{proof}
Every corner-chain face incident to the short side has an exterior arc of length at most $r$. Lemma~\ref{rd-lem:arc} gives
\begin{equation}\label{rd-eq:corner-girth}
 \sigma_E\le\frac{r}{1/2-2\lambda}\le\frac52r<3r.
\end{equation}
Inside the central disc, the rung connecting a point of $Z\cap\gamma_a$ to the initial point of $Z\cap\gamma_x$ has length less than $\lambda\sigma_E\le r/8$, if the corner chain is nonempty. If empty, the two points coincide. The same applies at $a$. Triangle inequality gives depths at most $r+r/8<\zeta_r$. External discs do not change this bound, since the short-side points are on the original geodesic $\gamma_a$.

For (ii), apply Lemma~\ref{rd-lem:corner}. The first different face lies in the whole $a$-corner chain and so satisfies \eqref{rd-eq:corner-girth}. Hypothesis (H3) bounds its overlap with $D$ by $3\varkappa r$, even if $D$ has much larger girth. When defining the recovery set, retain only candidates with girth at most $3r$. For (iii), the other copy may have arbitrarily large girth, but the smaller-girth bound in \eqref{rd-eq:intersection} is at most $4\varkappa r$.
\end{proof}

\subsection{Scale selection and triangle classes}
\label{rd-sec:scales}

For a given radius $r\ge1$, the index $j$ chosen in \eqref{rd-eq:scale} separates three ranges of components. Components below $j$ are \emph{small}: their diameters are at most $\Delta=A\sigma_{j-1}\le\log(r)/\log(2|S|)$, and balls of radius $\Delta$ contain at most $2r$ elements by \eqref{rd-eq:small-volume}. Small faces can therefore be absorbed into $\Delta$-neighbourhoods of the sides at polynomial cost, and triangles whose centre is small are thin. The component $\Theta_j$, when present, is the current component. Components above $j$ have girth greater than $100r$, so they meet the short side $\gamma_a$ only in arcs of length at most $r$, less than a hundredth of their girth. The cases of Lemma~\ref{rd-lem:diagrams} then split into the four classes below, according to the size of the central faces and where their long arcs lie. In classes (C) and (D), Lemma~\ref{rd-lem:classes} places the big central faces meeting the short side in the current component and gives the further bounds $\sigma_j\le4r$ in (C) and $\sigma_j<3r$ in (D). These bounds make the estimates of the next section polynomial in $r$.

Fix $r\ge1$ and choose the unique $j\ge1$ with
\begin{equation}\label{rd-eq:scale}
 e^{c\sigma_{j-1}}\le r<e^{c\sigma_j},\qquad \sigma_0=0.
\end{equation}
Call components with index less than $j$ small, and all the others big. Put $\Delta=A\sigma_{j-1}$. The elementary word-counting bound gives, in every stage and in the limit,
\begin{equation}\label{rd-eq:small-volume}
 |B(1,\Delta)|\le2(2|S|)^\Delta\le2r,
 \qquad \Delta\le\log(r)/\log(2|S|)\le r.
\end{equation}
All components with index greater than $j$ have girth greater than $100r$, by (H4). Components absent from a finite stage are simply omitted. Thus \eqref{rd-eq:scale} is still useful if $j>m$.

Use Lemma~\ref{rd-lem:diagrams} to assign every triple to one of the following classes, resolving any choices by a fixed rule.
\begin{enumerate}[label=(\Alph*)]
\item There is a vertex at distance at most $\Delta$ from all three sides. This includes a common vertex; a small central face; and a central V-configuration whose two faces meeting the short side are both small.
\item A big central face $Z$ has exterior arcs on $\beta$ and $\gamma_x$ each of length at least $\lambda\sigma_Z$.
\item A big central face has an exterior arc on $\beta$ or $\gamma_x$ of length less than $\lambda\sigma_Z$.
\item A V-configuration has at least one big central face meeting the short side.
\end{enumerate}

\begin{lemma}\label{rd-lem:classes}
In class (C), $Z$ belongs to $\Theta_j$, $\sigma_j\le4r$, and its other two exterior arcs have length at least $\sigma_j/4$. In class (D), any big central face meeting the short side belongs to $\Theta_j$ and $\sigma_j<3r$. Every big central V-face has both exterior arcs of length at least $\lambda$ times its own girth.
\end{lemma}
\begin{proof}
For class (C), write $p=|\partial Z|$. There are at most four interior arcs. After removing them and the exterior arc of length less than $\lambda\sigma_Z$, the other two arcs have total length greater than $p-5\lambda\sigma_Z$. Each is at most $p/2$, so each has length at least $(1/2-5\lambda)p\ge p/4$. One is on the short side, giving $\sigma_Z\le p\le4r$. Since later components have girth greater than $100r$, the big component must be current, namely $j$.

For a central V-face, Lemma~\ref{rd-lem:arc} with $k=3$ bounds both exterior arcs below by $(1/2-3\lambda)p\ge7p/20$. If one is on the short side, its girth is at most $20r/7<3r$, so a big such face is again current. The same lower bound exceeds $\lambda\sigma$ for every central V-face.
\end{proof}

This classification works at all radii.
Figure~\ref{rd-fig:classes} shows the data used in classes (B) and (C); class (A) is
Figure~\partref{rd-fig:shapes}{i}, or (ii) and (iii) with small central faces, and class (D)
is Figure~\partref{rd-fig:shapes}{iii} with $Z_{12}$ or $Z_{31}$ big.

\begin{figure}[ht]
\centering
\begin{tikzpicture}[scale=1.05,font=\footnotesize]
\begin{scope}
 \coordinate (o) at (0,0); \coordinate (a) at (0,1.6); \coordinate (x) at (5,0.8);
 \fill[gray!20] (0,0.35)--(0,1.25)--(0.9,1.456)--(3.0,1.12)--(3.0,0.48)--(0.9,0.144)--cycle;
 \draw[thick] (o)--(a)--(x)--cycle;
 \draw (0,1.25)--(0.9,1.456); \draw (0,0.35)--(0.9,0.144); \draw (3.0,1.12)--(3.0,0.48);
 \draw (0,0.18)--(0.45,0.072); \draw (0,1.43)--(0.45,1.528);
 \foreach \X in {3.6,4.2,4.6} \draw (\X,{1.6-0.16*\X})--(\X,{0.16*\X});
 \fill (0,1.25) circle (1.4pt); \fill (0.9,1.456) circle (1.4pt); \fill (0.9,0.144) circle (1.4pt);
 \node[left] at (0,1.25) {$a'$}; \node[above] at (0.9,1.456) {$P$};
 \node[below] at (0.9,0.144) {$\le \zeta_r$};
 \node at (1.8,0.8) {$Z$, copy $D$};
 \node[below left] at (o) {$1$}; \node[above left] at (a) {$a$}; \node[right] at (x) {$x$};
 \node[left] at (0,0.75) {$\gamma_a$}; \node[above] at (2.3,1.25) {$\beta$}; \node[below] at (2.3,0.35) {$\gamma_x$};
 \node at (2.5,-0.85) {(B) two long arcs};
\end{scope}
\begin{scope}[xshift=7.3cm]
 \coordinate (c) at (0,1);
 \coordinate (one) at ($(c)+(270:1.3)$); \coordinate (a) at ($(c)+(136:1.35)$);
 \coordinate (x) at (4.8,1.9);
 \coordinate (P1) at ($(c)+(126:1)$); \coordinate (P2) at ($(c)+(90:1)$);
 \coordinate (Q) at ($(c)+(55:1)$); \coordinate (A1) at ($(c)+(146:1)$);
 \fill[gray!20] (c) circle (1); \draw (c) circle (1);
 \draw[thick] (one)--($(c)+(255:1)$) arc[start angle=255,end angle=146,radius=1]--(a);
 \draw[thick] (a)--(P1) arc[start angle=126,end angle=90,radius=1];
 \draw[thick] (one)--($(c)+(285:1)$) arc[start angle=-75,end angle=55,radius=1];
 \draw[thick] (P2) .. controls (1.5,2.45) and (3.3,2.2) .. (x)
   coordinate[pos=0.3] (b1) coordinate[pos=0.55] (b2) coordinate[pos=0.8] (b3);
 \draw[thick] (Q) .. controls (1.5,1.4) and (3.3,1.6) .. (x)
   coordinate[pos=0.3] (k1) coordinate[pos=0.55] (k2) coordinate[pos=0.8] (k3);
 \foreach \i in {1,2,3} \draw (b\i)--(k\i);
 \foreach \P in {P1,P2,A1} \fill (\P) circle (1.4pt);
 \node[left] at (A1) {$a'$}; \node[above] at (P1) {$P_1$}; \node[above] at (P2) {$P_2$};
 \node[below] at (one) {$1$}; \node[left] at (a) {$a$}; \node[right] at (x) {$x$};
 \node at (c) {$Z$, copy $D$};
 \node[below] at ($(c)+(108:1)$) {$e$}; \node[above] at (b1) {$k$};
 \node[left] at ($(c)+(180:1)$) {$\gamma_a$}; \node[right] at ($(c)+(0:1)$) {$\gamma_x$};
 \node[above] at (b3) {$\beta$};
 \node at (2.0,-0.85) {(C) one short arc, on $\beta$};
\end{scope}
\end{tikzpicture}
\caption{Classes (B) and (C) of Section~\ref{rd-sec:scales}. (B) The central copy $D$ has
long arcs on $\beta$ and $\gamma_x$. The first point of $D\cap\gamma_x$ has depth at most
$\zeta_r$, so $x$ determines $D$; likewise $b$ determines $a^{-1}D$. The corner points
$a'\in\gamma_a$ and $P\in\beta$ recover each other with at most $\vartheta_r$ choices
(Figure~\ref{rd-fig:corner}). (C) The face $Z$ has long arcs on the short side $\gamma_a$ and
on $\gamma_x$, so $\sigma_j\le4r$, and a short arc on $\beta$. The coordinates of
\eqref{rd-eq:coordinates} are $b=\pi ek$, with $\pi$ from $a$ to $P_1$ ($|\pi|\le \zeta_r$),
$e=z_{v_1v_2}$ from $P_1$ to $P_2$ inside $D$ ($|e|\le A\sigma_j$), and $k$ from $P_2$ to
$x$. The estimate pairs $\widehat F$, built from $f$ at the $a$-corner, with $J$, built from
$h$ at the $x$-corner, through the compression $K_{\widehat G_k}$ on the copy $D$ of the
middle piece of $g$.}
\label{rd-fig:classes}
\end{figure}

\subsection{Estimates for the four classes}
\label{rd-sec:estimates}

Constants implicit in $\lesssim$ in this section depend only on $A,\varkappa,\lambda,|S|$ and are uniform in the stage, $j$, and $r$. Norms without subscripts on functions are $\ell^2$ norms.

\subsubsection{Thin triangles}
For arbitrary $x,y\in H$, set
\[
 C_0(x,y)=\mathcal N_\Delta(x\gamma_{x^{-1}y})\cup \mathcal N_\Delta(y\gamma_{y^{-1}x}).
\]
These sets are symmetric, equivariant, contain their endpoints, and have diameter at most $d(x,y)+2\Delta$. At radius $R=r+2\Delta$, their cardinalities near an endpoint are at most
\begin{equation}\label{rd-eq:thin-N}
 \varsigma_0=2(R+\Delta+1)|B(1,\Delta)|\lesssim(r+1)^2.
\end{equation}
Indeed, a point in the neighbourhood and in $B(x,R)$ lies within $\Delta$ of a geodesic vertex at depth at most $R+\Delta$. The centre in class (A) belongs to all three sets. Lemma~\ref{rd-lem:CR} gives
\begin{equation}\label{rd-eq:thin-bound}
 T_A\lesssim(r+1)^3\norm f\norm g\norm h.
\end{equation}

\subsubsection{Two long arcs}
Assign the central copy $D$ and $a'=\max_a(D\cap\gamma_a)$, meaning its last point towards $a$. Lemma~\ref{rd-lem:corner-bounds} and the long arc on $\gamma_x$ show that $x$ determines $D$ from an integer starting depth at most $\zeta_r$. At each depth there is at most one candidate, even across different component indices, by Lemma~\ref{rd-lem:intersection}. Thus $U\le \zeta_r+1$.

Translate by $a^{-1}$ and put $C=a^{-1}D$, $u=a^{-1}a'$. The long arc on $\gamma_b$ starts at depth at most $\zeta_r$, so $b$ gives at most $\zeta_r+1$ candidates for $C$. Once $C$ is specified, let $P$ be its first point on $\gamma_b$. The backward segment towards $1$ determines, by Lemma~\partref{rd-lem:corner-bounds}{ii}, a recovery set of at most $\vartheta_r$ points containing $u$. Include
\[
 (i(C),\pos_C(u),u^{-1})
\]
in $W(b)$ for all these candidates. The assigned tail is $u^{-1}=(a')^{-1}a$, so the data are correct. Therefore
\begin{equation}\label{rd-eq:long-UW}
 U\lesssim r+1,\qquad W\le(\zeta_r+1)\vartheta_r\lesssim(r+1)^2.
\end{equation}
Lemma~\ref{rd-lem:copy} and \eqref{rd-eq:bootstrap-input} yield
\begin{equation}\label{rd-eq:long-bound}
 T_B\lesssim(r+1)^{13/4}\sqrt{M_H(r)}\norm f\norm g\norm h.
\end{equation}
The supremum norm controls the direct sum of compression operators uniformly over the components and their translates.

\subsubsection{One short arc}
First suppose the short exterior arc is on $\beta$. Lemma~\ref{rd-lem:classes} gives a current copy $D$, with $\sigma_j\le4r$, and long arcs on $\gamma_a,\gamma_x$. Put $\Delta'=A\sigma_j\le4Ar$. For each current copy let $f_D(a)$ equal $f(a)$ if $\gamma_a$ contains an interval in $D$ of length at least $\lambda\sigma_j$, and zero otherwise. Let $h_D(x)$ equal $h(x)$ if $\gamma_x$ has such an interval starting at depth at most $\zeta_r$. Then
\begin{equation}\label{rd-eq:fh-partition}
 \sum_D\norm{f_D}^2\le(r+1)\norm f^2,
 \qquad \sum_D\norm{h_D}^2\le(\zeta_r+1)\norm h^2.
\end{equation}
The bounds follow by recording the first point of a qualifying interval; each integer depth identifies at most one copy. Every assigned triple is included in its $f_D,h_D$ terms.

Fix $D$. Let $P_1=D_{v_1}$ and $P_2=D_{v_2}$ be the first and last points of $D\cap\beta$, from $a$. Write $\tau_1,\tau_2$ for their depths on $\gamma_b$, and set
\begin{equation}\label{rd-eq:coordinates}
 \pi=\gamma_b(\tau_1),\quad
 e=\gamma_b(\tau_1)^{-1}\gamma_b(\tau_2)=z_{v_1v_2},\quad
 k=\gamma_b(\tau_2)^{-1}b.
\end{equation}
Then $b=\pi e k$, $a=D_{v_1}\pi^{-1}$, $x=D_{v_2}k$, with $|\pi|\le \zeta_r$ and $|e|\le \Delta'$. The latter follows from convexity and the diameter of $D$.

From $(a,D)$, recover $P_1$ using the last point of $D\cap\gamma_a$ towards $a$. Denote the resulting set by $\mathcal P(a,D)$; it has at most $\vartheta_r$ points. From $(x,D)$ similarly recover $P_2$ at the $x$-corner; call the set $\mathcal Q(x,D)$, of size at most $\vartheta'_r=1+4\varkappa r$. These sets depend only on the displayed data, by Lemma~\ref{rd-lem:corner-bounds}.

Define
\begin{align*}
 F(v,\pi)&=f_D(D_v\pi^{-1})\ind_{D_v\in\mathcal P(D_v\pi^{-1},D)},
 &\widehat F(v)^2&=\sum_\pi F(v,\pi)^2,\\
 J(v,k)&=h_D(D_vk)\ind_{D_v\in\mathcal Q(D_vk,D)}.
\end{align*}
For each $a$, a chosen point $D_v$ determines $\pi$ uniquely, and similarly for $x,k$. Consequently
\begin{equation}\label{rd-eq:FJ}
 \norm{\widehat F}^2\le \vartheta_r\norm{f_D}^2,
 \qquad \sum_{v,k}J(v,k)^2\le \vartheta'_r\norm{h_D}^2.
\end{equation}
For each $e,k$ define
\[
 \widehat G_k(e)^2=\sum_\pi g(\pi e k)^2
 \ind_{\pi=\gamma_{\pi e k}(|\pi|),\,
 \pi e=\gamma_{\pi e k}(|\pi|+|e|),\,
 |\pi|\le \zeta_r,\,|e|\le \Delta'}.
\]
Every $b$ has at most $(\zeta_r+1)(\lfloor \Delta'\rfloor+1)$ pairs of such cut depths. Therefore
\begin{equation}\label{rd-eq:G-cuts}
 \sum_{e,k}\widehat G_k(e)^2
 \le(\zeta_r+1)(\Delta'+1)\norm g^2.
\end{equation}
For fixed $D$, Cauchy--Schwarz in $\pi$ bounds its contribution by
\[
 \sum_{v_2,k}J(v_2,k)\sum_{v_1}\widehat F(v_1)
                         \widehat G_k(z_{v_1v_2}).
\]
The inner sum is $K_{\widehat G_k}^{\mathsf T}\widehat F$. Transposition preserves the operator norm, and each $\widehat G_k$ is supported in $B(1,\Delta')$. Theorem~\ref{rd-thm:compression}, followed by Cauchy--Schwarz in $(v_2,k)$, gives the fixed-copy bound
\begin{equation}\label{rd-eq:fixed-W1}
 12(\Delta'+1)^{3/2}
 \sqrt{\vartheta_r \vartheta'_r(\zeta_r+1)(\Delta'+1)}\,
 \norm{f_D}\norm g\norm{h_D}.
\end{equation}
Sum over $D$ using \eqref{rd-eq:fh-partition}. The coefficient is $O((r+1)^{9/2})$.

If the short arc is on $\gamma_x$, use the substitution $(a,b,x)\mapsto(a^{-1},x,b)$ and $(f,g,h)\mapsto(f^\vee,h,g)$. Translate the old triangle by $a^{-1}$ and reverse the short side. Its new short side is the deterministic geodesic $a^{-1}\overline{\gamma_a}$ depending only on $a^{-1}$. This deterministic choice supplies the geodesic data used in the preceding calculation. Its two long sides are the required translates of $\gamma_b,\gamma_x$. Thus the same estimate applies, and
\begin{equation}\label{rd-eq:W1-bound}
 T_C\lesssim(r+1)^{9/2}\norm f\norm g\norm h.
\end{equation}

\subsubsection{Vertex centres}
If this class is nonempty, Lemma~\ref{rd-lem:classes} gives $\sigma_j<3r$. Two distinct big copies are called adjacent along a geodesic if both meet it in intervals of length at least $\lambda$ times their respective girths, the intervals intersect, and at least one copy is current. For arbitrary endpoints define the symmetric equivariant set
\[
 C_V(x,y)=\bigcup_{\gamma\in\{x\gamma_{x^{-1}y},\,y\gamma_{y^{-1}x}\}}
 \left(\mathcal N_\Delta(\gamma)\cup
 \bigcup_{D,E\text{ adjacent along }\gamma}T(D,E)\right).
\]
The central vertex $t$ of a class (D) triangle is in all three sets. Indeed, on any side the two incident central faces either include a small face, which puts $t$ within $\Delta$ of that side, or are both big. In the latter case their exterior arcs are long and meet at the endpoint of their common rung. At least one of these faces meets the short side and is current. Their copies are distinct by graphical reduction, and $t$ lies in their overlap tree.

Put $\epsilon=\max(\Delta,\lambda\sigma_j)$ and $R=r+2\epsilon$. Every included overlap tree has diameter less than $\lambda\sigma_j$ and meets the geodesic. Thus $C_V(x,y)$ has diameter at most $d(x,y)+2\epsilon$. If such a tree meets $B(x,R)$, a point in the intersection of its two geodesic intervals is at depth at most
\[
 R_*=R+\lambda\sigma_j.
\]
The first point of each interval is therefore at depth at most $R_*$. At any starting depth at most one long copy is possible across all component indices. There are at most $(\lfloor R_*\rfloor+1)^2$ candidate pairs. Each contributes fewer than $\varkappa\sigma_j$ vertices. Including the neighbourhood terms and the two orientations gives
\begin{equation}\label{rd-eq:V-N}
 |C_V(x,y)\cap B(x,R)|
 \le2\bigl((R+\Delta+1)|B(1,\Delta)|+(R_*+1)^2\varkappa\sigma_j\bigr)
 \lesssim(r+1)^3.
\end{equation}
Lemma~\ref{rd-lem:CR} now gives
\begin{equation}\label{rd-eq:V-bound}
 T_D\lesssim(r+1)^{9/2}\norm f\norm g\norm h.
\end{equation}

\subsection{Closing the RD estimate and passing to the limit}
\label{rd-sec:assembly}

Combining \eqref{rd-eq:thin-bound}, \eqref{rd-eq:long-bound}, \eqref{rd-eq:W1-bound}, and \eqref{rd-eq:V-bound} gives constants $c_1,c_2$ independent of $H=G_m$ such that
\begin{equation}\label{rd-eq:final-bootstrap}
 M_H(r)\le c_1(r+1)^{13/4}\sqrt{M_H(r)}+c_2(r+1)^{9/2}.
\end{equation}
The quantity on the left is already finite by \eqref{rd-eq:M}. For nonnegative $M,p,q$, the inequality $M\le p\sqrt M+q$ implies $M\le p^2+2q$, since $p\sqrt M\le(M+p^2)/2$. Thus \eqref{rd-eq:final-bootstrap} implies \eqref{rd-eq:main}. The radius-zero case is immediate. Complex functions follow by pointwise domination of convolution by convolution of absolute values. This proves the RD assertion of Theorem~\ref{rd-thm:main}, directly for finite and infinite $m$.

The case $m=\infty$ of Theorem~\ref{rd-thm:main} was proved above together with
the finite cases. It also follows from the finite cases alone, by the following
general fact about direct limits of marked groups. Here it applies with
$P(r)=C(r+1)^{13/2}$ and $m_0(r)=1$.

\begin{proposition}[Passage to the limit]\label{rd-prop:limit}
Let $G_1\to G_2\to\cdots$ be surjective homomorphisms of groups marked by the
same finite alphabet $S$, with direct limit $G_\infty$, and let
$P\colon\N\to[0,\infty)$ be any function. Suppose that for every integer
$r\ge0$ there is an index $m_0(r)$ such that
\begin{equation}\label{rd-eq:prefix-bound}
 \norm{\lambda_{G_m}(f)}\le P(r)\norm f_2
 \qquad\text{for all }m\ge m_0(r)\text{ and all }f\in\C G_m
 \text{ with }\supp f\subseteq B_{G_m}(1,r).
\end{equation}
Then, for every integer $r\ge0$,
\begin{equation}\label{rd-eq:limit-bound}
 \norm{\lambda_{G_\infty}(f)}\le P(r)\norm f_2
 \qquad\text{for all }f\in\C G_\infty
 \text{ with }\supp f\subseteq B_{G_\infty}(1,r).
\end{equation}
In particular, if all sufficiently late $G_m$ have ball-form \RD{} with the same
constant and exponent, then so does $G_\infty$. The threshold $m_0(r)$ may
depend on $r$.
\end{proposition}
\begin{proof}
Fix $r$, a function $f$ on $G_\infty$ with support $E\subseteq B_{G_\infty}(1,r)$,
and a finitely supported function $g$ on $G_\infty$ with support $E'$. Choose for
each $x\in E$ a geodesic word $w_x$ in $S^{\pm1}$ representing it, so
$|w_x|\le r$, and for each $y\in E'$ any word $w_y$ representing it. Let
$\mathcal W$ be the finite set consisting of the words $w_x$, the words $w_y$,
and the products $w_xw_y$.

Two words of $\mathcal W$ that represent the same element of $G_\infty$
represent the same element of $G_m$ for all large $m$, because a relation in a
direct limit holds at a finite stage. Conversely, two words that represent the
same element of some $G_m$ represent the same element of $G_\infty$, since
$G_m$ maps to $G_\infty$. As $\mathcal W$ is finite, there is $m\ge m_0(r)$ such
that two words of $\mathcal W$ are equal in $G_m$ if and only if they are equal
in $G_\infty$.

For this $m$, let $\tilde x,\tilde y\in G_m$ be the elements represented by
$w_x,w_y$, and put $\tilde f(\tilde x)=f(x)$, $\tilde g(\tilde y)=g(y)$, and
zero elsewhere. By the choice of $m$ the maps $x\mapsto\tilde x$ and
$y\mapsto\tilde y$ are injective, and $\tilde x\tilde y=\tilde x'\tilde y'$ if
and only if $xy=x'y'$. Hence
\[
 \norm{\tilde f}_2=\norm f_2,\qquad\norm{\tilde g}_2=\norm g_2,\qquad
 \norm{\tilde f*\tilde g}_2=\norm{f*g}_2,
\]
the last because the two convolutions have the same coefficients, indexed by
corresponding products. Moreover $\supp\tilde f\subseteq B_{G_m}(1,r)$, because
$|w_x|\le r$. So \eqref{rd-eq:prefix-bound} gives
$\norm{f*g}_2\le P(r)\norm f_2\norm g_2$. Finitely supported $g$ are dense in
$\ell^2(G_\infty)$, and \eqref{rd-eq:limit-bound} follows.
\end{proof}

For lacunary hyperbolicity, apply \cite[Proposition 5.6]{Gru15}: the diameter--girth bound and the multiplicative sparseness implied by (H4) give a lacunary hyperbolic group. More precisely, the finite prefixes have $2A\sigma_m$-slim triangles and the shortest nontrivial word killed by the map to the full presentation has length at least $\sigma_{m+1}$; see Lemma~\ref{lem:graphical-radii}. Lemma~\ref{lem:graphical-radii} converts these bounds into the four-point constants and isometry radii used in the assembly construction.

\subsection{Realising the overlap condition by finite refinement}
\label{rd-sec:labelling}

We prove Theorem~\ref{rd-thm:construction} directly. Osajda's theorem \cite[Theorem 2.7]{Osa20} supplies, after a subsequence, a reduced strong girth-relative $C'(\lambda)$ labelling over a finite alphabet $S_0$. Arrange strictly increasing integer girths $\sigma_i\ge2$. We retain this labelling and add one independent finite coordinate. Because refinement only removes matching labelled paths, the original small cancellation condition is preserved.

We construct the extra coordinate by the local lemma.

\begin{lemma}\label{rd-lem:treecount}
In a graph of maximum degree $D$, the number of embedded subtrees with $k$ vertices containing a specified vertex is at most $(4D)^{k-1}$. The same bound holds for subtrees containing a specified edge. A prescribed rooted $k$-vertex tree has at most $|V(\Theta)|D^{k-1}$ embeddings in $\Theta$, where its root image is unrestricted.
\end{lemma}
\begin{proof}
Root a subtree at the specified vertex, choose an order on its children, and encode it by a rooted plane tree and neighbour choices. There are at most $4^{k-1}$ rooted plane trees with $k$ vertices, and at most $D^{k-1}$ choices for their edges. This overcounts the embedded subtrees and proves the first bound. An edge-containing subtree contains either specified endpoint, so the first bound already applies with that endpoint fixed. For the last assertion choose the root image and then the image of each child among at most $D$ neighbours.
\end{proof}

Fix the edge orientations so that the old labels are positive. Give each such edge a fresh independent letter from a set of size $q$. Reversing an edge reverses this fresh letter as a formal symbol. A tree isomorphism reversing a fixed edge cannot preserve these letters.

\begin{lemma}\label{rd-lem:equalities}
If $\phi:T_1\to T_2$ is an isomorphism between embedded $k$-vertex trees in one graph, with no fixed vertex, then the probability that independent fresh edge labels are preserved is at most $q^{-(k-1)/2}$.
\end{lemma}
\begin{proof}
A fixed unoriented edge would have to be reversed, making the event impossible. Otherwise, regard each of the $k-1$ edge-identification requirements as an edge in a graph whose vertices are the underlying independent edge variables. The in-degree and out-degree are at most one, because $\phi$ is injective. Thus its components are chains or cycles, with no loops. A chain of $l$ constraints imposes $l$ independent equalities. A cycle of $l\ge2$ constraints either is inconsistent with orientations or imposes $l-1\ge l/2$ independent equalities. Summing gives at least $(k-1)/2$ independent equalities.
\end{proof}

Choose constants
\begin{equation}\label{rd-eq:labelling-constants}
 b=4D,\quad B=8D^{A+2},\quad
 q\ge\max\{16B^4e^{10},\ D^{2(A+1)}(2b)^2e^{5/2}\},
\end{equation}
with $q$ integral. We construct the fresh labels inductively on $\Theta_n$, keeping earlier fresh labels fixed. Put $k_i=\sigma_i$. There are two kinds of bad event in the finite product probability space of the edges of $\Theta_n$.
\begin{enumerate}[label=(\roman*)]
\item An internal event specifies two embedded $k_n$-vertex trees in $\Theta_n$ and a fixed-point-free isomorphism between them; their fresh labels agree.
\item For each $i<n$ and each embedded $k_i$-vertex tree $T\subseteq\Theta_n$, a cross event says that its fresh labelling agrees under some rooted isomorphism with a tree in $\Theta_i$. Choose the root of $T$ once and for all; its image in $\Theta_i$ is unrestricted.
\end{enumerate}
We include all these events; avoiding them prevents the corresponding trees from having matching product labels. If a required tree does not exist, its class of events is empty.

The number of vertices in $\Theta_i$ is at most $D^{A\sigma_i+1}$. The internal-event probability is at most $q^{-(k_n-1)/2}$. By Lemma~\ref{rd-lem:treecount}, the number of internal events involving a fixed fresh edge is bounded by
\begin{equation}\label{rd-eq:X}
 X_n\le2b^{k_n-1}|V(\Theta_n)|D^{k_n-1}\le B^{k_n}.
\end{equation}
To see the factor two, distinguish whether the fixed edge is in the domain or image tree, and in the second case use the inverse isomorphism. For cross events of class $i$, the probability is at most
\begin{equation}\label{rd-eq:crossprob}
 |V(\Theta_i)|D^{k_i-1}q^{-(k_i-1)}
 \le\left(\frac{D^{2(A+1)}}q\right)^{k_i-1}.
\end{equation}
At most $b^{k_i}$ events in this class involve a fixed fresh edge.

Use local-lemma weights
\begin{equation}\label{rd-eq:LLLweights}
 u=(4B^{k_n})^{-1}\quad\text{for internal events},
 \qquad v_i=(2b)^{-k_i}\quad\text{for class }i.
\end{equation}
All weights are at most $1/2$. As the $k_i$ are distinct positive integers,
\[
 \sum_{i<n}b^{k_i}v_i=\sum_{i<n}2^{-k_i}\le1.
\]
Two events are dependent only if they share a fresh edge. An internal event uses at most $2(k_n-1)$ fresh edges. Using $\log(1-t)\ge-2t$ for $0\le t\le1/2$, its local-lemma right-hand side is at least
\begin{align*}
 u\exp\left(-4(k_n-1)B^{k_n}u
            -4(k_n-1)\sum_{i<n}b^{k_i}v_i\right)
 &\ge(4B^{k_n})^{-1}e^{-5(k_n-1)}\\
 &\ge(4B^2e^5)^{-(k_n-1)}.
\end{align*}
This is at least its probability by \eqref{rd-eq:labelling-constants}. A class-$i$ event uses $k_i-1$ fresh edges, so its right-hand side is at least
\begin{align*}
 v_i\exp\left(-2(k_i-1)B^{k_n}u
              -2(k_i-1)\sum_{h<n}b^{k_h}v_h\right)
 &\ge(2b)^{-k_i}e^{-(5/2)(k_i-1)}\\
 &\ge((2b)^2e^{5/2})^{-(k_i-1)},
\end{align*}
which dominates \eqref{rd-eq:crossprob}. The finite asymmetric Lov\'asz local lemma, for instance \cite[Lemma 2.1]{Osa20}, therefore gives a choice avoiding all events. This completes the induction with the same finite $q$ at every stage.

\begin{proposition}\label{rd-prop:refinement}
The product of the old labelling and the fresh coordinate satisfies (H3) with $\varkappa=1$.
\end{proposition}
\begin{proof}
Reducedness and (H2) are inherited. By Lemma~\ref{rd-lem:intersection}, an off-diagonal fibre-product component is a tree with embedded projections. These geometric and projection assertions follow from reducedness, convexity, and (H2), as shown in the proof of that lemma. If an overlap involving indices $i\le n$ had at least $k_i=\sigma_i$ vertices, it would contain a subtree with exactly $k_i$ vertices. For $i<n$, its two embedded projections would give an avoided cross event at stage $n$. For $i=n$, the projections give a fixed-point-free isomorphism: a fixed vertex would put the entire connected component on the diagonal, by reducedness. This is an avoided internal event. Both contradictions prove the claim.
\end{proof}

\begin{proof}[Proof of Theorem~\ref{rd-thm:construction}]
The product labelling has alphabet $S=S_0\times\{1,\dots,q\}$. It is reduced
and satisfies (H2) because the old labelling does, and it satisfies (H3) with
$\varkappa=1$ by Proposition~\ref{rd-prop:refinement}. The girths were arranged to
increase strictly, and the degree and diameter bounds are those of the given
graphs; this is (H1), and (a) follows. Deleting components removes pieces and
fibre-product components and creates none, so (b) holds. For (c), the alphabet
is now fixed, so $c=A\log(2|S|)$ is a constant, and since the girths tend to
infinity a subsequence satisfying \eqref{rd-eq:separation} can be chosen
recursively. A subsequence of an expander family is an expander family, and
each component embeds isometrically in the Cayley graph
\cite[Lemma~2.4]{ACGH}.
\end{proof}

\end{document}